\documentclass[10pt, 3p, fleqn]{elsarticle}
\usepackage[T1]{fontenc}
\usepackage{amssymb}
\usepackage{amsmath}
\usepackage{esint}
\usepackage{siunitx}
\usepackage{nicematrix}
\usepackage{caption}
\usepackage{subcaption}
\usepackage{float}
\usepackage{makecell}
\usepackage{mathrsfs}
\usepackage[x11names]{xcolor}

\AtBeginDocument{%
  \def\thefnote{\myfnsymbol{fnote}}}
\makeatletter
\def\myfnsymbol#1{\expandafter\@myfnsymbol\csname c@#1\endcsname}
\def\@myfnsymbol#1{\ifcase#1\or\textcolor{IllinoisDarkOrange}{$\ast$}\else\@ctrerr\fi}
\def\fntext[#1]#2{\g@addto@macro\@fnotes{%
   \refstepcounter{fnote}\elsLabel{#1}%
   \def\thefootnote{\thefnote}
   \global\setcounter{footnote}{\c@fnote-1}%
   \footnotetext{#2}}}
\makeatother

\biboptions{sort&compress}
\journal{Journal of Computational Physics}

\definecolor{IllinoisDarkOrange}{RGB}{200,65,19}
\usepackage{hyperref}
\hypersetup{
    pdfnewwindow=true,      
    colorlinks=true,        
    pdftitle={PCIC: Cylindrical Volume Moment Calculation and Interface Reconstruction for Sub-Grid Modeling in Volume of Fluid Methods},
    pdfauthor={Andrew Cahaly, Valentin Wasquel, Zonghao Zou, Olivier Desjardins, Fabien Evrard},
}
\makeatletter
    \AtBeginDocument{\def\@linkcolor{IllinoisDarkOrange}}
    \AtBeginDocument{\def\@anchorcolor{IllinoisDarkOrange}}
    \AtBeginDocument{\def\@citecolor{IllinoisDarkOrange}}
    \AtBeginDocument{\def\@filecolor{IllinoisDarkOrange}}
    \AtBeginDocument{\def\@urlcolor{IllinoisDarkOrange}}
    \AtBeginDocument{\def\@menucolor{IllinoisDarkOrange}}
    \AtBeginDocument{\def\@pagecolor{IllinoisDarkOrange}}
\makeatother

\usepackage{bbm}
\makeatletter
\DeclareFontFamily{U}{bbm}{}
\DeclareFontShape{U}{bbm}{m}{n}
 {  <5> <6> <7> <8> <9> <10> <12> gen * bbm
    <-> bbm10}{} 
\makeatother
\usepackage{xfrac}
\usepackage{multirow}
\usepackage{tikz,tikz-3dplot,tkz-euclide}
\usetikzlibrary{calc}
\usetikzlibrary{arrows}
\usetikzlibrary{shapes.geometric, arrows, intersections, through}
\usetikzlibrary{decorations.text, decorations.shapes, backgrounds}
\usetikzlibrary{arrows.meta}
\usetikzlibrary{colorbrewer}
\usetikzlibrary{patterns}
\usepackage{pgfplots} 
\usepackage{pgfplotstable}
\usepgfplotslibrary{colorbrewer}
\usepgfplotslibrary{colormaps}
\usepgfplotslibrary{patchplots}
\usetikzlibrary{colorbrewer}
\pgfplotsset{compat=newest}
\usepackage{ifthen}
\usepackage{setspace}
\usetikzlibrary{cd, arrows, matrix, intersections, math, calc}
\tikzset
  {midarrow/.style={decoration={markings,mark=at position 0.5 with
     {\arrow[xshift=1pt,scale=0.6]{angle 60[length=2pt]}}},postaction={decorate}}
  }
  \tikzset
  {revmidarrow/.style={decoration={markings,mark=at position 0.5 with
     {\arrowreversed[xshift=-1pt,scale=0.6]{angle 60[length=2pt]}}},postaction={decorate}}
  }
\usepackage[export]{adjustbox}
\usepackage{enumitem}
\usetikzlibrary{tikzmark}
\usetikzlibrary{decorations.pathreplacing,decorations.markings}

\usepackage[colorinlistoftodos]{todonotes}

\tikzset{
  on each segment/.style={
    decorate,
    decoration={
      show path construction,
      moveto code={},
      lineto code={
        \path [#1]
        (\tikzinputsegmentfirst) -- (\tikzinputsegmentlast);
      },
      curveto code={
        \path [#1] (\tikzinputsegmentfirst)
        .. controls
        (\tikzinputsegmentsupporta) and (\tikzinputsegmentsupportb)
        ..
        (\tikzinputsegmentlast);
      },
      closepath code={
        \path [#1]
        (\tikzinputsegmentfirst) -- (\tikzinputsegmentlast);
      },
    },
  },
  mid arrow/.style={postaction={decorate,decoration={
        markings,
        mark=at position .5 with {\arrow[#1]{stealth}}
      }}},
}
\tikzset{
    set arrow inside/.code={\pgfqkeys{/tikz/arrow inside}{#1}},
    set arrow inside={end/.initial=>, opt/.initial=},
    /pgf/decoration/Mark/.style={
        mark/.expanded=at position #1 with
        {
            \noexpand\arrow[\pgfkeysvalueof{/tikz/arrow inside/opt}]{\pgfkeysvalueof{/tikz/arrow inside/end}}
        }
    },
    arrow inside/.style 2 args={
        set arrow inside={#1},
        postaction={
            decorate,decoration={
                markings,Mark/.list={#2}
            }
        }
    },
}

\makeatletter
\renewcommand*\env@matrix[1][\arraystretch]{%
  \edef\arraystretch{#1}%
  \hskip -\arraycolsep
  \let\@ifnextchar\new@ifnextchar
  \array{*\c@MaxMatrixCols c}}
\makeatother

\newcommand{\cylradius}{r}
\newcommand{\diag}{\mathrm{diag}}

\begin{document}

\overfullrule=0mm

\begin{frontmatter}

\title{PCIC: Cylindrical Volume Moment Calculation and Interface Reconstruction for Sub-Grid Modeling in Volume of Fluid Methods}

\author[Cornell]{Andrew Cahaly\fnref{equal}}
\author[Supaero,Illinois]{Valentin Wasquel\fnref{equal}}
\author[Cornell]{Zonghao Zou}
\author[Cornell]{Olivier Desjardins}
\author[Illinois]{Fabien Evrard}
\fntext[equal]{\,These authors contributed equally to this work.}

\address[Cornell]{Sibley School of Mechanical and Aerospace Engineering, Cornell University, Ithaca, NY 14853, United States}
\address[Supaero]{ISAE-SUPAERO, Universit\'e de Toulouse, 10 Avenue Marc P\'elegrin, 31055 Toulouse, France}
\address[Illinois]{Department of Aerospace Engineering, University of Illinois Urbana-Champaign, Urbana, IL~61801, United States}

\begin{abstract}
The accurate modeling of topology changes remains a significant challenge in geometric Volume of Fluid (VOF) simulations. When using traditional single-plane reconstruction (PLIC), fluid structures smaller than the mesh size cannot be resolved and spurious numerical breakup is triggered; this impacts important flow statistics such as drop size distributions. Recent advances have introduced paraboloid and two-plane reconstructions, which have improved high-curvature performance and enabled sub-grid film reconstructions, respectively. However, sub-grid ligament reconstructions have remained elusive. In this work, a novel cylindrical interface reconstruction strategy called PCIC is introduced for sub-grid ligament modeling. PCIC is facilitated by deriving the analytical volume moments of quadratic cylinders clipping polyhedra; this allows for exact mass conservation during volume moment transport. From the transported moments, a straight circular cylinder can be defined in the center cell of a $5\times5\times5$ stencil. First, a quadratic principal curve is fitted to the normalized first-order moments in the stencil (the liquid barycenters), from which the cylinder's orientation and origin are approximated. The cylinder radius is then chosen to conserve volume. On-the-fly ligament detection is achieved using connected-component labeling and moments of inertia criteria, which allows for simulations to automatically choose between PLIC and PCIC in each interface cell at runtime. PCIC is demonstrated in multiphase flow test cases, where it exhibits robust reconstruction of sub-grid ligaments. This allows for relatively low-resolution PCIC simulations to provide comparable results to traditional high-resolution simulations.\\

\noindent © 2026. This manuscript version is made available under the CC-BY-NC-ND 4.0 license.\\
  \noindent \url{http://creativecommons.org/licenses/by-nc-nd/4.0}

\end{abstract}


\end{frontmatter}

\makeatletter
\renewcommand*{\@alph}[1]{\ifcase#1\or a\or b\or
   c\or d\or e\or f\or g\or h\or i\or j\or k\or l\or m\or n\or o\or p\or q\or r\or s\or t\or u\or v\or w\or x\or y\or z\or aa\or ab\or ac\or ad\or ae\or af\or ag\or ah\or ai\or aj\or
  ak\or al\or am\or an\or ao\or ap\or aq\or ar\or as\or at\or au\or
  av\or aw\or ax\or ay\or az%
  \or ba\or bb\or bc\or bd\or be\or bf\or bg\or bh\or bi\or bj\or
  bk\or bl\or bm\or bn\or bo\or bp\or bq\or br\or bs\or bt\or bu\or
  bv\or bw\or bx\or by\or bz%
  \or ca\or cb\or cc\or cd\or ce\or cf\or cg\or ch\or ci\or cj\or
  ck\or cl\or cm\or cn\or co\or cp\or cq\or cr\or cs\or ct\or cu\or
  cv\or cw\or cx\or cy\or cz%
  \or da\or db\or dc\or dd\or de\or df\or dg\or dh\or di\or dj\or
  dk\or dl\or dm\or dn\or do\or dp\or dq\or dr\or ds\or dt\or du\or
  dv\or dw\or dx\or dy\or dz%
  \or ea\or eb\or ec\or ed\or ee\or ef\or eg\or eh\or ei\or ej\or
  ek\or el\or em\or en\or eo\or ep\or eq\or er\or es\or et\or eu\or
  ev\or ew\or ex\or ey\or ez%
  \or fa\or fb\or fc\or fd\or fe\or ff\or fg\or fh\or fi\or fj\or
  fk\or fl\or fm\or fn\or fo\or fp\or fq\or fr\or fs\or ft\or fu\or
  fv\or fw\or fx\or fy\or fz%
  \or ga\or gb\or gc\or gd\or ge\or gf\or gg\or gh\or gi\or gj\or
  gk\or gl\or gm\or gn\or go\or gp\or gq\or gr\or gs\or gt\or gu\or
  gv\or gw\or gx\or gy\or gz%
  \else\@ctrerr\fi}
\makeatother
\renewcommand\thesubfigure{\alph{subfigure}}

\setlength\parindent{0pt}

\section{Introduction}
\label{Introduction}
The Volume of Fluid (VOF)~\cite{Hirt} method is a popular choice for capturing and transporting the immiscible fluid-fluid interface in multiphase flow simulations. By implicitly representing the interface with a phasic indicator function $\chi$, the VOF method can strictly enforce volume conservation and handle complex topology changes such as in liquid atomization. The indicator function $\chi$ is given a value of $0$ at locations in one fluid region (e.g., the gas) and a value of $1$ in the other region (e.g., the liquid). When used in a finite volume framework, the volume fraction field $\alpha$ is derived directly from $\chi$ for a given grid cell $\Omega_i$ with
\begin{equation} \label{eq:alpha}
\alpha_i=\frac{1}{V_i}\int_{\Omega_i}\chi\left(\mathbf{x}\right)\mathrm{d}V,
\end{equation}
where $V_i$ is the cell volume and $\mathbf{x}\in\mathbb{R}^3$ is the position. Thus, a cell with a volume fraction between $0$ and $1$ contains both fluid phases and the interface. The VOF method can directly advect the volume fraction field by integrating the material transport equation
\begin{equation} \label{eq:alpha_transport}
\frac{\partial \chi}{\partial t}+\mathbf{u}\cdot\nabla \chi=0,
\end{equation}
where $t$ is time and $\mathbf{u}$ is the fluid velocity, which is assumed to be continuous across the interface. A simple numerical integration of Eq.~\eqref{eq:alpha_transport} can lead to significant numerical diffusion errors, which effectively mix the phases and result in an artificial thickening of the interface. There exist several approaches for counteracting this behavior, such as algebraic VOF approaches~\cite{Hirt,Ubbink}, but the most effective solution is the geometric VOF method.

Geometric VOF approaches directly calculate the phasic volume fluxes into each grid cell; this requires knowledge of the sharp interface location for cutting the flux volumes into liquid and gas regions~\cite{Owkes}. Therefore, an approximate sharp interface must be reconstructed from the volume fraction field. The Piecewise Linear Interface Calculation (PLIC)~\cite{Youngs} has been the most widely used approach to produce such approximations. PLIC fits a plane in each mixed-phase grid cell such that the cell's volume fraction is conserved exactly. Some of the more popular PLIC algorithms include the LVIRA and ELVIRA methods of \citet{Pilliod}. These methods minimize a least-squares loss function between the actual volume fractions in a stencil of cells and the volume fractions formed by a plane cutting through the same stencil. Thus, these methods assume an approximately linear interface on the scale of the stencil. An alternative approach is the recently developed PLIC-Net method~\cite{Cahaly}, which uses a neural network to predict a plane's orientation in a cell given a stencil of volume fractions and their normalized first order moments (i.e., liquid and gas barycenters). Such barycenters are readily available in geometric VOF schemes since they can be directly transported alongside the volume fraction field and then updated to be consistent with the next interface reconstruction at negligible additional cost~\cite{Dyadechko,Chenadec,Owkes2}. Since paraboloid training data was used, PLIC-Net does not assume a linear interface in its stencil and outperforms ELVIRA and LVIRA in high-curvature accuracy, spurious interface errors, and computational cost.

Although PLIC-Net enables accurate single-plane reconstructions at scales approaching the mesh size, it still suffers from the fundamental flaw of all traditional geometric VOF methods: any interfacial features below the mesh size cannot be resolved, leading to spurious numerical breakup of the interface. Indeed, such ``automatic'' interface topology change, which at first appears to be a main advantage of VOF methods, in fact constitutes one of the largest remaining challenges for accurate VOF simulations. Although Adaptive Mesh Refinement (AMR) can supply increased resolution to regions with high curvatures and small features, this is not a solution to the overarching problem. Firstly, this only delays numerical breakup to a later point without addressing the underlying physics. Secondly, the cost can become prohibitive, since it is common for multiphase flows to have small structures with characteristic sizes of a few microns or below, which can be several orders of magnitude less than the characteristic sizes of the resolved structures~\cite{Villermaux}. For example, recent work by \citet{Ling} numerically studied the aerodynamic breakup of a millimeter-scale drop in the bag breakup regime. The researchers used an AMR approach with an equivalent uniform resolution of $2048$ cells across the initial drop diameter. This required a computational cost of 4.3 million core-hours, and the observed interface breakup was still due to a numerical cutoff length controlled by the grid size.

Recent research has started to address this challenge. The Reconstruction with 2 Planes (R2P) method~\cite{Chiodi, Han3} introduces a sub-grid representation for thin films by allowing two planes to co-exist in the same cell. The algorithm uses transported surface normal vectors to create initial guesses for plane orientations, and then minimizes a least-squares loss function on a stencil of cells that accounts for volume fractions, their liquid and gas barycenters, and surface area. Volume is conserved throughout the optimization process. This method has been demonstrated in several cases, such as drop breakup in turbulent cross-flows, where it has enabled accurate predictions of drop size distributions when coupled with physics-based breakup modeling~\cite{Han4}. Similarly, R2P was also successfully deployed in pressure swirl atomization cases, as shown in \citet{Giliberto}. Meanwhile, higher order methods such as the Piecewise Parabolic Interface Calculation (PPIC) have also been explored~\cite{Evrard}. By fitting a paraboloid in each interface cell, PPIC allows for improved reconstructions at higher curvatures with an increased order of accuracy~\cite{Evrard1}. Several algorithms have been suggested for the PPIC interface reconstruction, such as the approach of \citet{Jibben}, which uses a PLIC reconstruction to align a paraboloid reference frame and then minimizes the volume mismatch between the paraboloid and PLIC reconstruction in a stencil. A Moment of Fluid (MOF) method~\cite{Dyadechko} has also been proposed by \citet{Evrard2}, which directly minimizes a loss function with zeroth, first, and second-order volume moments from a single cell. Recent work by \citet{Lopez} expanded on the Jibben method by introducing an iterative solution, where the PLIC reconstruction is iteratively updated by the previous paraboloid reconstruction. Although PPIC methods are successful for reasonably resolved interfaces, they so far have failed to capture sub-grid elements such as fluid ligaments.

There have been several proposed approaches to address physical ligament dynamics and breakup in multiphase flow simulations. \citet{Chirco} proposed an algorithm using the second-order moments of a volume fraction indicator function to detect thin fluid structures. This algorithm was used to identify films and ligaments before they became thinner than the mesh size; thin films were then punctured at a chosen thickness threshold to induce Taylor--Culick retraction. However, ligament modeling and sub-grid representations were not directly addressed in the work. Meanwhile, \citet{Kim} argued that a theoretical model based on the Rayleigh--Plateau instability can be used to predict final droplet sizes from resolved ligaments without directly resolving sub-grid ligaments. By detecting the thin structures near the mesh size with an eccentricity metric and then converting them to Lagrangian droplets based on the theory, the researchers showed good agreement with drop size distributions produced by spray atomization experiments. However, they note that ligament breakup was still being triggered prematurely by the mesh size, and their approach was applied in the context of Newtonian flows with low liquid-gas viscosity ratios. Numerical studies by \citet{Farsoiya} demonstrate how the viscosity ratio can dramatically change breakup behavior for a drop in homogeneous isotropic turbulence. At high viscosity ratios, thin ligaments can form and be maintained well below typical mesh sizes prior to breakup, and the characteristic breakup times are impacted. The work by \citet{Farsoiya} relied on AMR simulations with an equivalent uniform mesh of up to $1024^3$. Non-Newtonian effects can further complicate the physics by introducing ligament stabilization and delaying breakup even longer, resulting in extremely thin and meandering ligaments~\cite{Datta,Chandra,Christanti}. Therefore, sub-grid representations of thin ligaments become imperative for the accuracy of these complex simulations. Work by \citet{Audoly} proposed a model for non-Newtonian ligaments by connecting 1D piecewise Lagrangian rods with internal representations of viscous stresses. Although this is a useful physics model, a full sub-grid interface reconstruction scheme for ligaments in VOF simulations remains to be demonstrated.

In this work, a novel cylindrical interface reconstruction is derived, implemented, and tested in complex geometric VOF simulations for sub-grid ligament modeling. This method is referred to as the Piecewise Cylindrical Interface Calculation, or PCIC, hereinafter. The first step to enable PCIC is to calculate the volume moments of quadratic cylinders clipping arbitrary polyhedra. Such volume moments must be accurately and robustly calculated at low cost during the VOF advection process. Traditional methods for such calculations include Monte-Carlo methods~\cite{Evans,Hahn}, recursive octree subdivision~\cite{Kudela,Divi}, or surface triangulation techniques~\cite{Tolle}. While versatile, these approaches tend to be computationally expensive when reaching the high accuracy required in VOF simulations. Instead, closed-form expressions are derived by using the divergence theorem to reduce the three-dimensional moment integrals into a collection of simpler one-dimensional integrals~\cite{Kromer,Kromer2,Antolin}. This approach builds on the work of \citet{Evrard}, which recently derived similar expressions for the volume moments of a paraboloid clipping a polyhedron. The resulting calculation of the volume moments achieves machine precision while being several orders of magnitude faster than the other existing approaches. With the ability to calculate the volume moments produced by cylinder-clipping, the reconstruction problem can then be addressed. The orientation and location of a cylinder in a given grid cell are determined by fitting a quadratic principal curve~\cite{Hastie} to the phasic barycenters in a $5\times5\times5$ stencil of cells. The principal curve is fitted iteratively by calculating control points from volume-weighted averages of the barycenters. The point on the resulting quadratic curve closest to the center of the stencil is then selected to serve as the cylinder origin, and the tangent vector at the same point is the cylinder orientation. The radius is set to conserve volume. This method is deployed in multiphase flow simulations by using on-the-fly ligament detection. Thin ligaments approaching the grid size are detected by using connected-component labeling (CCL)~\cite{He} with moments of inertia criteria in a $5\times5\times5$ stencil of cells. When a ligament is detected in a cell, PCIC is automatically selected as its reconstruction method. Detection and modeling of ligament tips is also carried out to create clean and realistic reconstructions. When tested in simulations, PCIC enables accurate thin ligament representations at the sub-grid scale on relatively coarse meshes and produces results comparable to well-resolved AMR approaches. 

Section~\ref{Methods1} describes the method used to derive the volume moments of a polyhedron clipped by a quadratic cylinder. Section~\ref{Methods2} discusses the methods for PCIC reconstruction, ligament detection, and deployment in VOF simulations. Section~\ref{Results1} presents a range of test cases for the volume moment calculations. Section~\ref{Results2} demonstrates the PCIC reconstruction in static, pure advection, and multiphase simulation tests and discusses its performance. Section~\ref{Conclusion} concludes this work.

\section{Method for Calculating Volume Moments of Cylinders Clipping Polyhedra}
\label{Methods1}
In this section, the techniques used to derive the volume moments of quadratic cylinders clipping polyhedra are introduced. Section~\ref{sec:fondation} summarizes important notation and results from previous work on paraboloid volume moments~\cite{Evrard}, Section~\ref{sec:new_work} details the methods used to derive the new analytical volume moment expressions for quadratic cylinders, and Section~\ref{sec:robustness} discusses approaches to robustly handle ill-posed edge cases.

\subsection{Mathematical Foundation}\label{sec:fondation}
For the sake of completeness, this section revisits key notations, definitions, and steps of the approach introduced by \citet{Evrard}; a more comprehensive description is available in \cite{Evrard} directly.

\begin{figure}[b!]\centering\vspace{-1.5em}
    \includegraphics{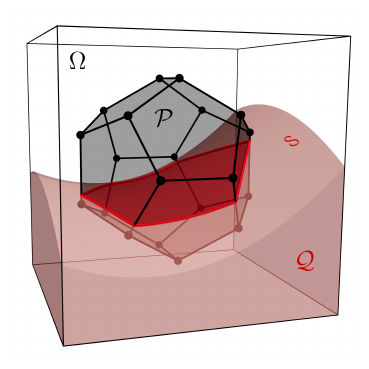}\includegraphics{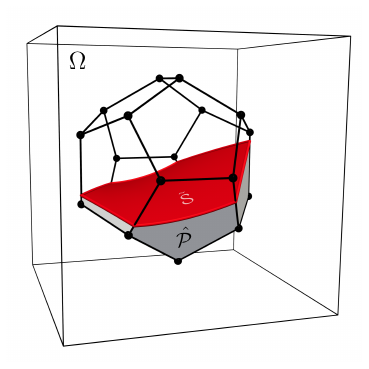}
    \vspace{-1.5em}
    \caption{Notation summary: $\Omega$ is a subset of $\mathbb{R}^3$, here illustrated as a cube. $\mathcal{P} \subset \Omega$ is a polyhedron, here illustrated as a regular dodecahedron, and  $\mathcal{Q} \subset \Omega$ is the clipping region bounded inside $\Omega$ by the quadratic surface $\mathcal{S}$, here illustrated as a hyperbolic paraboloid. Inside $\Omega$, the surface $\mathcal{S}$ can be parameterized by a single-valued function of the coordinates $x$ and $y$. The portion of $\mathcal{S}$ inside $\mathcal{P}$ is denoted as $\tilde{\mathcal{S}} \equiv \mathcal{S} \cap \mathcal{P}$. The intersection of $\mathcal{P}$ with $\mathcal{Q}$ is denoted as $\hat{\mathcal{P}} \equiv \mathcal{P} \cap \mathcal{Q}$.}\label{fig:intro}
\end{figure}

\begin{figure}
    \centering\vspace{-1em}
    \begin{tabular}{ccc}
    \subfloat[\centering
    Integration domain for $\boldsymbol{\mathcal{M}}^{\hat{\mathcal{P}}_1}_i$]{\adjincludegraphics[width=0.3\textwidth,trim=0 0 0 0,clip=true]{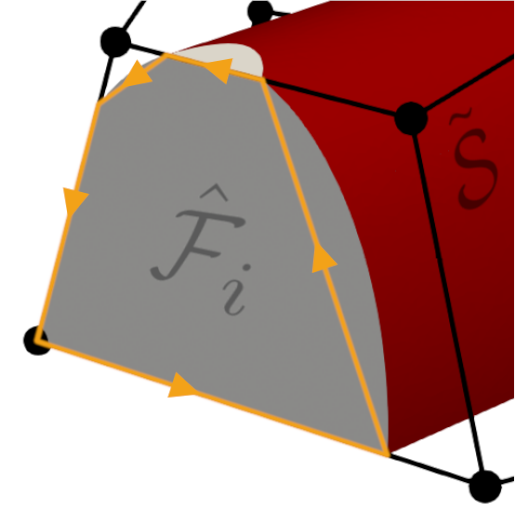}} \quad
    \subfloat[\centering 
    Integration domain for $\boldsymbol{\mathcal{M}}^{\hat{\mathcal{P}}_2}_i$]{\adjincludegraphics[width=0.3\textwidth,trim=0 0 0 0,clip=true]{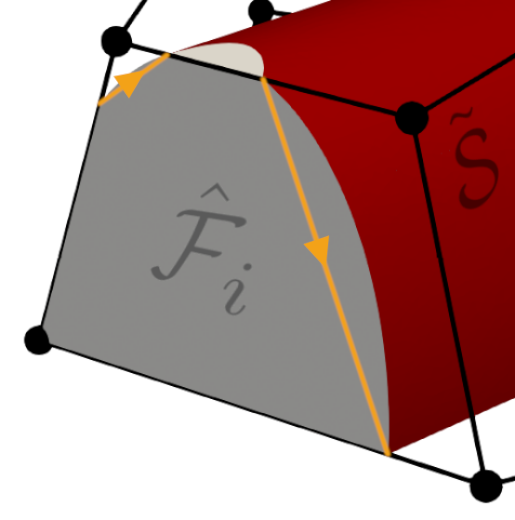}} \quad
    \subfloat[\centering
    Integration domain for $\boldsymbol{\mathcal{M}}^{\hat{\mathcal{P}}_3}_i$]{\adjincludegraphics[width=0.3\textwidth,trim=0 0 0 0,clip=true]{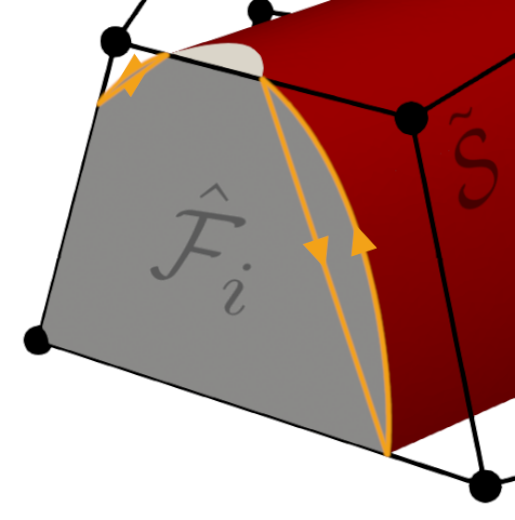}}
    \end{tabular}
    \caption{Integration domains for the contribution of face $\hat{\mathcal{F}}_i$ to $\boldsymbol{\mathcal{M}}^{\hat{\mathcal{P}}_1}$, $\boldsymbol{\mathcal{M}}^{\hat{\mathcal{P}}_2}$ and $\boldsymbol{\mathcal{M}}^{\hat{\mathcal{P}}_3}$}
    \label{fig:M123}
\end{figure}

Let the study be restricted to some subset $\Omega$ of $\mathbb{R}^3$. Consider the polyhedron $\smash{\mathcal{P} \subset \Omega}$, bounded by the planar polygonal faces $\smash{\mathcal{F}_i, i \in \{1,\ldots,n_\mathcal{F}\}}$, and the ``clipping region'' $\smash{\mathcal{Q}\subset \Omega}$. The region $\smash{\mathcal{Q}}$ is bounded by the quadratic surface $\smash{\mathcal{S}}$, which is required to be described inside $\Omega$ by the parametric representation\begin{equation}
    \mathbf{x}_\mathcal{S}(x,y) = 
   \begin{bmatrix}
    x &
    y &
    z_{\mathcal{S}}(x,y)
    \end{bmatrix}^\intercal \, , \label{def:paramSquadric}
\end{equation}
with $\smash{z_{\mathcal{S}}}$ a single-valued function.
The regions $\smash{\hat{\mathcal{P}} \equiv \mathcal{P} \cap \mathcal{Q}}$ and $\smash{\hat{\mathcal{F}}_i \equiv \mathcal{F}_i \cap \mathcal{Q}}$ correspond to the intersections of the polyhedron $\smash{\mathcal{P}}$ and face $\smash{\mathcal{F}_i}$ with the clipping region $\smash{\mathcal{Q}}$, respectively; $\smash{\tilde{\mathcal{S}} \equiv \mathcal{S} \cap \mathcal{P}}$ is the intersection of the surface $\smash{\mathcal{S}}$ with the polyhedron~$\smash{\mathcal{P}}$. These notations are illustrated in Fig.~\ref{fig:intro}.
The objective is to derive analytical expressions for the first geometric moments of the clipped polyhedron $\smash{\hat{\mathcal{P}}}$, defined as
\begin{equation}
    \boldsymbol{\mathcal{M}}^{\hat{\mathcal{P}}} = \int_{\hat{\mathcal{P}}}  \boldsymbol{\Upsilon}(\mathbf{x}) \ \mathrm{d}\mathbf{x} \, , \quad \text{ with } \boldsymbol{\Upsilon} : \left\{\begin{array}{lcl} 
        \mathbb{R}^3 & \!\to & \mathbb{R}^4 \\
        \mathbf{x} & \!\mapsto & \begin{bmatrix}  1 & x & y & z \end{bmatrix}^\intercal \end{array}\right.  \, ,
\label{eq:moments}
\end{equation}
and with $\smash{\mathbf{x} = \begin{bmatrix} x & y & z \end{bmatrix}^\intercal}$ the position vector in the Cartesian coordinate system equipped with the orthonormal basis $\smash{(\mathbf{e}_x,\mathbf{e}_y,\mathbf{e}_z)}$. 
Upon applying the divergence theorem, 
$\smash{\boldsymbol{\mathcal{M}}^{\hat{\mathcal{P}}}}$ can be expressed as a boundary integral over $\smash{\partial{} \hat{\mathcal{P}}}$,
\begin{equation}
	\boldsymbol{\mathcal{M}}^{\hat{\mathcal{P}}} 
    = \int_{\partial\hat{\mathcal{P}}} \boldsymbol{\Phi}(\mathbf{x}) (\mathbf{n}_{\partial\hat{\mathcal{P}}} \cdot \mathbf{e}_z) \, \mathrm{d} \mathbf{a}  \, , \quad \text{ with } \boldsymbol{\Phi} : \left\{\begin{array}{lcl} 
        \mathbb{R}^3 & \to & \mathbb{R}^4 \\
        \mathbf{x} & \mapsto & {\displaystyle \int_0^z} \boldsymbol{\Upsilon}(\mathbf{x}) \ \mathrm{d}z \end{array}\right. ,
\label{eq:vol_after_div1}
\end{equation}
where $\smash{\mathrm{d}\mathbf{a}}$ is an infinitesimal surface element on $\smash{\partial\hat{\mathcal{P}}}$ and $\smash{\mathbf{n}_{\partial\hat{\mathcal{P}}}}$ is the outward-pointing normal to $\smash{\partial\hat{\mathcal{P}}}$. 
This integral can then be split into the contribution of the clipped surface $\tilde{\mathcal{S}}$, and the combined contributions of each planar face $\hat{\mathcal{F}}_i$ of the clipped polyhedron.
Expressing the $z$-component of each face of $\smash{\partial\hat{\mathcal{P}}}$ as a function of $x$ and $y$, $\smash{\boldsymbol{\mathcal{M}}^{\hat{\mathcal{P}}}}$ is recast as an integral on $\tilde{\mathcal{S}}^\perp$ and $\smash{\cup_i \hat{\mathcal{F}}_i^\perp}$, the projections of the boundary faces onto the $xy$-plane, i.e.,
\begin{equation}
    \boldsymbol{\mathcal{M}}^{\hat{\mathcal{P}}} = \int_{\tilde{\mathcal{S}}^\perp}\boldsymbol{\Phi}_{\mathcal{S}}(x,y) \,\mathrm{d}x\, \mathrm{d}y + \sum_{i=1}^{n_{\mathcal{F}}}\text{sign}(\mathbf{n}_i \cdot \mathbf{e}_z)\int_{\hat{\mathcal{F}}_i^\perp}\boldsymbol{\Phi}_{\mathcal{F}_i}(x,y) \,\mathrm{d}x\, \mathrm{d}y \ , \label{eq:divperp}
\end{equation}
with $\smash{\boldsymbol{\Phi}_{\mathcal{S}}(x, y) \equiv \boldsymbol{\Phi}(x, y, z_{\mathcal{S}}(x, y))}$ and $\smash{\boldsymbol{\Phi}_{\mathcal{F}_i}(x, y) \equiv \boldsymbol{\Phi}(x, y, z_{\mathcal{F}_i}(x, y))}$, in which $\smash{z_{\mathcal{F}_i}}$ is the parameterization of the $z$-component of $\mathcal{F}_i$ as a function of $x$ and $y$, and with $\smash{\mathbf{n}_i}$ the outward-pointing normal to the planar face $\smash{\mathcal{F}_i}$. 
Upon applying the divergence theorem a second time, 
$\smash{\boldsymbol{\mathcal{M}}^{\hat{\mathcal{P}}}}$ can be formulated as a sum of one-dimensional integrals over each edge of the projected faces\footnote{That is because each edge of $\smash{\tilde{\mathcal{S}}^\perp}$ corresponds to exactly one edge of the projected faces $\smash{\cup_i\hat{\mathcal{F}}^\perp_i}$, hence the boundary of $\smash{\tilde{\mathcal{S}}^\perp}$ is itself a subset of the union of the clipped planar face boundaries.} 
$\hat{\mathcal{F}}^\perp_i$. The moments are then given as
\begin{align}
	\boldsymbol{\mathcal{M}}^{\hat{\mathcal{P}}} = \sum\limits_{i=1}^{n_{\!\mathcal{F_{\phantom{k}\!\!}}}} \sum\limits_{j = 1}^{n_{\partial\hat{\mathcal{F}}_i}}  \int_{0}^{1} \left(\boldsymbol{\Psi}_{\mathcal{F}_i}(x_{i,j}(t),y_{i,j}(t)) - \mathbbm{1}^{\partial\tilde{\mathcal{S}}}_{i,j} \, \boldsymbol{\Psi}_{\mathcal{S}}(x_{i,j}(t),y_{i,j}(t)) \right) y^\prime_{i,j}(t) \ \mathrm{d}t \, , \label{eq:vol_after_div7}
\end{align}
with 
\begin{equation}
    \boldsymbol{\Psi}_\square : \left\{\begin{array}{lcl} 
        \mathbb{R}^2 & \to & \mathbb{R}^4 \\
        (x,y) & \mapsto & {\displaystyle \int_0^x} \boldsymbol{\Phi}_\square(x,y) \ \mathrm{d}x\end{array}\right. ,
\end{equation}
and with $\smash{\mathbf{x}_{i,j}(t) = x_{i,j}(t)\mathbf{e}_x + y_{i,j}(t)\mathbf{e}_y +z_{i,j}(t) \mathbf{e}_z}$ some arbitrary parameterization of the $j$th edge of the $i$th clipped face $\smash{\hat{\mathcal{F}}_i}$ for the parameter $\smash{t \in [0,1]}$. The number of edges that delimit the clipped face $\smash{\hat{\mathcal{F}}_i}$ is denoted by $\smash{n_{\partial \hat{\mathcal{F}}_i}}$, and the coefficient $\smash{\mathbbm{1}^{\partial\tilde{\mathcal{S}}}_{i,j}}$ is defined~as\begin{equation}
    \mathbbm{1}^{\partial\tilde{\mathcal{S}}}_{i,j} = \left\{\begin{array}{cl} 
        1 & \text{if the $j$th edge of $\hat{\mathcal{F}}_i$ also belongs to $\partial\tilde{\mathcal{S}}$} \\ 0 & \text{otherwise}\end{array}\right. .
\end{equation}
To further simplify the expression of the first moments, the linear parameterization denoted $\smash{\bar{\mathbf{x}}_{i,j}(t) = \bar{x}_{i,j}(t) \mathbf{e}_x}\linebreak\smash{+\;\bar{y}_{i,j}(t) \mathbf{e}_y + \bar{z}_{i,j}(t) \mathbf{e}_z}$, linking the start-point of the $j$th edge of $\hat{\mathcal{F}}_i$ to its end-point by a straight line, is introduced. This parameterization is only exact for the edges of $\smash{\hat{\mathcal{P}}}$ that do not belong to the clipped surface $\smash{\tilde{\mathcal{S}}}$, i.e., originate from the edges of the polyhedron $\mathcal{P}$.
For a correct parameterization of the edges of $\smash{\tilde{\mathcal{S}}}$, which consist of conic section arcs owing to the quadratic nature of $\mathcal{S}$, a rational B\'ezier curve formulation is employed. Such a parameterization reproduces conic section arcs exactly~\cite{Farin}, and is defined by two end-points, a control point~$\mathbf{x}^*$, and a scalar weight~$w$. 
Following the introduction of these two parameterizations, the moment integral is then split into three contributions:
\begin{itemize}[leftmargin=4em]
    \item[$\boldsymbol{\mathcal{M}}^{\hat{\mathcal{P}}_1}$:] This first term is an approximation of the integral of $\boldsymbol{\Psi}_{\mathcal{F}_i}$ (the planar face contribution) using the linear parameterization $\bar{\mathbf{x}}$. 
    \item[$\boldsymbol{\mathcal{M}}^{\hat{\mathcal{P}}_2}$:] This second term is an approximation of the integral of $\boldsymbol{\Psi}_{\mathcal{S}}$ (the quadratic surface contribution) using the linear parameterization $\bar{\mathbf{x}}$.
    \item[$\boldsymbol{\mathcal{M}}^{\hat{\mathcal{P}}_3}$:] This last corrective term accounts for the difference between a correct parameterization of the edge and its linear approximation.
\end{itemize}
This splitting merely aims to simplify the integrals and the resulting closed-form expressions of the moments, as well as to minimize numerical round-off errors when estimating these quantities. It does not impact the end result, i.e., the analytical derivation of the moments, which remains exact. 
The integration domains corresponding to each contribution to the moments are illustrated in Fig.~\ref{fig:M123}.
Ultimately, the moments $\smash{\boldsymbol{\mathcal{M}}^{\hat{\mathcal{P}}}}$ are expressed as:
\begin{equation}
\boldsymbol{\mathcal{M}}^{\hat{\mathcal{P}}} = \boldsymbol{\mathcal{M}}^{\hat{\mathcal{P}}_1} + \boldsymbol{\mathcal{M}}^{\hat{\mathcal{P}}_2} + \boldsymbol{\mathcal{M}}^{\hat{\mathcal{P}}_3} = \sum\limits_{i=1}^{n_{\!\mathcal{F_{\phantom{k}\!\!}}}} \left( \boldsymbol{\mathcal{M}}^{\hat{\mathcal{P}}_1}_i + \boldsymbol{\mathcal{M}}^{\hat{\mathcal{P}}_2}_i + \boldsymbol{\mathcal{M}}^{\hat{\mathcal{P}}_3}_i \right) \, .
\end{equation}
After integrating and simplifying each term, the facewise contributions to the moments can be written as
\begin{align}
    \boldsymbol{\mathcal{M}}^{\hat{\mathcal{P}}_1}_i &= \sum_{j=1}^{n_{\hat{\mathcal{F}}_i}} \mathcal{A}(\mathbf{x}_{i, j, 0}, \mathbf{x}_{i, j, 1}, \mathbf{x}_{i, \mathrm{ref}})\boldsymbol{\mathcal{B}}^{(1)}(\mathbf{x}_{i, j, 0}, \mathbf{x}_{i, j, 1}, \mathbf{x}_{i, ref}) \, ,\label{eq:mp1i} \\
    \boldsymbol{\mathcal{M}}^{\hat{\mathcal{P}}_2}_i &= \sum_{j=1}^{n_{\hat{\mathcal{F}}_i}} \mathbbm{1} ^{\partial \tilde{\mathcal{S}}}_{i,j} \mathcal{A}(\mathbf{x}_{i, j, 0}, \mathbf{x}_{i, j, 1}, \mathbf{0})\boldsymbol{\mathcal{B}}^{(2)}(\mathbf{x}_{i, j, 0}, \mathbf{x}_{i, j, 1}) \, ,\\
    \boldsymbol{\mathcal{M}}^{\hat{\mathcal{P}}_3}_i &= \sum_{j=1}^{n_{\hat{\mathcal{F}}_i}} \mathbbm{1} ^{\partial \tilde{\mathcal{S}}}_{i,j} \mathcal{A}(\mathbf{x}_{i, j, 0}, \mathbf{x}_{i, j, 1}, \mathbf{x}_{i, j}^*) \boldsymbol{\mathcal{B}}^{(3)}(w_{i, j}, \mathbf{x}_{i, j, 0}, \mathbf{x}_{i, j, 1}, \mathbf{x}_{i, j}^*) \, .
\end{align}
In these expressions, $\smash{\mathbf{x}_{i,\mathrm{ref}}}$ is a reference point required to belong to the face $\smash{\hat{\mathcal{F}}_i}$, which is arbitrarily chosen as the first vertex of that face. The functions $\smash{\mathcal{A}}$ and $\smash{\boldsymbol{\mathcal{B}}^{(\square)}}$ are derived in \citet{Evrard} for the case of a paraboloid surface $\mathcal{S}$, and therefore are not repeated here. They consist of real multivariate polynomial functions of the coordinates of the edge end-points, except for the dependency on the rational B\'ezier weight $w_{i, j}$ which involves trigonometric functions and rational polynomials. 

\subsection{Extension to Quadratic Cylinders} \label{sec:new_work}

The strategy outlined in Section~\ref{sec:fondation} is now applied to the case of a cylindrical quadratic surface $\mathcal{S}$, in lieu of the paraboloid surface considered in \citet{Evrard}. While some intermediate results from \cite{Evrard} still hold, several definitions and equations require adaptation to this new type of quadratic surface. The following sections address these differences and derive new closed-form expressions for the first moments of a polyhedron clipped by a quadratic cylinder.

\subsubsection{Problem Definition}\label{def:halfspace} Consider the same arbitrary polyhedron $\mathcal{P}$ as in Section~\ref{sec:fondation} with planar polygonal faces $\mathcal{F}_i$, however $\mathcal{Q}$ is now the region that is `inside' the infinite cylinder $\mathcal{S}$, both implicitly defined as
\begin{align}
    \mathcal{Q} & = \{ \mathbf{x} \in \mathbb{R}^3 : \phi(\mathbf{x}) \le 0\} \, , \\
    \mathcal{S} & = \{ \mathbf{x} \in \mathbb{R}^3 : \phi(\mathbf{x}) = 0  \} \, , \label{eq:defS}
\end{align}
with
\begin{equation}
    \phi : \left\{\begin{array}{lcl} 
        \mathbb{R}^3 & \!\to & \mathbb{R} \\
        \mathbf{x} & \!\mapsto & \beta y^2 + z^2 - \cylradius^2 
    \end{array}\right. \ , (\cylradius,\beta) \in (\mathbb{R}^+_* \times \mathbb{R}_*) \, . \label{def:phi} 
\end{equation}
Based on the value of $\beta$, the surface $\mathcal{S}$ can either be an elliptic cylinder ($\beta > 0$),  or a hyperbolic cylinder ($\beta < 0$).
This analysis is restricted to polyhedra that lie above the $xy$-plane, i.e., the domain $\Omega$ under consideration is the $z \ge 0$ halfspace.
Under this restriction, the $z$-coordinate of every point of $\tilde{\mathcal{S}}$ is positive, hence the surface can be parametrized by a single-valued function of $x$ and $y$.
\subsubsection{First Divergence Theorem Application and Projection}
The definition of $\mathcal{S}$ provided by Eqs.~\eqref{eq:defS} and~\eqref{def:phi} does not affect the results of the first application of the divergence theorem detailed in Section~\ref{sec:fondation}. 
Similarly, the resulting two-dimensional integral can be split into the contribution of the clipped planar faces $\smash{\cup_i\hat{\mathcal{F}}_i}$ and of the clipped cylinder surface $\tilde{\mathcal{S}}$, and the first moments are then given by
\begin{equation}
    \boldsymbol{\mathcal{M}}^{\hat{\mathcal{P}}} = \int_{\tilde{\mathcal{S}}} \boldsymbol{\Phi}(\mathbf{x}) (\mathbf{n}_{\tilde{\mathcal{S}}} \cdot \mathbf{e}_z) \, \mathrm{d}\mathbf{a} + \sum_{i=1}^{n_{\mathcal{F}}} \int_{\hat{\mathcal{F}}_i} \boldsymbol{\Phi}(\mathbf{x}) (\mathbf{n}_{i} \cdot \mathbf{e}_z) \,\mathrm{d}\mathbf{a} \, , \label{eq:div1}
\end{equation}
with $\smash{\mathbf{n}_{\tilde{\mathcal{S}}}}$ the outward-pointing normal to $\smash{\tilde{\mathcal{S}}}$.
To conform to the new type of quadratic surface under consideration, the parameterization of ${\tilde{\mathcal{S}}}$ becomes
\begin{equation}
    \tilde{\mathcal{S}} = 
    \left\{ \begin{bmatrix}[1.05] 
    x \\
    y \\
    \sqrt{\cylradius^2 - \beta y^2}
    \end{bmatrix} , \quad (x,y) \in \tilde{\mathcal{S}}^{\perp}
    \right\} \, . \label{def:paramS}
\end{equation}
The rest of the reasoning summarized in Section~\ref{sec:fondation} and detailed in \cite{Evrard} still holds. In particular, the term $\|\nabla\phi(\mathbf{x)\|}$ being the determinant of the parameterization~\eqref{def:paramS} of $\tilde{\mathcal{S}}$, it cancels out the denominator of $(\mathbf{n}_{\tilde{\mathcal{S}}} \cdot \mathbf{e}_z)$ leading to a non-rational integrand in the first integral of Eq.~\eqref{eq:div1}. The first moments are thus given as in Eq.~\eqref{eq:divperp} with $\boldsymbol{\Phi}_{\mathcal{S}}$ defined as
\begin{equation}
    \boldsymbol{\Phi}_{\!\mathcal{S}}(x,y)  = 
    \begin{bmatrix}[1.1]
    \sqrt{\cylradius^2 - \beta y^2} \\  
    x \sqrt{\cylradius^2 - \beta y^2} \\  
    y \sqrt{\cylradius^2 - \beta y^2} \\  
    \frac{1}{2}\left(\cylradius^2 - \beta  y^2\right)
    \end{bmatrix}\label{eq:primitives} \, .
\end{equation}

\subsubsection{Second Divergence Theorem Application}\label{sec:cylinder_case}
Applying the divergence theorem a second time, one recovers Eq.~\eqref{eq:vol_after_div7} with
\begin{equation}
    \boldsymbol{\Psi}_{\!\mathcal{S}}(x,y)  = 
    \begin{bmatrix}[1.1]
    x\sqrt{\cylradius^2 - \beta y^2} \\  
    \frac{1}{2}x^2 \sqrt{\cylradius^2 - \beta y^2} \\  
    xy \sqrt{\cylradius^2 - \beta y^2} \\  
    \frac{1}{2}x\left(\cylradius^2 - \beta y^2\right)
    \end{bmatrix}\label{eq:psy_cylinder}
\end{equation}
The same linear and rational B\'ezier parameterizations as introduced in Section~\ref{sec:fondation} are used to reformulate $\smash{\boldsymbol{\mathcal{M}}^{\hat{\mathcal{P}}}}$ into the sum of three contributions $\smash{\boldsymbol{\mathcal{M}}^{\hat{\mathcal{P}}_1}}$, $\smash{\boldsymbol{\mathcal{M}}^{\hat{\mathcal{P}}_2}}$ and $\smash{\boldsymbol{\mathcal{M}}^{\hat{\mathcal{P}}_3}}$.
The term $\smash{\boldsymbol{\mathcal{M}}^{\hat{\mathcal{P}}_1}}$ remains unchanged from \citet{Evrard}, since its expression only depends on the integrands~$\boldsymbol{\Psi}_{\mathcal{F}_i}$. It is therefore given by Eq.~\eqref{eq:mp1i} with $\smash{\mathcal{A}}$ and $\smash{\boldsymbol{\mathcal{B}}^{(1)}}$ defined as
\begin{equation}
    \mathcal{A} : \left\{ \begin{array}{ccl} \mathbb{R}^3\times\mathbb{R}^3\times\mathbb{R}^3 & \to & \mathbb{R} \\ \left(\mathbf{x}_{a},\mathbf{x}_{b},\mathbf{x}_{c}\right) & \mapsto & \frac{1}{2} \left(x_a(y_b-y_c) + x_b(y_c-y_a) + x_c(y_a-y_b) \right) \end{array}\right. \, , \label{eq:triangle_area}
\end{equation}
and
\begin{equation}
    \boldsymbol{\mathcal{B}}^{(1)} : \left\{ \begin{array}{ccl} \mathbb{R}^3\times\mathbb{R}^3\times\mathbb{R}^3 & \!\!\!\! \to \!\!\!\! & \mathbb{R}^4 \\ 
    \left(\mathbf{x}_{a},\mathbf{x}_{b},\mathbf{x}_{c}\right) & \!\!\!\! \mapsto \!\!\!\! &  \frac{1}{12} {\text{\footnotesize$
\begin{bmatrix} 4 \left( z_{a} + z_{b} + z_{c} \right) \\  
    \left( z_{a} + z_{b} + z_{c} \right) \left( x_{a} + x_{b} + x_{c} \right) + x_{a} z_{a} + x_{b} z_{b} + x_{c} z_{c} \\
     \left( z_{a} + z_{b} + z_{c} \right) \left( y_{a} + y_{b} + y_{c} \right) + y_{a} z_{a} + y_{b} z_{b} + y_{c} z_{c} \\ z_{a}^2 + z_{b}^2 +
  z_{c}^2 + z_{b} z_{c} +
  z_{a} z_{b} + z_{a} z_{c}\end{bmatrix}$}} \end{array}\right. \!\! .
\end{equation}
However, the new expression of $\boldsymbol{\Psi}_{\!\mathcal{S}}$ given in Eq.~\eqref{eq:psy_cylinder} and derived for the case of a quadratic cylinder does not enable the direct derivation of $\smash{\boldsymbol{\mathcal{M}}^{\hat{\mathcal{P}}_2}}$ and $\smash{\boldsymbol{\mathcal{M}}^{\hat{\mathcal{P}}_3}}$ as done in \cite{Evrard}. These terms require reformulation, as detailed in the following two subsections.

\subsubsection{New Second Contribution \texorpdfstring{$\boldsymbol{\mathcal{M}}^{\hat{\mathcal{P}}_2}$}{MP2}
}
For the case of a paraboloid surface $\mathcal{S}$ in \cite{Evrard}, the moment contribution $\smash{\boldsymbol{\mathcal{M}}^{\hat{\mathcal{P}}_2}}$ is defined as
\begin{equation}
\boldsymbol{\mathcal{M}}^{\hat{\mathcal{P}}_2} = \sum_{i=1}^{n_\mathcal{F}} \boldsymbol{\mathcal{M}}^{\hat{\mathcal{P}}_2}_i = \sum_{i=1}^{n_\mathcal{F}} \sum_{j=1}^{n_{\partial \hat{\mathcal{F}}_i}} -\mathbbm{1}_{i,j}^{\partial \tilde{\mathcal{S}}}\int_0^1  \, \boldsymbol{\Psi}_{\mathcal{S}}(\overline{x}_{i,j}(t), \overline{y}_{i,j}(t)) (y_{i,j,1} - y_{i,j,0}) \, \mathrm{d}t \, ,
\end{equation}
with $\smash{\mathbf{x}_{i,j,0} \equiv \mathbf{x}_{i,j}(0)}$ and $\smash{\mathbf{x}_{i,j,1} \equiv \mathbf{x}_{i,j}(1)}$. 
In the present work, which considers quadratic cylinders, this would lead the volume contribution of that integral (i.e., the contribution of the first component of $\smash{\boldsymbol{\Psi}_{\mathcal{S}}}$) to read as
\begin{equation*}
\int_0^1 \bar{x}_{i,j}(t) \sqrt{\cylradius^2-\beta ((1-t) y_{i, j, 0} + t y_{i,j,1})^2} \; (y_{i,j,1}-y_{i, j, 0}) \, \mathrm{d}t \, ,
\end{equation*} 
in which the presence of the square root nullifies the benefit of using a linear approximation of the curved integration domain.
Instead, the linear interpolation of $y_{i,j}$ is substituted by a linear parameterization of $\smash{t \mapsto (\cylradius^2-\beta y_{i,j}^2(t))^{\sfrac{1}{2}}}$. The previous integral thus becomes
\begin{equation*}
\int_0^1 \bar{x}_{i,j}(t) \left( (1-t)\sqrt{\cylradius^2-\beta y_{i, j, 0}^2} + t \sqrt{\cylradius^2-\beta y_{i, j, 1}^2}\right) (y_{i,j,1}-y_{i, j, 0}) \, \mathrm{d} t \ .
\end{equation*}
Using the fact that
\begin{equation}
-\mathbbm{1}^{\partial\tilde{\mathcal{S}}}_{i,j} \neq 0 \ \Rightarrow \ \left\{\begin{matrix}
    z_{i,j,0} = \sqrt{\cylradius^2-\beta \, y_{i,j,0}^2} \\
    z_{i,j,1} = \sqrt{\cylradius^2-\beta \, y_{i,j,1}^2} 
\end{matrix}\right. \, ,
\end{equation}
this integral can be rewritten as
\begin{equation}
\int_0^1 \bar{x}_{i,j}(t) ((1-t) z_{i, j, 0} + tz_{i,j,1}) (y_{i,j,1}-y_{i, j, 0}) \, \mathrm{d}t = \int_0^1\bar{x}_{i,j}(t)\bar{z}_{i,j}(t)\bar{y}_{i,j}^\prime(t) \, \mathrm{d}t \, .
\end{equation}
Following the same logic for the other moment components, $\boldsymbol{\mathcal{M}}^{\hat{\mathcal{P}}_2}$ is redefined as
\begin{equation}
\boldsymbol{\mathcal{M}}^{\hat{\mathcal{P}}_2} = \sum_{i=1}^{n_\mathcal{F}} \boldsymbol{\mathcal{M}}^{\hat{\mathcal{P}}_2}_i = \sum_{i=1}^{n_\mathcal{F}} \sum_{j=1}^{n_{\partial \hat{\mathcal{F}}_i}} -\mathbbm{1}_{i,j}^{\partial \tilde{\mathcal{S}}} \int_0^1  \boldsymbol{\Psi}_{\mathcal{S}}^{\mathcal{R}}(\overline{\mathbf{x}}_{i,j}, t) \overline{y}_{i,j}^\prime(t) \, \mathrm{d}t \, , \label{eq:newm2}
\end{equation}
where
\begin{equation}
\boldsymbol{\Psi}_{\mathcal{S}}^{\mathcal{R}} : \left\{\begin{array}{lcl} 
        ([0,1] \!\to \mathbb{R}^3) \times [0, 1] & \!\to & \mathbb{R}^4 \\
        (\bar{\mathbf{x}}, t) & \!\mapsto & \begin{bmatrix}
            \bar{x}(t) \, \bar{z}(t) \\
            \frac{1}{2}\bar{x}(t)^2 \, \bar{z}(t) \\
            \bar{x}(t) \, \bar{y}(t) \, \bar{z}(t) \\
            \frac{1}{2} \bar{x}(t) \, \bar{z}(t)^2
        \end{bmatrix} \end{array}\right. \, .
\end{equation}
With this new definition and in the case of a quadratic cylinder $\mathcal{S}$ as defined by Eqs.~\eqref{eq:defS} and~\eqref{def:phi}, $\smash{\boldsymbol{\mathcal{M}}^{\hat{\mathcal{P}}_2}_{i}}$ can be expressed as
\begin{align}
    \boldsymbol{\mathcal{M}}^{\hat{\mathcal{P}}_2}_{i} & = \sum\limits_{j = 1}^{n_{\partial\hat{\mathcal{F}}_i}} \mathbbm{1}^{\partial\tilde{\mathcal{S}}}_{i,j} (y_{i,j,0}-y_{i,j,1})\boldsymbol{\mathcal{B}}^{(2)}\left({\mathbf{x}}_{i,j,0},{\mathbf{x}}_{i,j,1}\right) \, ,
\end{align}
where the operator $\smash{\boldsymbol{\mathcal{B}}^{(2)} : \mathbb{R}^3\times\mathbb{R}^3 \to \mathbb{R}^4}$ reads as
\begin{equation}
    \boldsymbol{\mathcal{B}}^{(2)}\left(\mathbf{x}_{a},\mathbf{x}_{b}\right) = \frac{1}{24}{\text{\small$
\begin{bmatrix} 4 \left( x_{a} (2z_a+z_b) + x_{b} (z_a+2z_b) \right) \\  
     \left((z_a+ z_b)(x_a+x_b)^2
    +2x_a^2z_a +2x_b^2z_b\right)
     \\
    2 \left((x_a+x_b)(y_a+y_b)(z_a+ z_b)
    +2x_ay_az_a +2x_by_bz_b\right)
    \\ 
     \left((z_a+ z_b)^2(x_a+x_b)
    +2x_az^2_a +2x_bz_b^2\right)
    \end{bmatrix}$}} \, .
\end{equation}

\subsubsection{New Third Contribution \texorpdfstring{$\boldsymbol{\mathcal{M}}^{\hat{\mathcal{P}}_3}$}{MP3}
}
The new definition of $\boldsymbol{\mathcal{M}}^{\hat{\mathcal{P}}_2}$ introduced in Eq.~\eqref{eq:newm2} implies a new definition of the third contribution $\boldsymbol{\mathcal{M}}^{\hat{\mathcal{P}}_3}$, which then reads as
\begin{align}
    \boldsymbol{\mathcal{M}}^{\hat{\mathcal{P}}_3}_{i} &= \sum\limits_{i = 1}^{n_{\mathcal{F}}}\int_{0}^{1} \sum\limits_{j = 1}^{n_{\partial\hat{\mathcal{F}}_i}} -\mathbbm{1}^{\partial\tilde{\mathcal{S}}}_{i,j} \left[\left(\boldsymbol{\Psi}_{\mathcal{F}_i}(\bar{x}_{i,j}(t),\bar{y}_{i,j}(t))-\boldsymbol{\Psi}^{\mathcal{R}}_{\mathcal{S}}(\bar{\mathbf{x}}_{i,j}, t)\right) \bar{y}^\prime_{i,j}(t) \right. \\ 
    & \quad\quad\quad\quad\quad\quad \left. + \left( \phantom{\boldsymbol{\Psi}^{\mathcal{R}}_{\mathcal{S}}\!\!\!\!\!\!\!\!\!\!}\boldsymbol{\Psi}_{\mathcal{S}}(x_{i,j}(t),y_{i,j}(t)) - \boldsymbol{\Psi}_{\mathcal{F}_i}(x_{i,j}(t),y_{i,j}(t)) \right)y^\prime_{i,j}(t)\right] \ \mathrm{d}t \, .\nonumber
\end{align}
The methodology to derive a closed-form expression remains otherwise similar.
After simplification, an expression similar to that of \citet{Evrard} is found,
\begin{align}
    \boldsymbol{\mathcal{M}}^{\hat{\mathcal{P}}_3}_{i} & = \sum\limits_{j = 1}^{n_{\partial\hat{\mathcal{F}}_i}} -\mathbbm{1}^{\partial\tilde{\mathcal{S}}}_{i,j} \mathcal{A}^\dagger\left({\mathbf{x}}_{i,j,0},{\mathbf{x}}_{i,j,1},\mathbf{x}_{i,j}^{\star}\right)\boldsymbol{\mathcal{B}}^{(3)}\left(w_{i,j},{\mathbf{x}}_{i,j,0},{\mathbf{x}}_{i,j,1},\mathbf{x}_{i,j}^{\star}\right) \, , \label{eq:m3k}
\end{align}
where $\mathcal{A}^\dagger$ is a variation of the operator $\mathcal{A}$ introduced in \cite{Evrard} and given in~Eq.~\eqref{eq:triangle_area}, which computes the signed projected area of a triangle, but with the projection now occurring onto the $yz$-plane instead of the $xy$-plane, i.e.,
\begin{equation}
    \mathcal{A}^\dagger : \left\{ \begin{array}{ccl} \mathbb{R}^3\times\mathbb{R}^3\times\mathbb{R}^3 & \to & \mathbb{R} \\ \left(\mathbf{x}_{a},\mathbf{x}_{b},\mathbf{x}_{c}\right) & \mapsto & \frac{1}{2} \left(y_a(z_b-z_c) + y_b(z_c-z_a) + y_c(z_a-z_b) \right) \end{array}\right. \, . \label{eq:triangle_area2}
\end{equation}
The operator $\smash{\boldsymbol{\mathcal{B}}^{(3)}: \mathbb{R}^{+} \times \mathbb{R}^3 \times \mathbb{R}^3 \times \mathbb{R}^3 \to \mathbb{R}^4}$ is given in full in~\ref{apdx:M3}.

\subsubsection{Generalization to the Full \texorpdfstring{$\mathbb{R}^3$}{R3} Space}\label{sec:generalization}
\begin{figure}[h] \centering
        \includegraphics[width=\linewidth,trim=0 700 0 700,clip=true]{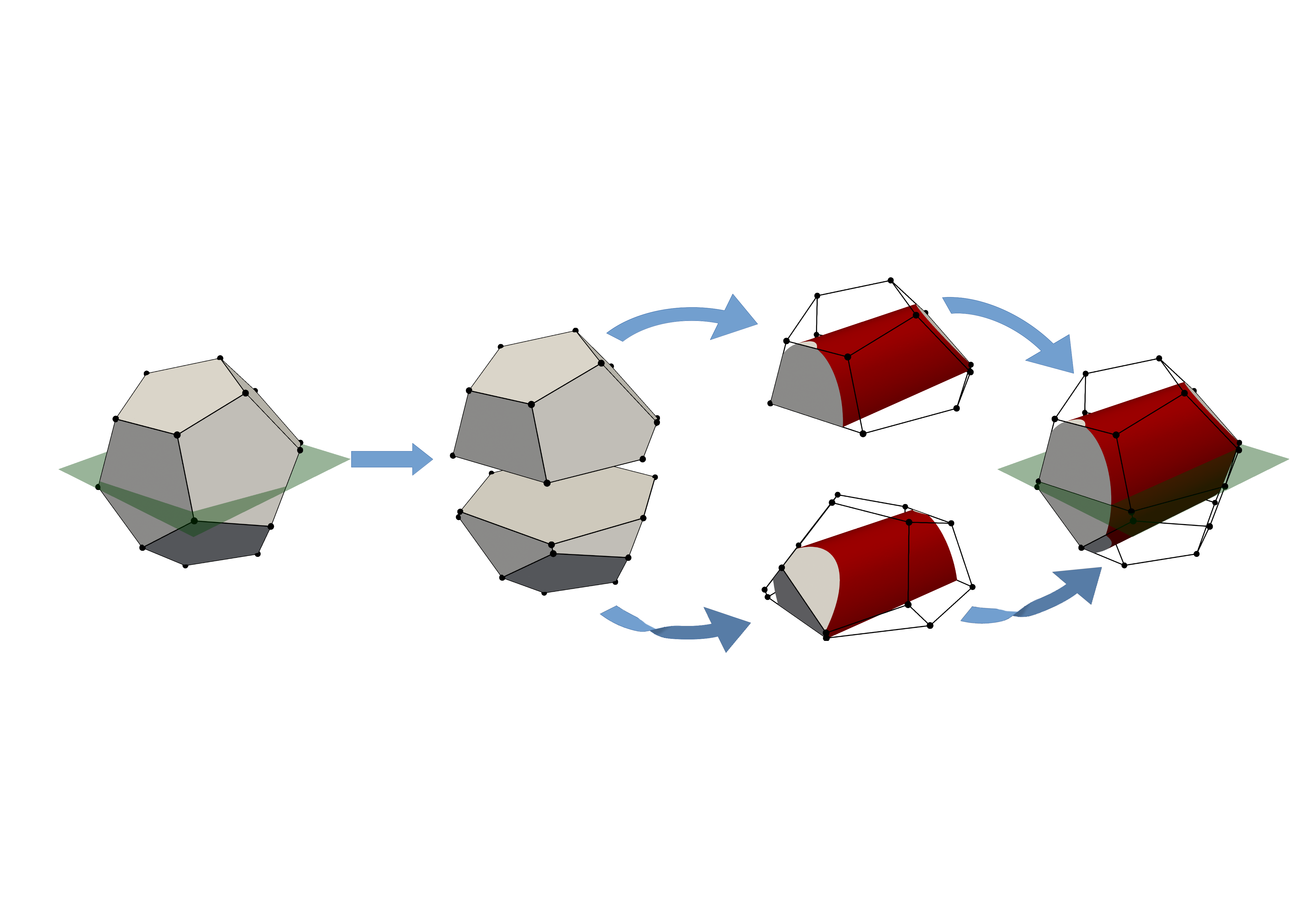}
    \caption{When a polyhedron that intersects the $xy$-plane is to be clipped by a quadratic cylinder centered around $\mathbf{e}_x$, the polyhedron is first clipped by the $xy$-plane and then each half is considered independently according to Section~\ref{def:halfspace}.}
    \label{fig:merging}
\end{figure}

In Section~\ref{def:halfspace}, the analysis was restricted to polyhedra that lie above the $xy$-plane. To compute the first moments of a polyhedron clipped by a quadratic cylinder when the polyhedron does not entirely lie above the $xy$-plane, the polyhedron is first clipped by the $xy$-plane and the approach described above is applied to both `top' and `bottom' clipped regions independently, as illustrated in Fig.~\ref{fig:merging}. Note that the bottom clipped region must be rotated by $180$ degrees about $\mathbf{e}_x$ prior to using the closed-form expressions derived in this section.

\subsection{On Floating-Point Arithmetic and Robustness}\label{sec:robustness}
As is the case in \cite{Evrard}, some intersection configurations can be ill-posed for calculating the moments of the clipped polyhedron. This includes cases in which the surface $\mathcal{S}$ is tangential to one or several edges of $\mathcal{P}$, or contains at least one vertex of $\mathcal{P}$. The perturbation procedure introduced in \cite{Evrard} is applied to handle such degenerate cases, using the same values of the distance and angular thresholds ($\smash{\epsilon_{\text{corner}} = 10^{2}\times\epsilon_{64}}$, $\smash{\epsilon_{\text{tangent}} = 10^6\times\epsilon_{64}}$), and nudging distance ($\smash{\epsilon_{\text{nudge}} = 10^{10}\times\epsilon_{128}}$)\footnote{The constants $\epsilon_{64} = 2^{-52}$ and $\epsilon_{128} = 2^{-112}$ correspond to the relative precision of double and quadruple precision floating-point arithmetics, respectively.}. For all test-cases considered in this paper, this perturbation procedure was found to entirely mitigate the issue of ill-posed geometry and intersection configurations.
 
\section{Cylindrical Interface Reconstruction Method}
\label{Methods2}
In this section, the methods for reconstructing and detecting ligament-like structures are introduced in detail, and their implementation into the NGA2~\cite{Desjardins} flow solver is discussed. Section~\ref{Methods:PCIC} describes the PCIC reconstruction algorithm, Section~\ref{Methods:lig_det} introduces ligament detection, Section~\ref{Methods:tip_det} details ligament tip detection, and Section~\ref{Methods:NGA2} discusses flow solver specifics. Note that for simplicity these subsections describe PCIC's application to liquid ligaments; however, the proposed methods are easily extended to gas ligaments as well. 

\subsection{PCIC Principal Curve Fitting}
\label{Methods:PCIC}
PCIC only uses straight circular cylinders ($\beta = 1$) to approximate fluid ligaments. Such a cylinder is defined by a radius $r$, an origin $\mathbf{p}$, and an orientation vector $\mathbf{v}$. To inform the reconstruction, a $5\times5\times5$ stencil of liquid phase barycenters is considered. Note that only those barycenters that belong to the same connected-component labeling (CCL)~\cite{He} ligament structure as the stencil's center cell are used in the following process (See Section~\ref{Methods:lig_det} for CCL information). From these barycenters, it is necessary to extract $\mathbf{p}$ and $\mathbf{v}$ for the cylinder in the stencil's center cell. A second-order principal curve is fitted to the barycenter data using an iterative approach based on the method of \citet{Hastie}. The process is initialized with a Principal Component Analysis (PCA) of the barycenter data. The principal eigenvector is taken as an initial linear fit to the data. Three control points are then placed on this initial curve: one at each end point (defined by the range of barycenter projections onto the line), and one at the midpoint. All of the barycenter data are then projected onto the curve and grouped by closest control point. Each group of barycenters is averaged with a volume fraction weighting
\begin{equation}
\begin{aligned}
\mathbf{y}_i^{m+1}=\frac{1}{\sum_{j=1}^{n_i^m}{\alpha_{ij}^m}}\sum_{j=1}^{n_i^m}\mathbf{x}_{ij}^m{\alpha_{ij}^m},
\end{aligned}
\end{equation}
where $\mathbf{y}_i^{m+1}$ is the $i_{th}$ control point at iteration $m+1$, $n_i^m$ is the number of barycenters projecting to $\mathbf{y}_i^m$, $\mathbf{\alpha}_{i}^m$ are the volume fractions of the barycenters projecting to $\mathbf{y}_i^m$, and $\mathbf{x}_i^m$ are the barycenters projecting to $\mathbf{y}_i^m$. This results in three updated control points $\mathbf{y}_1^{m+1}$, $\mathbf{y}_2^{m+1}$, and $\mathbf{y}_3^{m+1}$ defining a quadratic curve. This process is repeated until the control points converge to within a maximum squared Euclidean distance tolerance of \num{e-5} (normalized by the square of the mesh size) or after a maximum of $10$ iterations. 

With the three final control points, the parameterized expression for the quadratic curve can be defined as
\begin{align} \label{eq:quad}
\mathbf{z} & =\mathbf{A}t^2+\mathbf{B}t+\mathbf{C} \\
\mathbf{A} & =\frac{\mathbf{y}_1}{\left(t_1-t_2\right)\left(t_1-t_3\right)}+\frac{\mathbf{y}_2}{\left(t_2-t_1\right)\left(t_2-t_3\right)}+\frac{\mathbf{y}_3}{\left(t_3-t_1\right)\left(t_3-t_2\right)} \\
\mathbf{B} & =-\frac{\mathbf{y}_1\left(t_2+t_3\right)}{\left(t_1-t_2\right)\left(t_1-t_3\right)}-\frac{\mathbf{y}_2\left(t_1+t_3\right)}{\left(t_2-t_1\right)\left(t_2-t_3\right)}-\frac{\mathbf{y}_3\left(t_1+t_2\right)}{\left(t_3-t_1\right)\left(t_3-t_2\right)} \\
\mathbf{C} & =\frac{\mathbf{y}_1t_2t_3}{\left(t_1-t_2\right)\left(t_1-t_3\right)}+\frac{\mathbf{y}_2t_1t_3}{\left(t_2-t_1\right)\left(t_2-t_3\right)}+\frac{\mathbf{y}_3t_1t_2}{\left(t_3-t_1\right)\left(t_3-t_2\right)},
\end{align}
where $\mathbf{z}$ is a point on the curve, $t$ is a parameter between $0$ and $1$ representing position along the curve, and $t_1$, $t_2$, and $t_3$ are the parameter values at control points $\mathbf{y_1}$, $\mathbf{y_2}$, and $\mathbf{y_3}$, respectively. The point on the curve closest to the stencil's center can then be found analytically. Let $\mathbf{x}_c$ be the stencil center. Then, the closest point on the curve to $\mathbf{x}_c$ is simply at the solution of the minimization
\begin{equation}
\begin{aligned}
t_{min}=\arg \min_t\left|\left|\left(\mathbf{A}t^2+\mathbf{B}t+\mathbf{C}\right)-\mathbf{x}_c\right|\right|^2,
\end{aligned}
\end{equation}
which can be written as the cubic equation
\begin{align}
0 & =at_{min}^3+bt_{min}^2+ct_{min}+d \\
a & =2\mathbf{A}\cdot\mathbf{A} \\
b & =3\mathbf{A}\cdot\mathbf{B} \\
c & =2\mathbf{A}\cdot\left(\mathbf{C}-\mathbf{x}_c\right)+\mathbf{B}\cdot\mathbf{B} \\
d & =\mathbf{B}\cdot\left(\mathbf{C}-\mathbf{x}_c\right),
\end{align}
where $t_{min}$ is the parameter at the desired point. The cubic equation is solved analytically for $t_{min}$, which can provide up to three real solutions and is guaranteed to have at least one; the solution resulting in a point closest to $\mathbf{x}_c$ is chosen as $t_{min}$. The origin of the cylinder being reconstructed in the stencil's center cell, $\mathbf{p}$, is set to the point on the curve at $t_{min}$ with
\begin{equation}
\begin{aligned}
\mathbf{p}=\mathbf{A}t_{min}^2+\mathbf{B}t_{min}+\mathbf{C}.
\end{aligned}
\end{equation}
The derivative of Eq.~\eqref{eq:quad} with respect to $t$ is then calculated and evaluated at $t_{min}$, which gives the tangent vector to the curve at $\mathbf{p}$. This tangent vector is used for the cylinder's orientation (i.e., the longitudinal direction of the cylinder) $\mathbf{v}$, which is defined as
\begin{equation}
\begin{aligned}
\mathbf{v}=\frac{2\mathbf{A}t_{min}+\mathbf{B}}{\left|\left|2\mathbf{A}t_{min}+\mathbf{B}\right|\right|}.
\end{aligned}
\end{equation}

For a schematic of the preceding process, refer to Fig.~\ref{PCIC_fit} in Section~\ref{sec:deform}. With the origin and orientation defined, only the radius is left to be determined. An iterative bisection method is used to select a cylinder radius that conserves the volume fraction in the stencil's center cell to within a tolerance of \num{e-14}. The search space is defined from a radius of $0$ to a radius that corresponds to a volume fraction of $1$ in the center cell, based on the given $\mathbf{p}$ and $\mathbf{v}$. This search space is then progressively halved by evaluating the volume fraction at the midpoint of the current search space (using the closed-form expressions derived in Section~\ref{sec:cylinder_case}) and comparing it to the target volume fraction. If the target is greater than the current value, the half with radii greater than the current midpoint is chosen as the new search space. If the target is less than the current value, the other half is used.

\subsection{Ligament Detection}
\label{Methods:lig_det}
PCIC, as designed in the preceding subsection, should only be used when reconstructing interfacial geometries that are approximated well by sub-grid or near sub-grid straight circular cylinders. Thin liquid ligaments are an excellent example of such a geometry and are the intended use case for PCIC, making it necessary to automatically identify these structures near the mesh size in VOF simulations. This is done with an algorithm based on connected-component labeling (CCL)~\cite{He}, which detects and labels connected regions, or clusters, within an array of given binary values. CCL maintains a list of the cells in each separate contiguous region. This approach was used successfully for the similar task of thin film detection in work by \citet{Han3}. 

The first step is to define a binary array based on a thickness criterion. A representation of thickness is calculated as twice the liquid volume in a $7\times7\times7$ stencil divided by the interfacial surface area in the same stencil
\begin{equation}
\begin{aligned}
h_i=\frac{2\sum_{j=1}^{7^3}\alpha_{ij}\mathscr{V}_{ij}}{\sum_{j=1}^{7^3}\mathscr{A}_{ij}},
\end{aligned}
\end{equation}
where $h_i$ is the stencil's thickness representation, $\mathbf{\alpha}_{i}$ are the volume fractions in the stencil cells, $\mathbf{\mathscr{V}}_i$ are the volumes of the stencil cells, and $\mathbf{\mathscr{A}}_i$ are the interface surface areas contained within the stencil cells. These surface areas are calculated with advected interface polygons from the prior timestep; this process is described in \citet{Han3}. With such a definition, if the volume and surface moments correspond to those of a straight circular cylinder, then $h_i$ is the radius of the cylinder. The CCL binary array is then defined by
\begin{equation}
\mathscr{X}_i =
\begin{cases}
    1, & h_i \leq \Delta x \\
    0, & h_i > \Delta x
\end{cases}
\end{equation}
where $\mathscr{X}_i$ is a cell's CCL binary value and $\Delta x$ is the minimum grid size in a simulation. Note that this process only creates contiguous labeled regions for general thin structures and does not yet differentiate between films and ligaments. Also note that ligaments with a diameter of up to $2\Delta x$ can be detected. This is a desired feature, so that a clean transition can be made between a PLIC and a PCIC reconstruction.

The liquid moments of inertia are also calculated in a $5\times5\times5$ stencil. The eigenvalues of the corresponding inertia tensor are calculated to determine the principal moments of inertia, providing information on the shape of the liquid structure contained in the stencil. To determine if a center cell contains a "ligament-like" structure, a criteria is used defined as
\begin{equation}
\ell_i =
\begin{cases}
    1, & \left(1.5\lambda_{i2} > \lambda_{i3}\right)\text{ and }\left(\lambda_{i2} > 1.5\lambda_{i1}\right) \\
    0, & \text{otherwise}
\end{cases}
\end{equation}
where $\ell_i$ is a binary indicator for ligament-likeness, and $\lambda_{i3}$, $\lambda_{i2}$, and $\lambda_{i1}$, are a cell's principal liquid volume moments of inertia, in descending order.

Finally, each detected CCL region is analyzed as follows. The minimum thickness of each structure is found from its list of $h_i$ values, and the percentage of "ligament-like" cells it contains is calculated. If a structure has a minimum $h_i$ value of at most $0.8\Delta x$ and contains at least $90\%$ ligament-like cells, then the structure is considered for possible PCIC reconstruction. Cells in the structure with $h_i$ values of at most $0.8\Delta x$ and $\ell_i$ values of $1$ are marked for PCIC reconstruction.

\subsection{Ligament Tip Detection}
\label{Methods:tip_det}
As described in the preceding subsection, the CCL based algorithm identifies cells that are ligament-like and should be handled with PCIC. However, special care is needed when reconstructing ligament tips. If PCIC is naively applied to a ligament tip, a spuriously thin ligament can develop and rapidly grow as volume fractions are advected into neighboring cells. Thus, it is preferable to cap the end of a ligament with a non-cylindrical reconstruction. PPIC would, in theory, work well for this and will be explored in future work, but in this work a PLIC reconstruction is invoked at ligament tips. To detect such tips, a simple algorithm is used based on the asymmetry of liquid barycenters in a $3\times3\times3$ stencil. Only cells belonging to the center cell's CCL structure are considered in the following process. First, the cell centers are subtracted from their respective barycenters and the results are normalized by the mesh size. From these normalized barycenters, the phasic geometric center of the stencil is calculated (a stencil liquid barycenter). Then, a volume-weighted average of the squared Euclidean distances between the geometric center and each normalized liquid barycenter, named $\ell_{\mathrm{tip}_i}$, is found by
\begin{equation}
\begin{aligned}
\ell_{\mathrm{tip}_i}=\frac{1}{\sum_{j=1}^g\alpha_{ij}}\sum_{j=1}^g\left|\left|\mathbf{x}_{n_{ij}}-\bar{\mathbf{x}}_{n_{i}}\right|\right|^2\alpha_{ij},
\end{aligned}
\end{equation}
where $g$ is the number of cells in the stencil that belong to the center cell's CCL structure, $\mathbf{x}_{n_{i}}$ are the normalized barycenters, and $\bar{\mathbf{x}}_{n_{i}}$ is the stencil barycenter. If $\ell_{\mathrm{tip}_i}$ is less than or equal to $0.2$, a ligament tip is detected and the reconstruction method is reverted back to PLIC.

\subsection{Flow Solver Implementation}
\label{Methods:NGA2}
The PCIC reconstruction method has been implemented in the open-source Interface Reconstruction Library\footnote{Available at \url{https://github.com/cmfr-ae-illinois/interface-reconstruction-library}} (IRL)~\cite{Chiodi2,Evrard}, and simple translation and deformation test cases are run within IRL. The open-source flow solver NGA2\footnote{Available at \url{https://github.com/desjardi/NGA2}}~\cite{Desjardins} is used to run multiphase flow simulations with IRL. The ligament detection and ligament tip detection algorithms are implemented within NGA2, while IRL handles the actual moment calculations and interface reconstructions. The process proceeds as follows. First, the CCL ligament detection is carried out, which results in each cell being labeled as either containing no interface, a non-ligament-like interface, or a ligament-like interface. Then, ligament tip detection is done on the identified ligaments; cells containing ligament tips have their ligament-like label removed. The reconstruction process is then executed, with PLIC reconstructing the non-ligament-like cells and PCIC reconstructing the ligament-like cells. PLIC-Net~\cite{Cahaly} is used for the PLIC reconstruction in this work. It is also conceivable to add R2P reconstructions to this process by detecting thin films with the method described in \citet{Han3} and \citet{Han4} and using two-plane reconstructions where appropriate.

NGA2 uses a one-fluid formulation to solve the two-phase Navier–Stokes equations. The flow solver is second-order in time and space and discretely mass, momentum, and kinetic energy conserving away from the interface. However, at the interface, local discontinuities slightly degrade these conservation properties; although mass is still discretely conserved, momentum is approximately conserved and kinetic energy conservation is not maintained. The surface tension force is calculated using the Continuous Surface Force (CSF) model of \citet{Brackbill}, which uses an interface curvature for PLIC interfaces determined with parabolic surface fits~\cite{Jibben}. This curvature estimation provides a desirable balance between accuracy and computational cost~\cite{Han2}. Note that in this work, no surface tension modeling is done on PCIC interfaces, as discussed in Section~\ref{results: hit}. The development of such a model will be addressed in future work. Volume fraction field advection is performed with the unsplit geometric advection method of \mbox{\citet{Owkes}}.

\section{Test Cases for Calculating Volume Moments of Cylinders Clipping Polyhedra}
\label{Results1}

As mentioned in Section~\ref{Methods:NGA2}, the closed-form expressions derived in Section~\ref{sec:cylinder_case} have been implemented in the \texttt{C++} Interface Reconstruction Library~(IRL)~\cite{Chiodi2}. In this section, the speed, accuracy, and robustness of this implementation are tested on a wide variety of configurations with comparisons to analytical expressions or to estimations obtained using the same Adaptive Mesh Refinement (AMR) approach\footnote{The AMR approach has been adapted to consider quadratic cylinders instead of paraboloid surfaces. The code corresponding to the new implementation can be found in the file \href{https://github.com/cmfr-ae-illinois/interface-reconstruction-library/blob/quadratic_cutting/irl/generic_cutting/cylinder_intersection/cylinder_intersection_amr.tpp}{\texttt{irl/generic\_cutting/cylinder\_intersection/{\linebreak}cylinder\_intersection\_amr.tpp}}. 
} as introduced in~\cite{Evrard}. Although only straight circular cylinders are used in Sections~\ref{Methods2} and~\ref{Results2} for interface reconstruction, both elliptic and hyperbolic cylinders with a wide range of $\beta$ values are tested in this section. 

Section~\ref{sec:unit_translation} presents a canonical case of a translating unit cube being cut by a cylinder. Section~\ref{sec:paramsweep} includes a parameter sweep of several different test geometries at randomized locations and rotations and cut by randomized cylinders; Section~\ref{sec:paramnudge} repeats the parameter sweep exercise but for the ill-posed cases addressed in Section~\ref{sec:robustness}. Section~\ref{sec:timings} contains timing results for the cutting operations.

\subsection{Unit Cube Translating Along \texorpdfstring{$\mathbf{e}_z$}{ez}}\label{sec:unit_translation}
\begin{figure}[b!] \centering\vspace{-1em}
    \begin{tabular}{ccc}
    \subfloat[$k=1$]{\adjincludegraphics[width=0.32\textwidth,trim=0 0 0 0,clip=true]{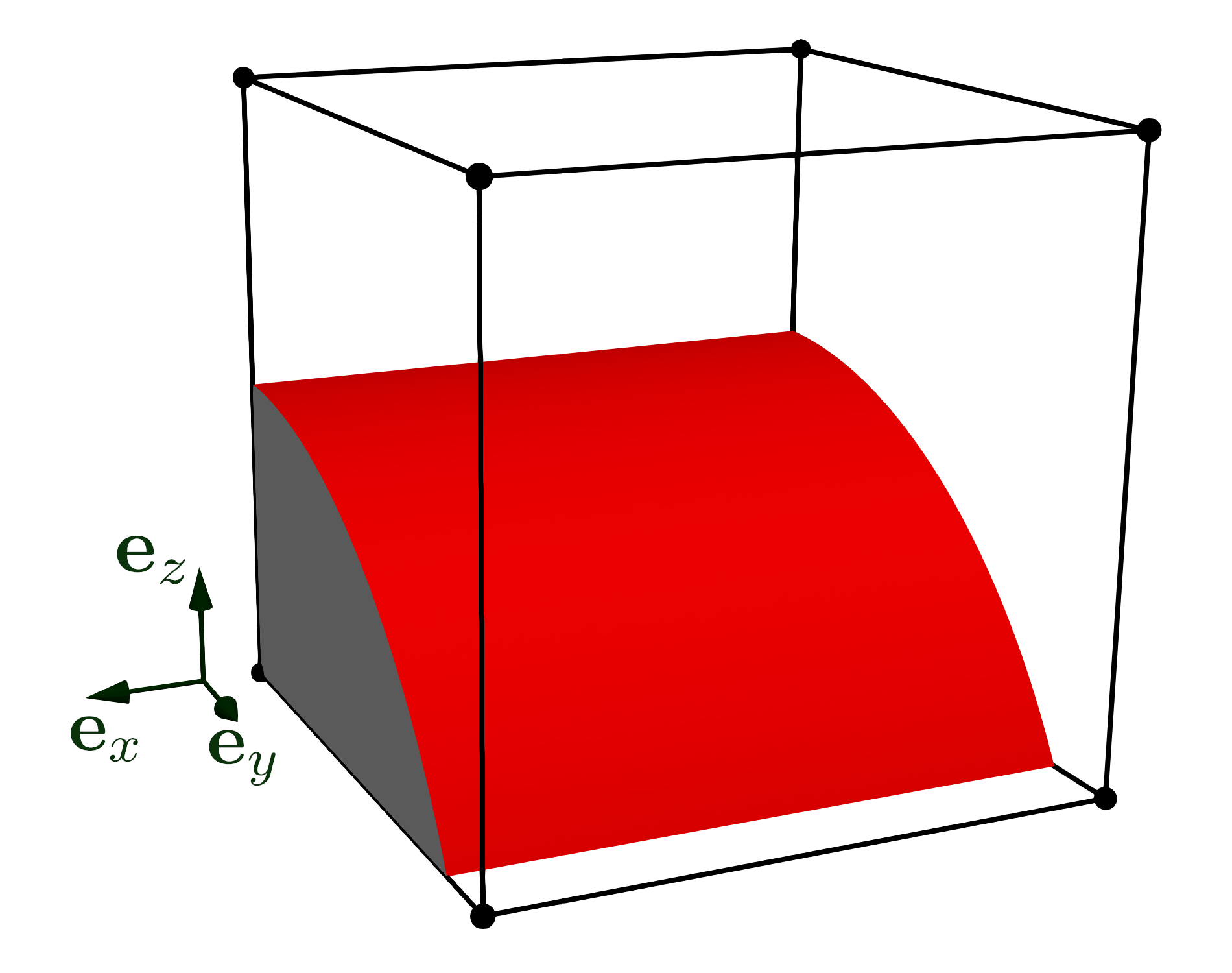}} \quad
    \subfloat[$k=2$]{\adjincludegraphics[width=0.32\textwidth,trim=0 0 0 0,clip=true]{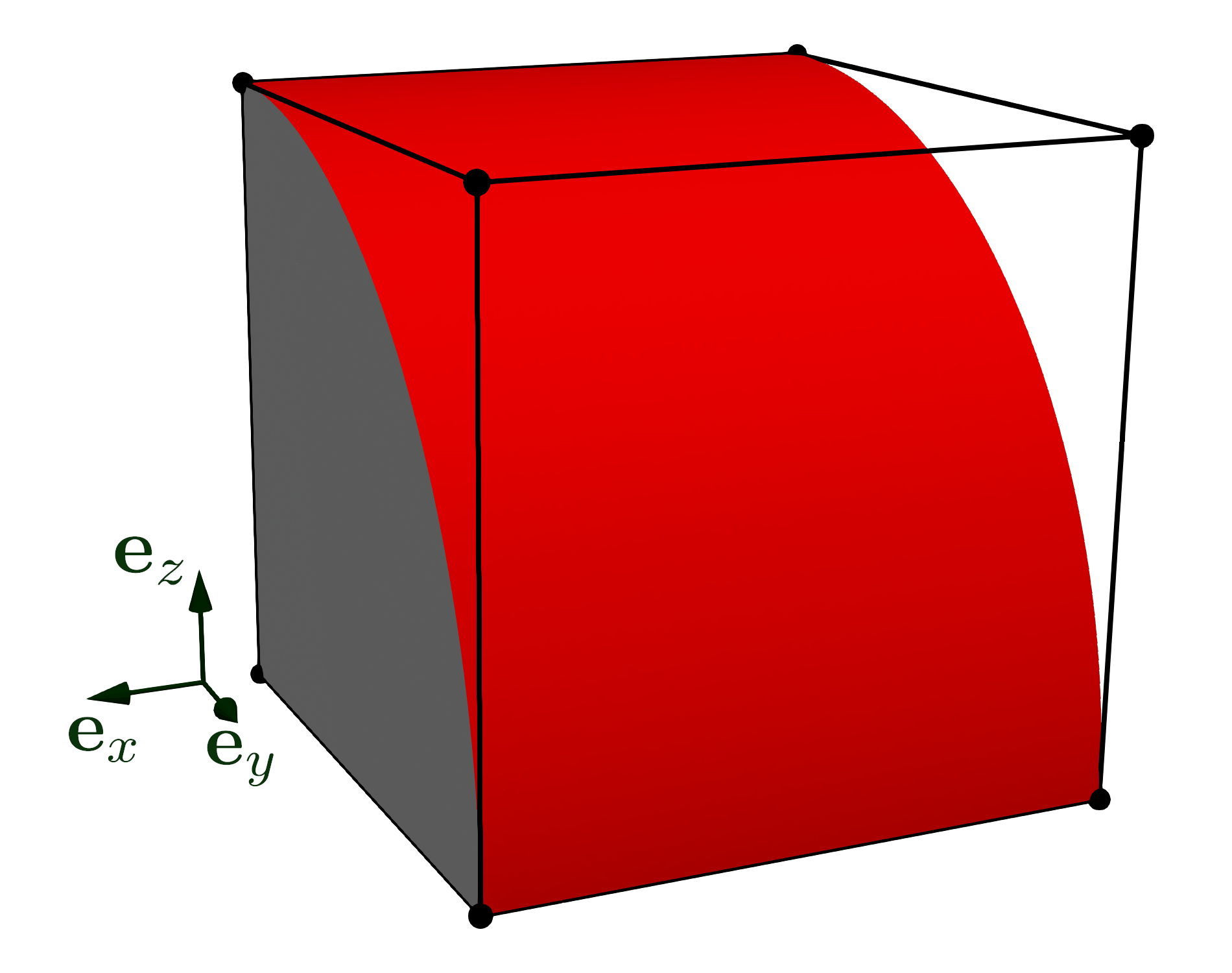}} \quad
    \subfloat[$k=3$]{\adjincludegraphics[width=0.32\textwidth,trim=0 0 0 0,clip=true]{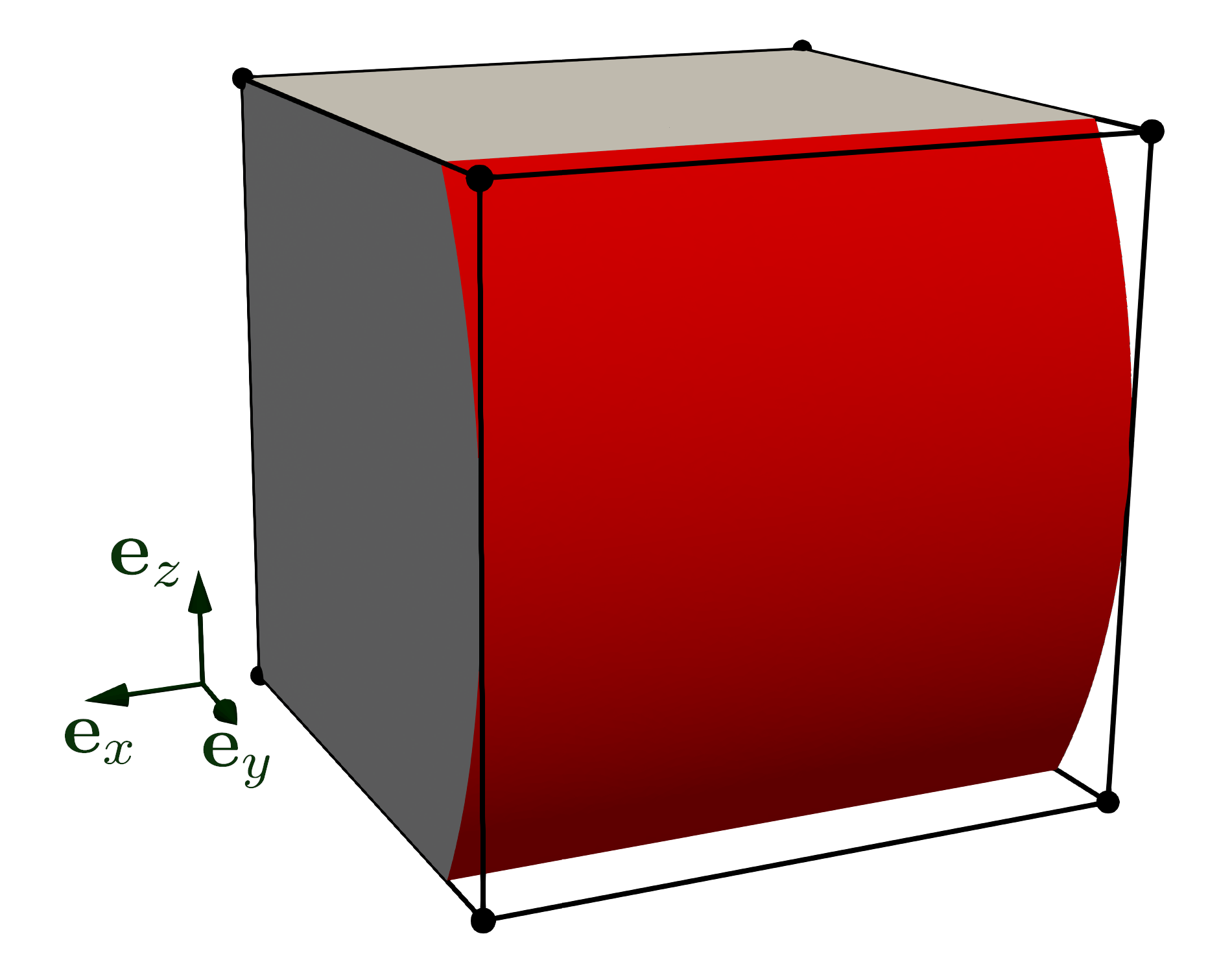}}
    \end{tabular}
    \caption{Unit cube centered at $\smash{\begin{bmatrix} 1/2 & 1/2 & 3/2 - k/2 \end{bmatrix}^\intercal}$ clipped by the elliptic cylinder parametrically defined as $y^2 + z^2 = 1$.}
    \label{fig:testcube}
\end{figure}

This first test considers the translation of a unit cube along the $\mathbf{e}_z$ axis, clipped by the canonical unit quadratic cylinder defined as in Eqs.~\eqref{eq:defS} and~\eqref{def:phi} with $\cylradius = \beta = 1$, as illustrated in Fig.~\ref{fig:testcube}. The unit cube is centered at the location $\smash{\mathbf{x}_c=\begin{bmatrix}1/2 & 1/2 & 3/2 - k/2\end{bmatrix}^\intercal}$, wherein the parameter $k$ is uniformly sampled on $[0,3]$ with the spacing $\Delta k = 10^{-3}$. For each value of $k$, the moments estimated with the expressions derived in Section~\ref{sec:cylinder_case} are compared to their exact expected value which, for this canonical test case, can be expressed as an analytic function of $k$.

\begin{figure}[t!] \centering\vspace{-2em}
    \includegraphics{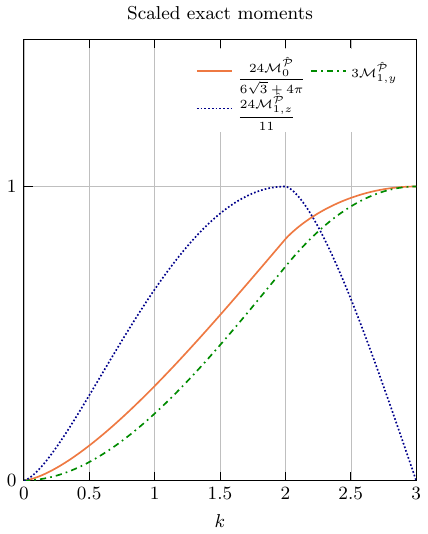}\includegraphics{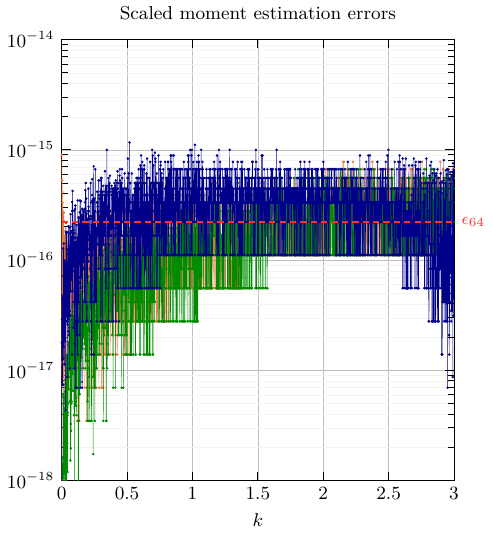}
    \caption{Moments of the unit cube centered at $\smash{\begin{bmatrix} 1/2 & 1/2 & 3/2 - k/2 \end{bmatrix}^\intercal}$ clipped by the elliptic cylinder defined as $y^2 + z^2 = 1$ (left); Errors in the estimation of the moments in $64$-bit floating-point format (right). The geometric moments and the corresponding estimation errors are scaled by their maximum value. The moments are computed using the closed-form expressions derived in~Section~\ref{sec:cylinder_case}. The $64$-bit machine epsilon $\epsilon_{64} = 2^{-52}$ is shown as the dashed red line.}\vspace{-1em}
    \label{fig:resultscube}
\end{figure}

Fig.~\ref{fig:resultscube} shows the exact moments of the clipped polyhedron $\hat{\mathcal{P}}$ (left graph) and the errors associated with their estimation using the proposed closed-form expressions (right graph), scaled by the maximum value of each exact moment. The estimation errors of all moment components are found to stay within an order of magnitude of $\epsilon_{64}$, the relative precision in 64-bit floating-point arithmetics.

\subsection{Parameter Sweep for Several Geometries} \label{sec:paramsweep}
In this second test, a selection of convex polyhedra (a regular tetrahedron, a cube, and a regular dodecahedron) and a nonconvex polyhedron (a ``hollow'' cube) are considered; they are constructed so that their volume is equal to unity. These polyhedra are centered at $\smash{\mathbf{x}_c = \begin{bmatrix}x_{c} & y_{c} & z_{c}\end{bmatrix}^\intercal}$ located within the cube $\smash{\begin{bmatrix}-1/2,1/2\end{bmatrix}^3}$ and rotated by $\smash{\{\theta_x, \theta_y, \theta_z\} \in [-\pi,\pi]^3}$ about the three main axes, respectively.
Each polyhedron is then intersected with a cylinder defined as in Eqs.~\eqref{eq:defS} and~\eqref{def:phi} with coefficients $\cylradius \in [0, 6/5]$ and $\beta \in [-10, 10]$. Finally, the geometric moments of the clipped polyhedron estimated using the expressions derived in Section~\ref{sec:cylinder_case} are compared with those obtained through AMR.

\newcolumntype{R}{>{$}r<{$}} %

Two series of sampling of this parameter space are conducted. First, all parameters are randomly sampled uniformly in their respective intervals. Then, a graded sweep is conducted wherein the parameters are chosen from a hand-selected discrete set of values that favors ambiguous and/or ill-posed intersection configurations\footnote{The graded sampling of the parameter space is chosen to create ambiguous configurations for the cube. In this graded sweep, the geometries other than the cube are not scaled to reach a unit volume so as to improve the chances of generating ambiguous/ill-posed intersection configurations.}. This discrete sampling is given as
\begin{align}
    \mathbf{x}_{c} &\in \left\{-\tfrac{1}{2},-\tfrac{1}{4},0,\tfrac{1}{4},\tfrac{1}{2}\right\}^3 \nonumber \\
    (\theta_x, \theta_y, \theta_z) &\in \{-\pi,-\tfrac{\pi}{2},0,\tfrac{\pi}{2},\pi\}^3 \nonumber \\
    \beta &\in \{\tfrac{9}{10},\ 1,\ \tfrac{16}{9},\ 2,\ \tfrac{9}{4},\ 4, -\tfrac{3}{4},\ -1,\ -\tfrac{5}{4}\} \nonumber \\
    \cylradius & \in \{\tfrac{1}{4}, \tfrac{1}{2}, \tfrac{1}{\sqrt{2}}, \tfrac{3}{4}, 1\} \nonumber
\end{align}
In total, $8$ million tests are conducted with random sampling of the parameter space (equally distributed between all geometries, with one half of the tests considering elliptic cylinders and the other half considering hyperbolic cylinders) and $5^3\times5^3\times9\times5 = 703,125$ distinct realizations are considered for each geometry in the graded parameter sweep. 

\begin{figure}[tbhp] \centering
    \begin{tabular}{c}
        \subfloat[\scriptsize Tetrahedron\label{fig:tet}]{\begin{minipage}{0.32\textwidth}\centering\adjincludegraphics[width=3.8cm,trim=0 0 0 0,clip=true]{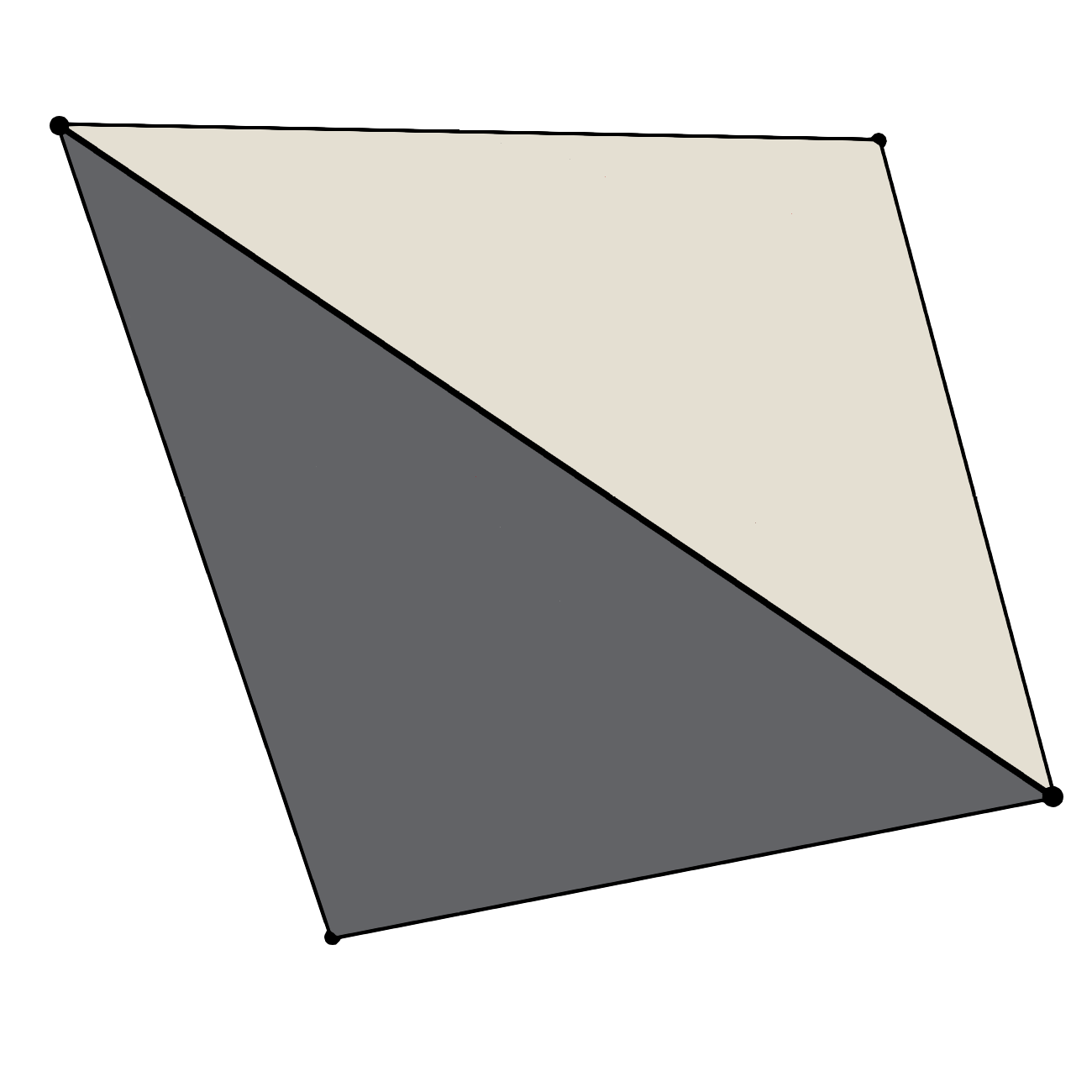}\end{minipage}}
        \subfloat[\scriptsize Clipped tetrahedron]{\begin{minipage}{0.32\textwidth}\centering\adjincludegraphics[width=3.8cm,trim=0 0 0 0,clip=true]{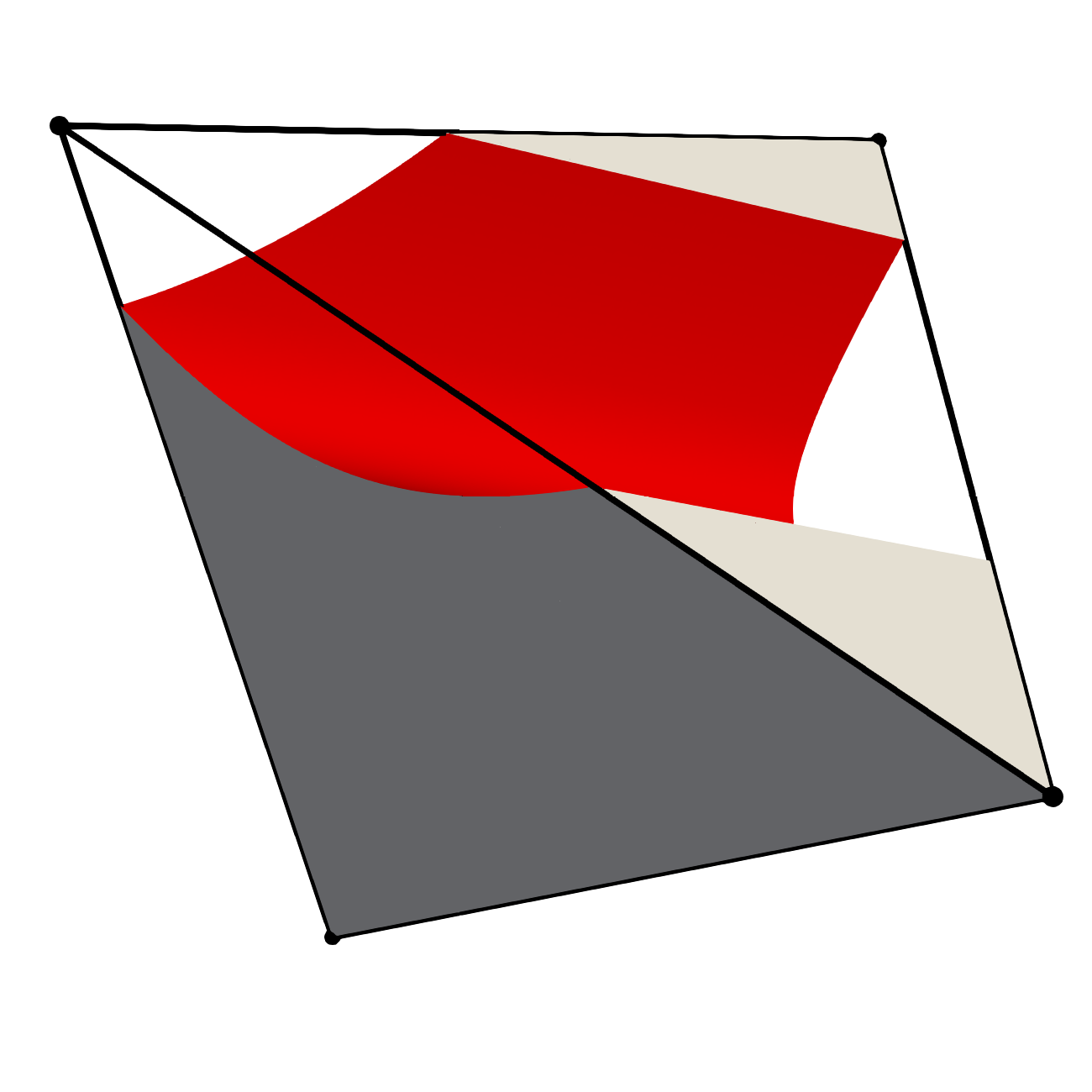}\end{minipage}}
        \subfloat[\scriptsize Clipped tetrahedron (AMR ref.)]{\begin{minipage}{0.35\textwidth}\centering\adjincludegraphics[width=3.8cm,trim=0 0 0 0,clip=true]{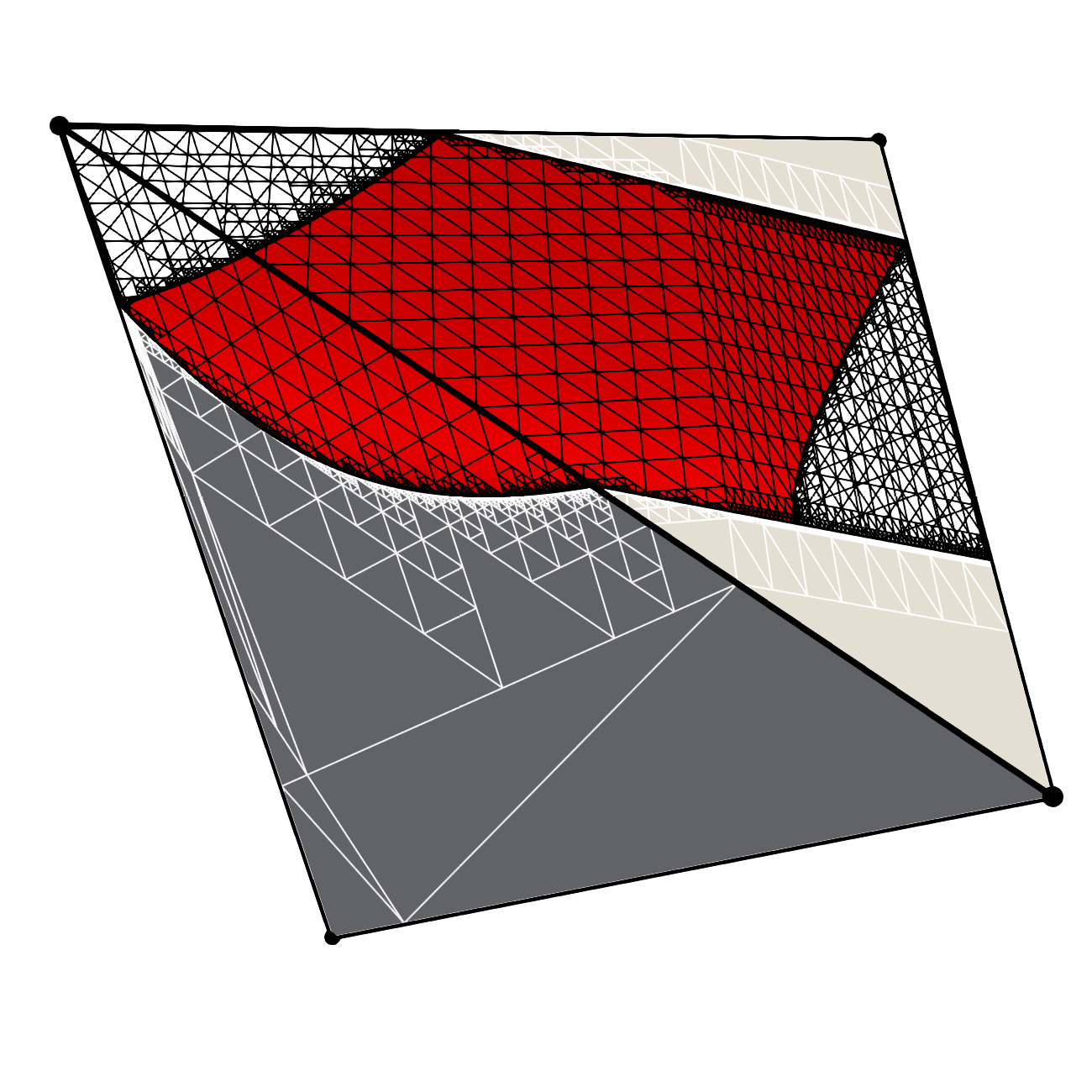}\end{minipage}}\\
        \subfloat[\scriptsize Cube\label{fig:cube}]{\begin{minipage}{0.32\textwidth}\centering\adjincludegraphics[width=3.8cm,trim=0 0 0 0,clip=true]{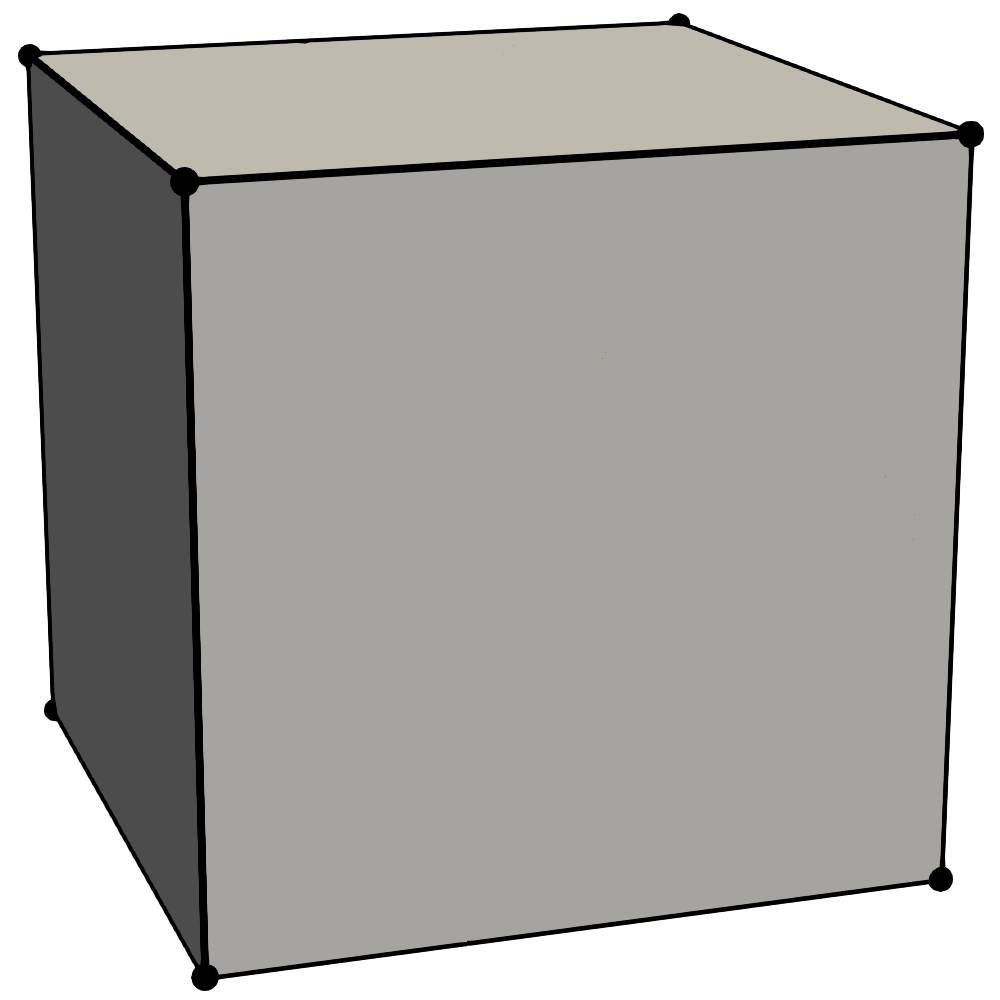}\end{minipage}}
        \subfloat[\scriptsize Clipped cube]{\begin{minipage}{0.32\textwidth}\centering\adjincludegraphics[width=3.8cm,trim=0 0 0 0,clip=true]{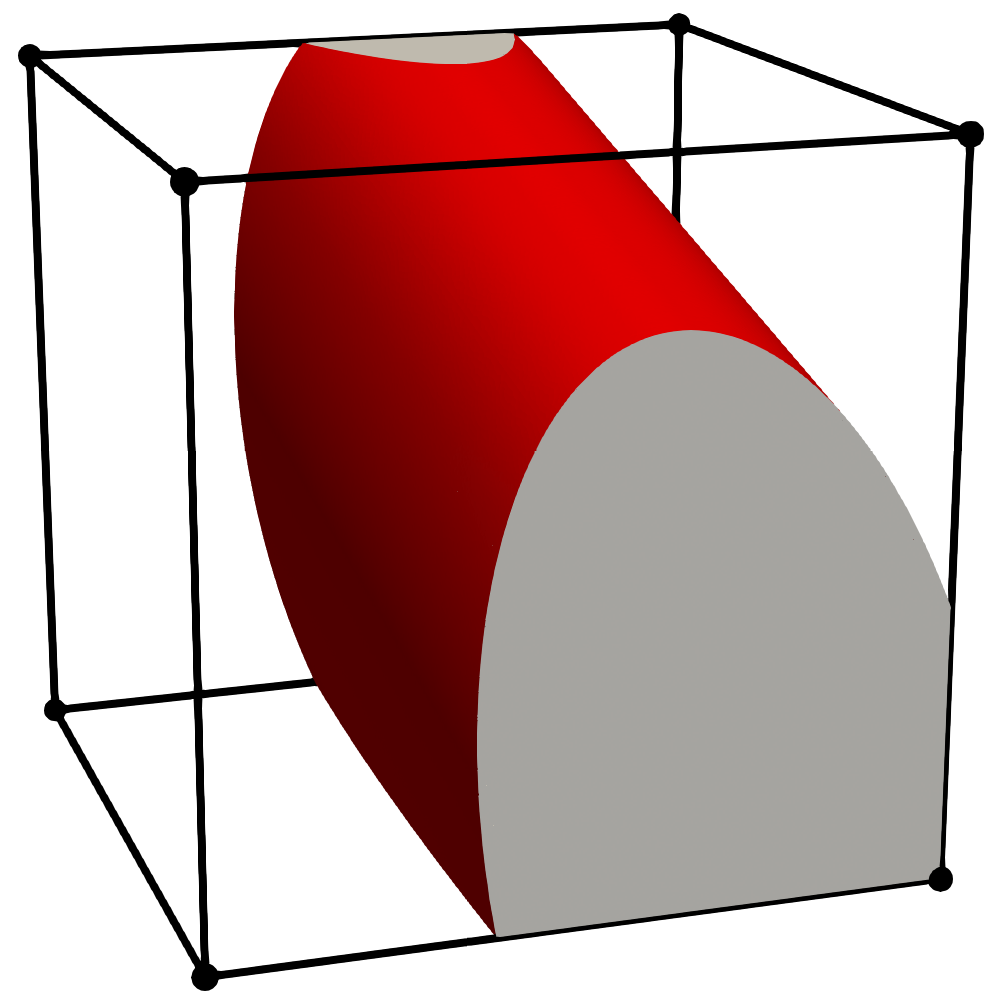}\end{minipage}}
        \subfloat[\scriptsize Clipped cube (AMR ref.)]{\begin{minipage}{0.35\textwidth}\centering\adjincludegraphics[width=3.8cm,trim=0 0 0 0,clip=true]{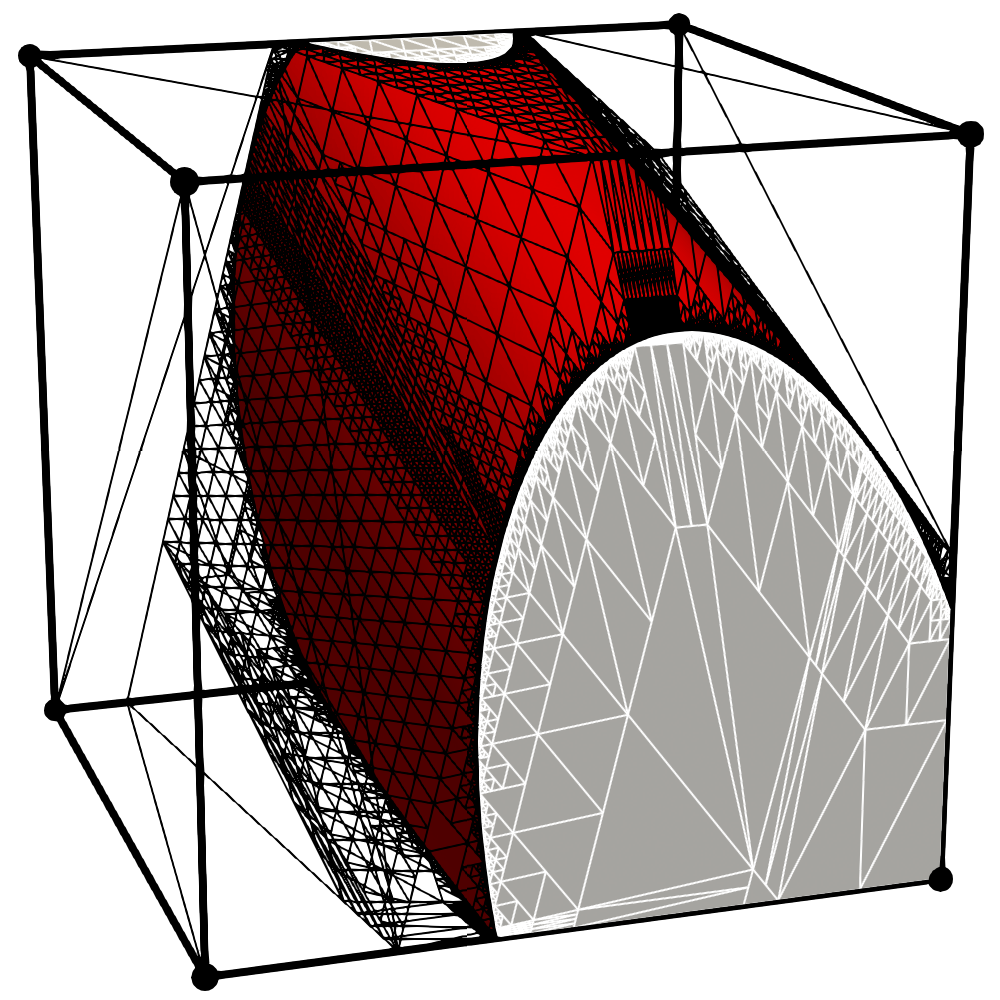}\end{minipage}}\\
        \subfloat[\scriptsize Dodecahedron\label{fig:dodeca}]{\begin{minipage}{0.32\textwidth}\centering\adjincludegraphics[width=3.8cm,trim=0 0 0 0,clip=true]{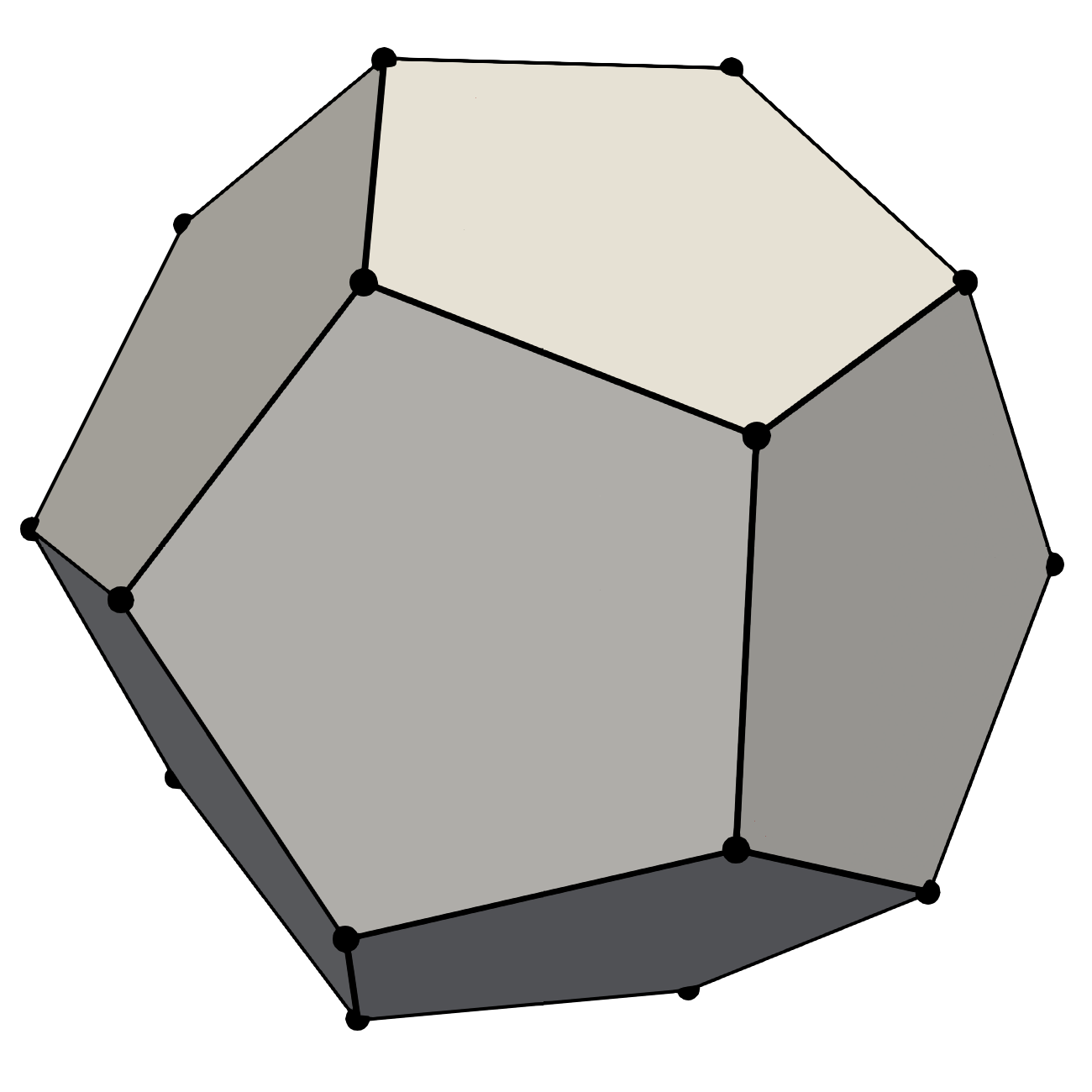}\end{minipage}}
        \subfloat[\scriptsize Clipped dodecahedron]{\begin{minipage}{0.32\textwidth}\centering\adjincludegraphics[width=3.8cm,trim=0 0 0 0,clip=true]{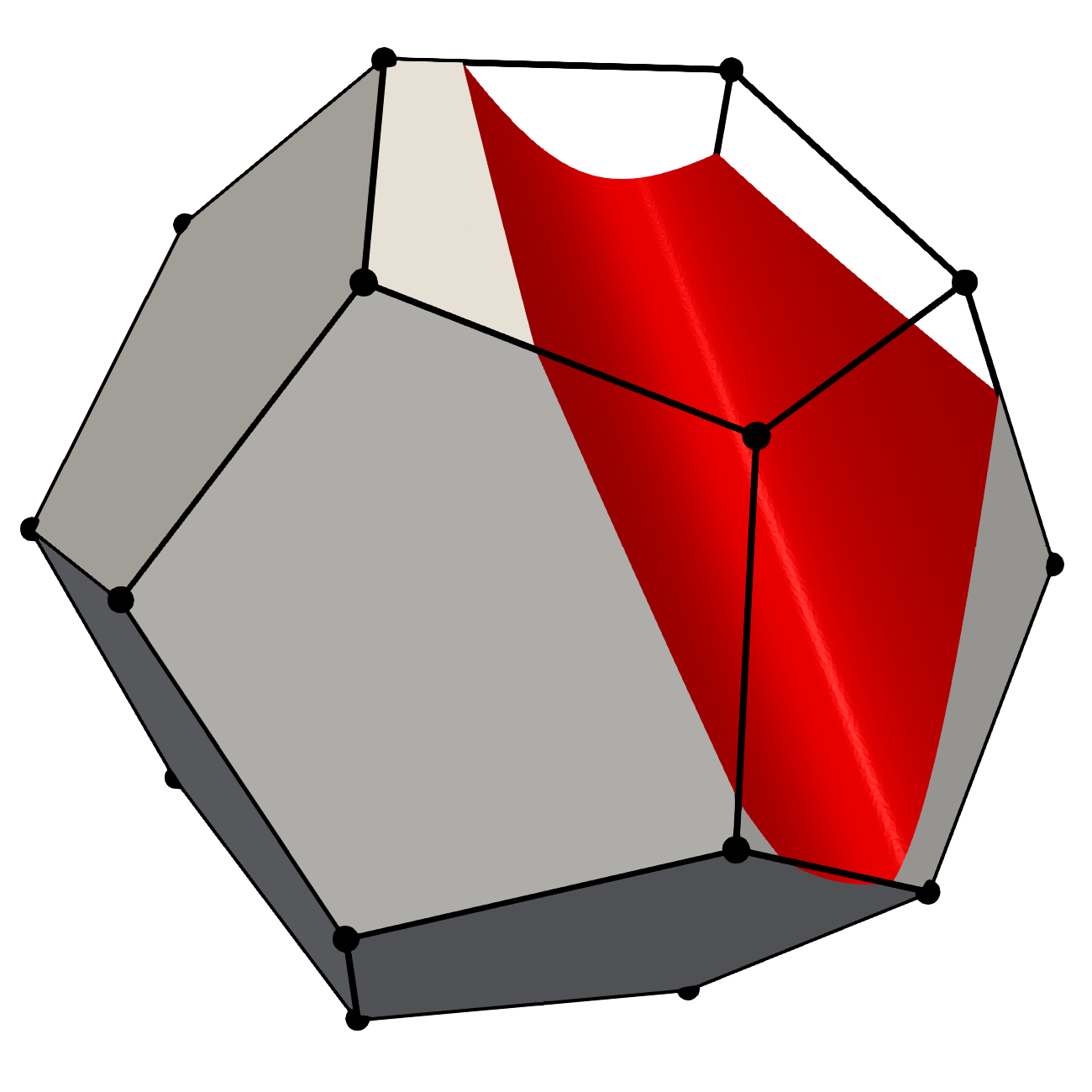}\end{minipage}}
        \subfloat[\scriptsize Clipped dodecahedron (AMR ref.)]{\begin{minipage}{0.35\textwidth}\centering\adjincludegraphics[width=3.8cm,trim=0 0 0 0,clip=true]{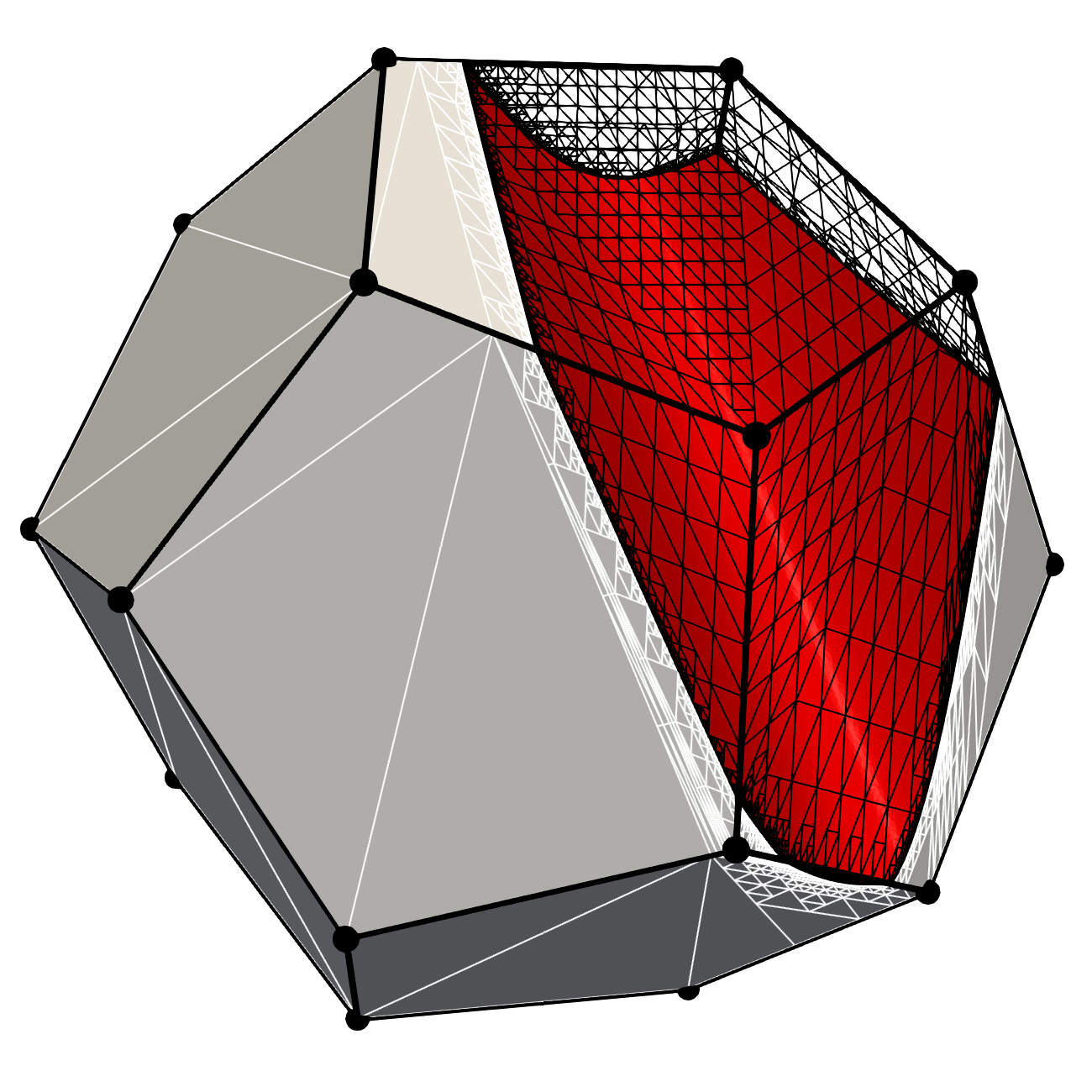}\end{minipage}}\\
        \subfloat[\scriptsize Hollow cube\label{fig:cubehole}]{\begin{minipage}{0.32\textwidth}\centering\adjincludegraphics[width=3.8cm,trim=0 0 0 0,clip=true]{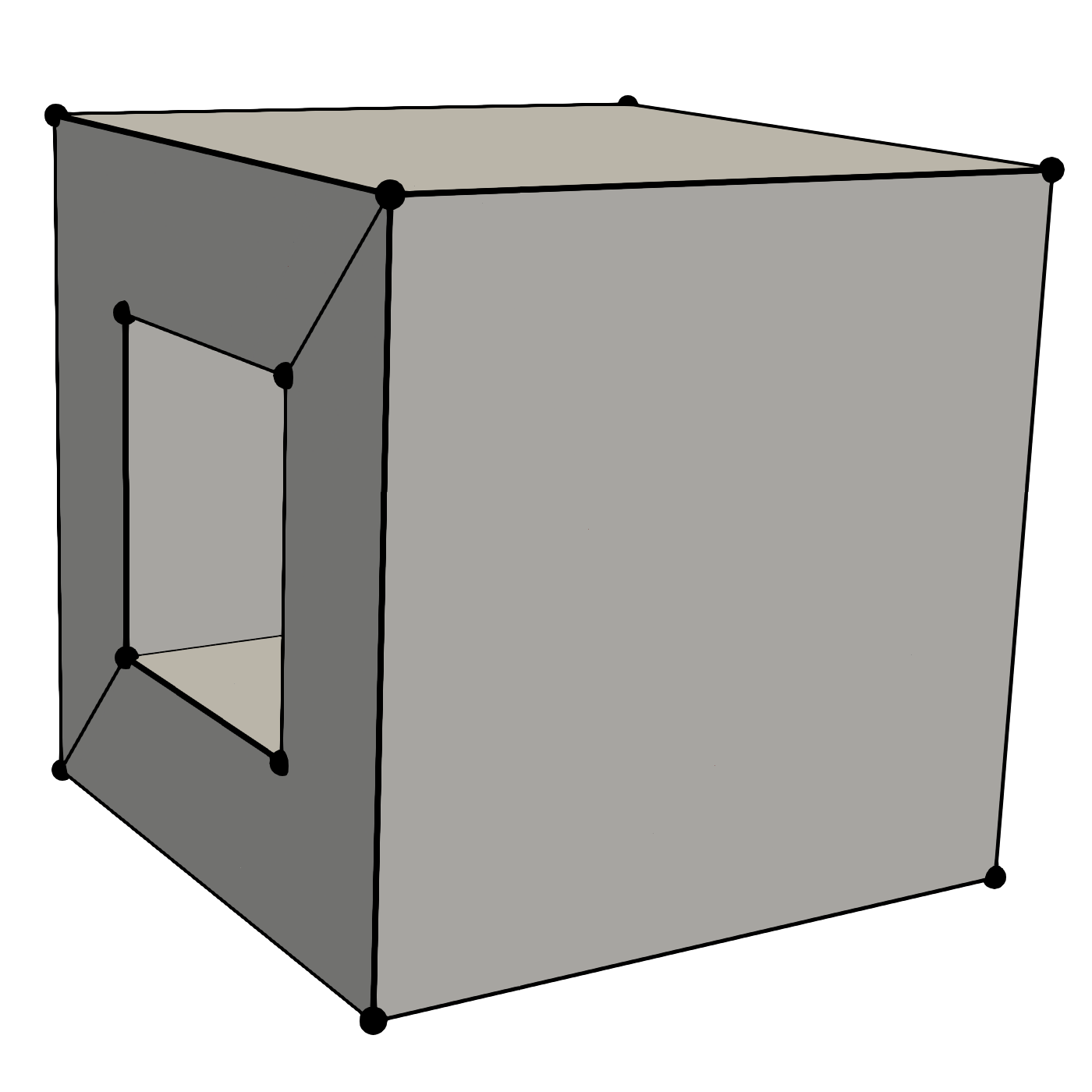}\end{minipage}}
        \subfloat[\scriptsize Clipped hollow cube]{\begin{minipage}{0.32\textwidth}\centering\adjincludegraphics[width=3.8cm,trim=0 0 0 0,clip=true]{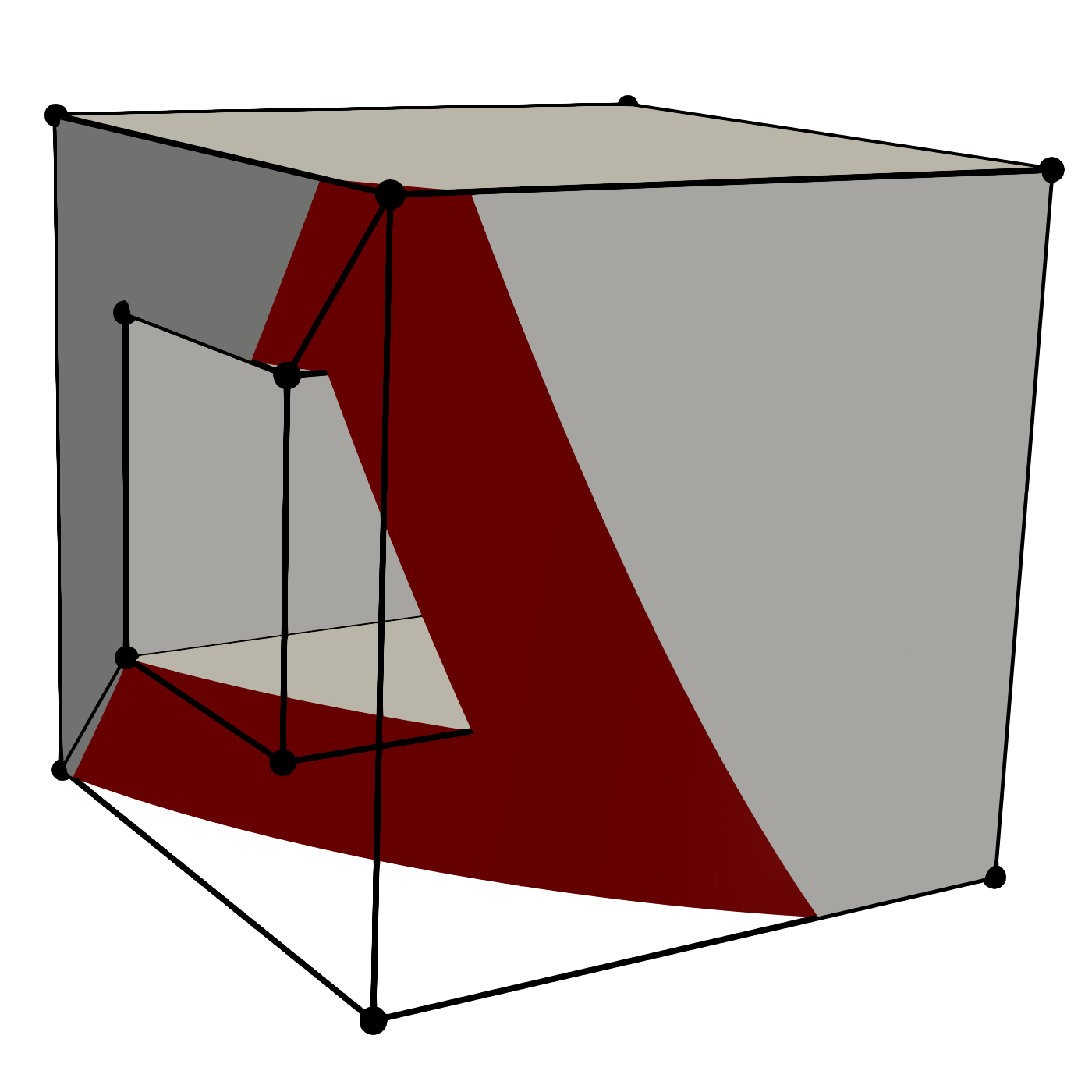}\end{minipage}}
        \subfloat[\scriptsize Clipped hollow cube (AMR ref.)]{\begin{minipage}{0.35\textwidth}\centering\adjincludegraphics[width=3.8cm,trim=0 0 0 0,clip=true]{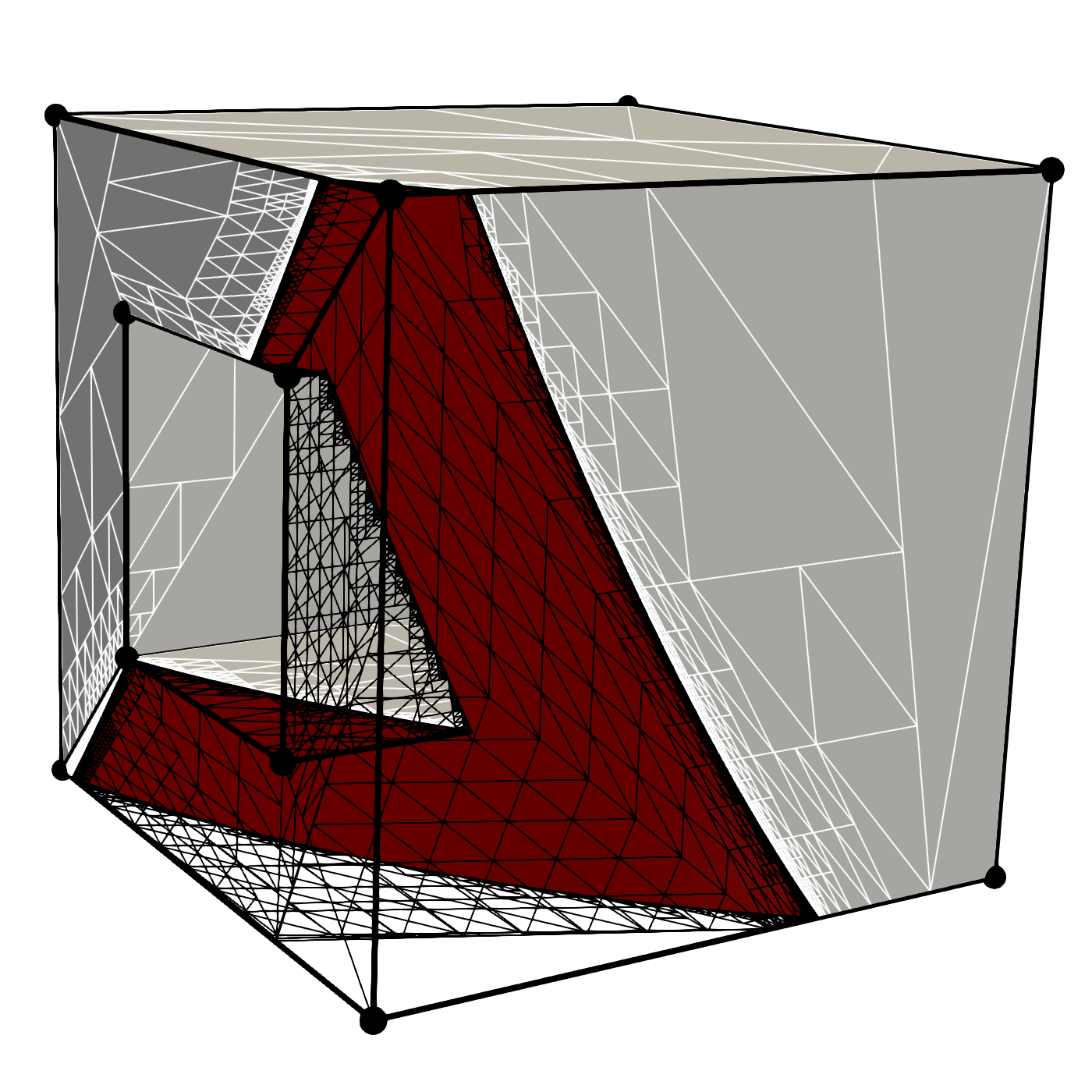}\end{minipage}}
    \end{tabular}
    \caption{Examples of random intersections of the polyhedra under consideration with a quadratic cylinder. The left column shows the original polyhedra; the central column shows them clipped by a cylinder; the right column shows the edges of the AMR used to calculate the reference moments.}
    \label{fig:testcases}
\end{figure}

\begin{figure}[tbhp] \centering
    \begin{tabular}{c}
        \subfloat[]{\begin{minipage}{0.24\textwidth}\centering\adjincludegraphics[width=\textwidth,trim=0 0 0 0,clip=true]{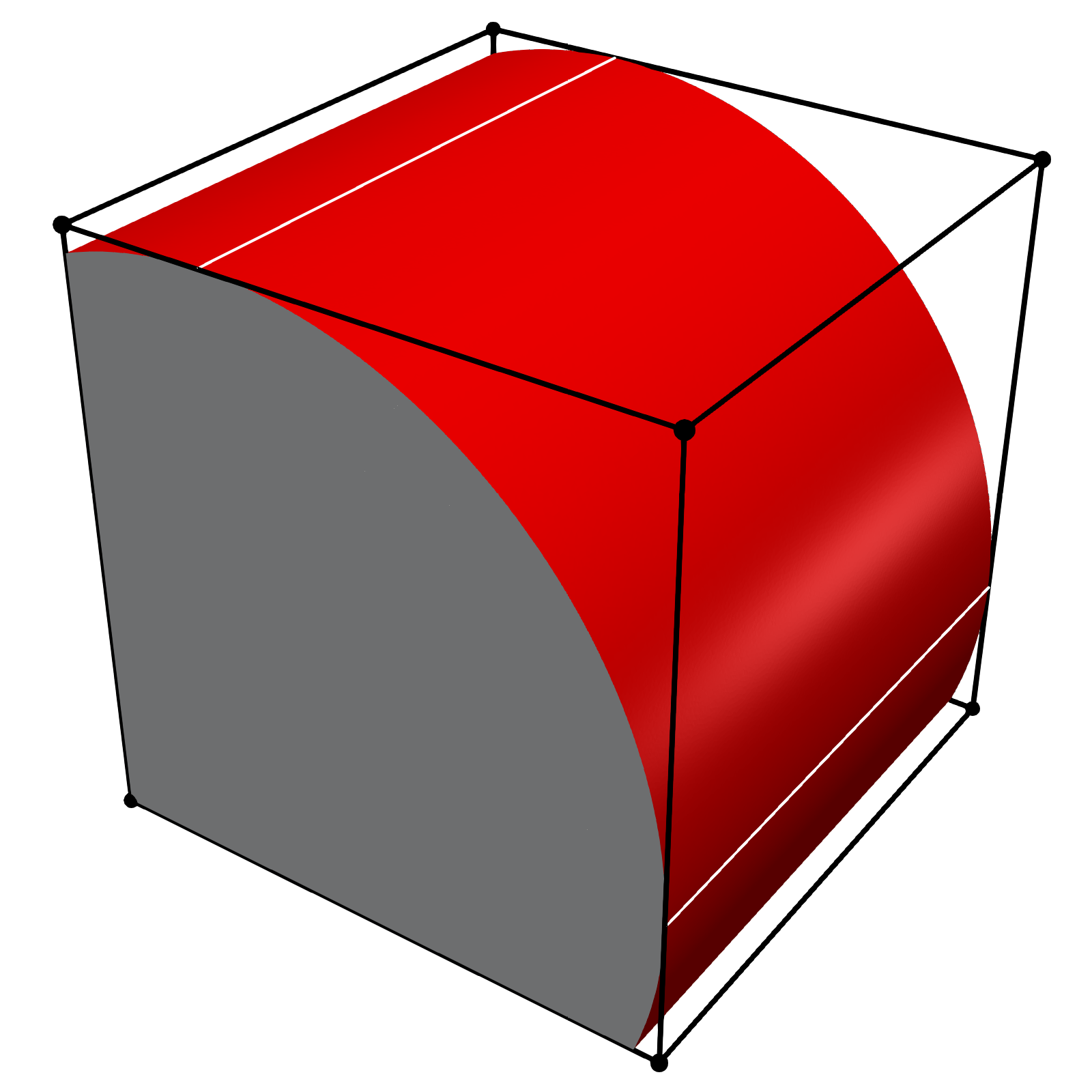}\end{minipage}}
        \subfloat[]{\begin{minipage}{0.24\textwidth}\centering\adjincludegraphics[width=\textwidth,trim=0 0 0 0,clip=true]{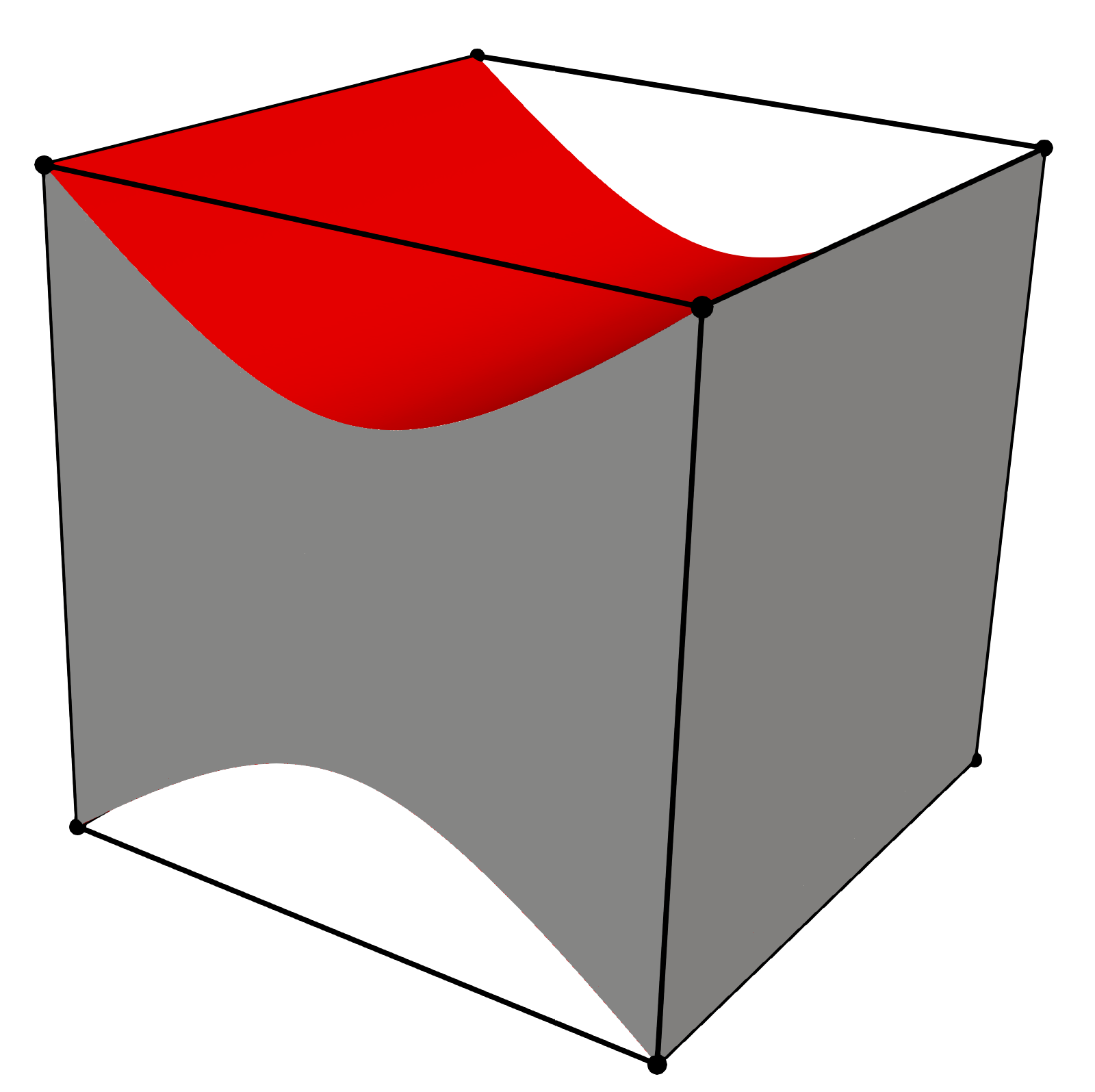}\end{minipage}}
        \subfloat[]{\begin{minipage}{0.24\textwidth}\centering\adjincludegraphics[width=\textwidth,trim=0 0 0 0,clip=true]{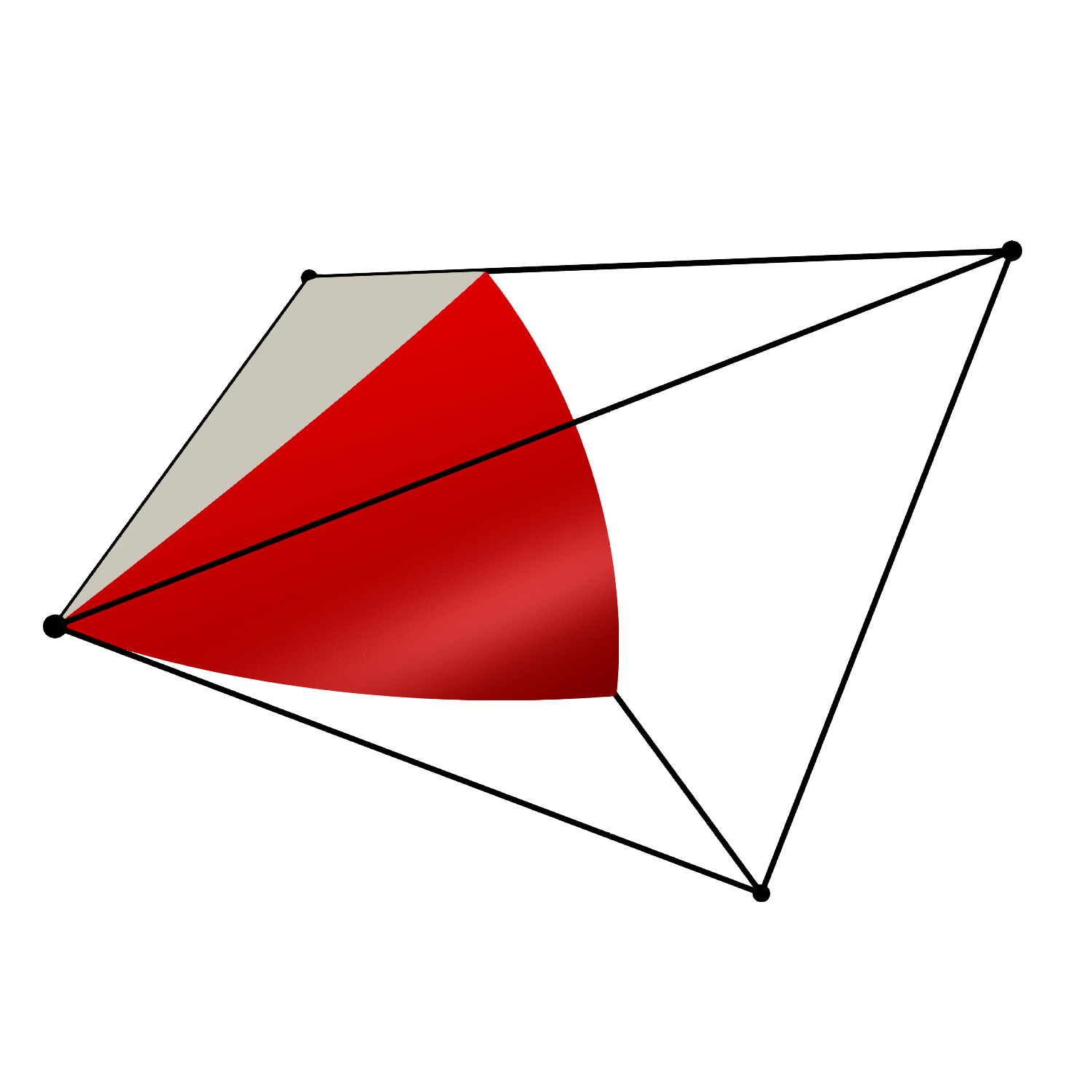}\end{minipage}}
        \subfloat[]{\begin{minipage}{0.24\textwidth}\centering\adjincludegraphics[width=\textwidth,trim=0 0 0 0,clip=true]{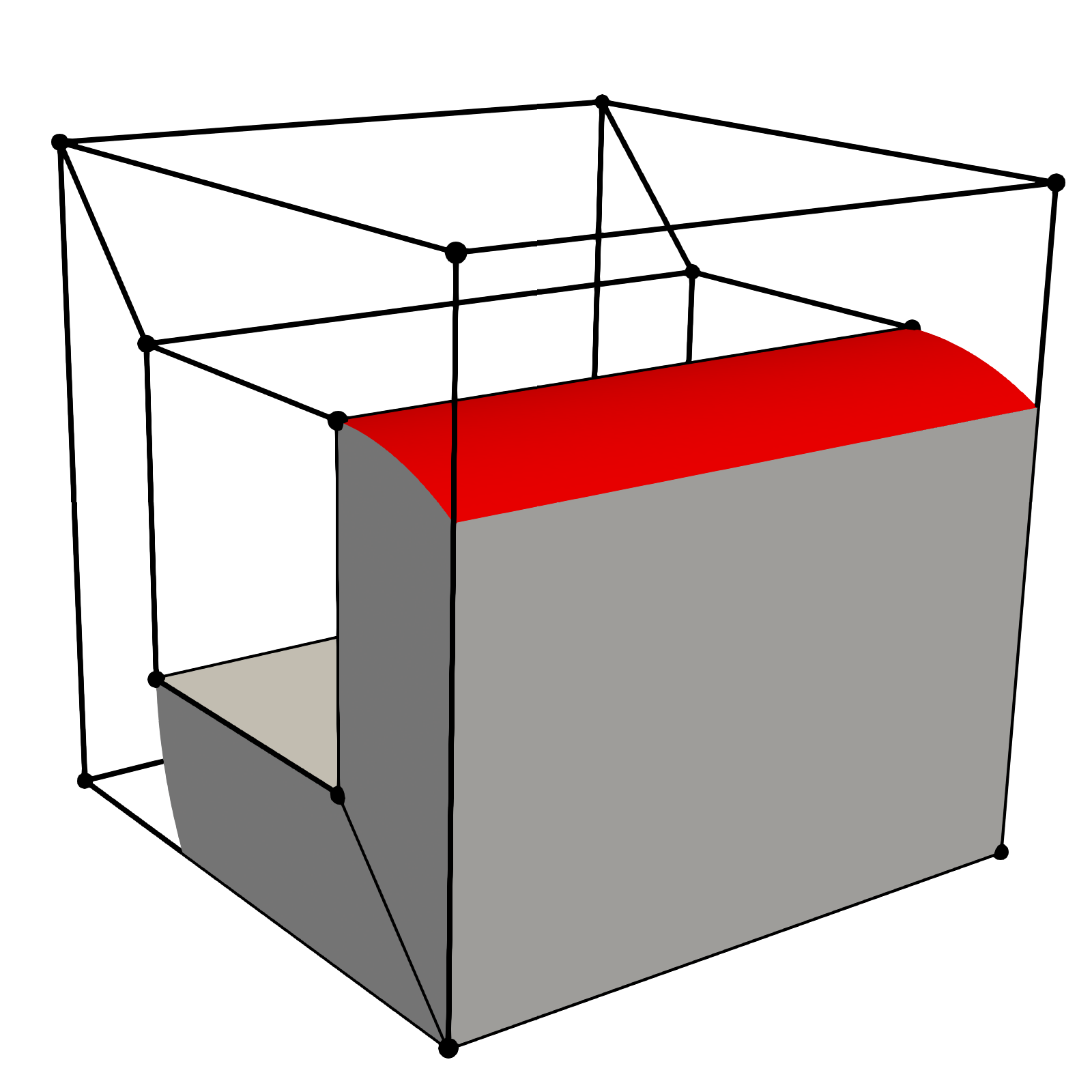}\end{minipage}}
    \end{tabular}
    \caption{Examples of intersections with degenerate cylinders and/or yielding ambiguous topologies, raised during the graded parameter sweep.}
    \label{fig:degeneratecases}
\end{figure}

Examples of random intersection configurations and of the corresponding AMR approximation are shown in Fig.~\ref{fig:testcases}. Figure~\ref{fig:degeneratecases} shows examples of degenerate intersection configurations encountered in the graded parameter sweep.
The maximum and average moment estimation errors are shown in Table~\ref{tab:randomsweep} for the random parameter sweep and in Table~\ref{tab:regularsweep} for the graded parameter sweep. 

\begin{table}[tbhp]
    \renewcommand\arraystretch{1.4}
    \scriptsize
    \captionsetup{position=top} 
    \caption{Random parameter sweep results. For each geometry, the following are provided: the number of tests conducted, the number of recursive levels used for the AMR reference moment calculation, the average and maximum errors in the estimation of the zeroth and first moments.}\label{tab:randomsweep}
    \begin{center}
    Elliptic cylinder results\\\vspace{0.25em}
        \begin{tabular}{|r|R|R|R|R|R|R|} \hline
            \multicolumn{1}{|c|}{\multirow{2}{*}{\quad\;\textbf{Geometry}\;\quad}} & \multicolumn{1}{c|}{\multirow{2}{*}{\shortstack{\textbf{Number}\\\textbf{of tests}}}} & \multicolumn{1}{c|}{\multirow{2}{*}{\shortstack{\textbf{AMR}\\\textbf{levels}}}} & \multicolumn{2}{c|}{\textbf{Zeroth moment error}} &  \multicolumn{2}{c|}{\textbf{First moments error}} \\ \cline{4-7}
            & & & \multicolumn{1}{c|}{Average} & \multicolumn{1}{c|}{Maximum} &\multicolumn{1}{c|}{Average} & \multicolumn{1}{c|}{Maximum} \\ \hline
        Tetrahedron &  10^6 & 15 & 1.3 \times 10^{-16} & 3.6 \times 10^{-15} & 4.0 \times 10^{-17} & 2.3 \times 10^{-15} \\ \hline
    Cube & 10^6 & 15 &  1.0 \times 10^{-16} & 1.4 \times 10^{-15} & 2.9 \times 10^{-17} & 8.0 \times 10^{-16} \\ \hline
    Dodecahedron & 10^6 & 15 &  1.4 \times 10^{-16} & 1.2 \times 10^{-15} & 3.3 \times 10^{-17} & 7.8 \times 10^{-16} \\ \hline
    Hollow cube & 10^6 & 15 & 1.2 \times 10^{-16} & 2.2 \times 10^{-15} & 3.3 \times 10^{-17} & 1.6 \times 10^{-15} \\ \hline
    \end{tabular}
    \vspace{1.5em}
    
    Hyperbolic cylinder results\\\vspace{0.25em}
        \begin{tabular}{|r|R|R|R|R|R|R|} \hline
            \multicolumn{1}{|c|}{\multirow{2}{*}{\quad\;\textbf{Geometry}\;\quad}} & \multicolumn{1}{c|}{\multirow{2}{*}{\shortstack{\textbf{Number}\\\textbf{of tests}}}} & \multicolumn{1}{c|}{\multirow{2}{*}{\shortstack{\textbf{AMR}\\\textbf{levels}}}} & \multicolumn{2}{c|}{\textbf{Zeroth moment error}} &  \multicolumn{2}{c|}{\textbf{First moments error}} \\ \cline{4-7}
            & & & \multicolumn{1}{c|}{Average} & \multicolumn{1}{c|}{Maximum} &\multicolumn{1}{c|}{Average} & \multicolumn{1}{c|}{Maximum} \\ \hline
    Tetrahedron & 10^6 & 15 & 2.4 \times 10^{-16} & 9.0 \times 10^{-15} & 9.6 \times 10^{-17} & 2.4 \times 10^{-14} \\ \hline
    Cube & 10^6 & 15 &  1.8 \times 10^{-16} & 7.1 \times 10^{-14} & 5.9 \times 10^{-17} & 5.0 \times 10^{-14} \\ \hline
    Dodecahedron & 10^6 & 15 &  2.4 \times 10^{-16} & 6.2 \times 10^{-14} & 7.2 \times 10^{-17} & 1.9 \times 10^{-14} \\ \hline
    Hollow cube & 10^6 & 15 & 2.4 \times 10^{-16} & 3.3 \times 10^{-15} & 8.2 \times 10^{-17} & 1.1 \times 10^{-12} \\ \hline
    \end{tabular}
    \end{center}
\end{table}

\begin{table}[tbhp]
        \renewcommand\arraystretch{1.4}
        \scriptsize
        \captionsetup{position=top} 
        \caption{Graded parameter sweep results. For each geometry, the following are provided: the number of tests conducted, the number of recursive levels used for the AMR reference moment calculation, the average and maximum errors in the estimation of the zeroth and first moments.}\label{tab:regularsweep}
        \begin{center}
        \begin{tabular}{|r|R|R|R|R|R|R|} \hline
        \multicolumn{1}{|c|}{\multirow{2}{*}{\quad\;\textbf{Geometry}\;\quad}} & \multicolumn{1}{c|}{\multirow{2}{*}{\shortstack{\textbf{Number}\\\textbf{of tests}}}} & \multicolumn{1}{c|}{\multirow{2}{*}{\shortstack{\textbf{AMR}\\\textbf{levels}}}} & \multicolumn{2}{c|}{\textbf{Zeroth moment error}} &  \multicolumn{2}{c|}{\textbf{First moments error}} \\ \cline{4-7}
        & & & \multicolumn{1}{c|}{Average} & \multicolumn{1}{c|}{Maximum} &\multicolumn{1}{c|}{Average} & \multicolumn{1}{c|}{Maximum} \\ \hline
        Tetrahedron &  703,125 & 15 & 3.2 \times 10^{-17} & 3.6 \times 10^{-15} & 1.1 \times 10^{-17} & 2.6 \times 10^{-15} \\ \hline
        Cube & 703,125 & 15 & 1.7 \times 10^{-16} & 5.0 \times 10^{-15} & 6.0 \times 10^{-17} & 3.4 \times 10^{-15} \\ \hline
        Dodecahedron & 703,
        125 & 15 & 1.6 \times 10^{-16} & 1.3 \times 10^{-15} & 5.9 \times 10^{-17} & 7.8 \times 10^{-16} \\ \hline
        Hollow cube &  703,125 & 15 & 1.9 \times 10^{-16} & 4.8 \times 10^{-15} & 9.6 \times 10^{-17} & 3.6 \times 10^{-15} \\ \hline
        \end{tabular}
        \end{center}
\end{table}

In every case of the random and graded parameter sweeps, the average error for the estimation of the volume of the clipped polyhedron is on the order of $10^{-16}$, and for the estimation of the first moments is on the order of $10^{-17}$. The maximum moment estimation error is around $10^{-14}$, except for the case of the hollow cube intersected with a hyperbolic cylinder, where the maximum error reaches $10^{-12}$ for two cases in which $\cylradius$ is particularly small.

\subsection{Parameter Sweep with Nudge}\label{sec:paramnudge} 
Another way to specifically test for ambiguous intersection configuration and assess the robustness of the algorithm is to consider random intersection cases but \textit{a posteriori} enforce that one vertex of the polyhedron $\mathcal{P}$ lands exactly on the surface of the quadratic cylinder $\mathcal{S}$. This is done by translating the polyhedron so as to position one of its vertices at the coordinate $\smash{[0 \quad 0 \quad \cylradius]^\intercal}$. All of these configurations require at least one iteration of the nudging procedure mentioned in Section~\ref{sec:robustness} and described in~\cite{Evrard}. The results of this random parameter sweep are given in Table~\ref{tab:ambiguouscases}. They do not display any significant increase of error compared to the random parameter sweep presented in Section~\ref{sec:paramsweep}. 

\begin{table}[tbhp]
    \renewcommand\arraystretch{1.4}
    \scriptsize
    \captionsetup{position=top} 
    \caption{Random parameter sweep with one vertex of the polyhedron lying exactly on the quadratic cylinder. The following are provided: the number of tests conducted, the number of recursive levels used for the AMR reference moment calculation, the average and maximum errors in the estimation of the zeroth and first moments.}\label{tab:ambiguouscases}
    \begin{center}
        \begin{tabular}{|r|R|R|R|R|R|R|} \hline
            \multicolumn{1}{|c|}{\multirow{2}{*}{\shortstack{\quad\;\textbf{Geometry}\;\quad\\(One vertex on $\mathcal{S}$)}}}& \multicolumn{1}{c|}{\multirow{2}{*}{\shortstack{\textbf{Number}\\\textbf{of tests}}}} & \multicolumn{1}{c|}{\multirow{2}{*}{\shortstack{\textbf{AMR}\\\textbf{levels}}}} & \multicolumn{2}{c|}{\textbf{Zeroth moment error}} &  \multicolumn{2}{c|}{\textbf{First moments error}} \\ \cline{4-7}
            & & & \multicolumn{1}{c|}{Average} & \multicolumn{1}{c|}{Maximum} &\multicolumn{1}{c|}{Average} & \multicolumn{1}{c|}{Maximum} \\ \hline
        Tetrahedron &  10^6 & 15 & 1.1 \times 10^{-16} & 3.0 \times 10^{-15} & 8.5 \times 10^{-17} & 7.0 \times 10^{-15} \\ \hline
        Cube & 10^6 & 15 &  1.1 \times 10^{-16} & 2.4 \times 10^{-15} & 5.9 \times 10^{-17} & 8.1 \times 10^{-14} \\ \hline
        Dodecahedron & 10^6 & 15 &  1.5 \times 10^{-16} & 1.7 \times 10^{-15} & 7.6 \times 10^{-17} & 1.8 \times 10^{-15} \\ \hline
        Hollow cube & 10^6 & 15 & 1.6 \times 10^{-16} & 2.2 \times 10^{-15} & 9.0 \times 10^{-17} & 1.8 \times 10^{-14} \\ \hline
    \end{tabular}
    \end{center}
\end{table}

\subsection{Computational Cost} \label{sec:timings}
Finally, the estimation of the moments in the random parameter sweep configuration is timed using the \texttt{C++} standard function \texttt{chrono::high\_resolution\_clock::now()}. The \texttt{C++} implementation of the closed-form expression was compiled with the GNU 14.2 suite of compilers~\cite{GCC14}, on a workstation whose characteristics are summarized in Table~\ref{tab:workstation}.

\begin{table}[tbhp]
    \renewcommand\arraystretch{1.4}
    \scriptsize
    \captionsetup{position=top} 
    \caption{Characteristics of the workstation used for the timing results presented in Section~\ref{sec:timings}.}\label{tab:workstation}
    \begin{center}
        \begin{tabular}{|c|l|l|} \hline
            \multirow{9}{*}{\textbf{CPU}}   & vendor\_id & GenuineIntel \\ \cline{2-3}
            & CPU family & 6 \\ \cline{2-3}
            & Model & 94 \\ \cline{2-3}
            & Model name & Intel\textregistered \, Core\texttrademark \, i5-6600 CPU @ 3.30GHz \\ \cline{2-3}
            & Stepping & 3 \\ \cline{2-3}
            & Microcode & 0xf0 \\ \cline{2-3}
            & Min/max clock CPU frequency & 800 MHz -- 3.90 GHz \\ \cline{2-3}
            & CPU asserted frequency & 3.60 GHz \\ \cline{2-3}
            & Cache size & 6144 KB \\ \hline
            \multirow{3}{*}{\textbf{Compiler}}   & Suite & GNU \\ \cline{2-3}
            & Version & 14.2.0 \\ \cline{2-3}
            & Flags & -O3 -DNDEBUG -DNDEBUG\_PERF$^\ast$ \\ \hline
        \end{tabular}\\\vspace{1mm}
    $^\ast$ -DNDEBUG\_PERF is a flag that disables IRL-specific debugging assertions~\cite{Chiodi2,Evrard}.
    \end{center}
\end{table}

The measured timings for the clipping of a polyhedron by a quadratic cylinder are summarized in Table~\ref{tab:timings}, along with timings of paraboloid- and plane-clipping remeasured on the same workstation, for the sake of comparison.
\begin{table}[t!]
    \renewcommand\arraystretch{1.4}
    \scriptsize
    \captionsetup{position=top} 
    \caption{Timings for the random parameter sweeps presented in Sections~\ref{sec:paramsweep} and~\ref{sec:paramnudge}.}\label{tab:timings}
    \begin{center}
        \begin{tabular}{|r|R|R|R|} \hline
            \multicolumn{1}{|c|}{\multirow{3}{*}{\quad\;\textbf{Geometry}\;\quad}} & \multicolumn{1}{c|}{\multirow{3}{*}{\shortstack{\textbf{Number} \\ \textbf{of tests}}}} & \multicolumn{2}{c|}{\textbf{Average moment calculation time}} \\ \cline{3-4}
            & & \multicolumn{1}{c|}{Zeroth moment only} & \multicolumn{1}{c|}{Zeroth and first moments} \\
            & & \multicolumn{1}{c|}{$\mu$s/test} & \multicolumn{1}{c|}{$\mu$s/test} \\ \hline
    Tetrahedron & 5 \times 10^7 & 6.26 & 7.16 \\ \hline
    Cube & 5 \times 10^7 & 9.20 & 10.13 \\ \hline
    Cube (one vertex on $\mathcal{S}$) & 5 \times 10^7 & 419.05 & 540.48 \\ \hline
    Dodecahedron & 5 \times 10^7 & 16.24 & 17.27 \\ \hline
    Hollow cube & 5 \times 10^7 & 19.94 & 22.07 \\ \hline
     \multicolumn{4}{l}{\textit{Comparison to paraboloid and plane clipping:}} \\
    \hline
    Cube (paraboloid clipping)~\cite{Evrard} & 5 \times 10^{7} & 2.02 & 3.22 \\ \hline
    Cube (paraboloid clipping)~\cite{Evrard} & 5 \times 10^{7} & 285.93 & 651.91 \vspace{-0.25em}\\ 
    (one vertex on $\mathcal{S}$) & & & \\ \hline
    Cube (plane clipping)~\cite{Chiodi2} & 5 \times 10^{7} & 0.33 & 0.35 \\ \hline
    \end{tabular}
    \end{center}
\end{table}
Because the closed-form expressions derived in Section~\ref{sec:cylinder_case} are only valid for polyhedra in the top halfspace, and the generalization presented in Section~\ref{sec:generalization} requires cutting the polyhedron in half, the program needs to consider more faces than the original polyhedron contains. Therefore, it was decided to not show the timings per face as in~\cite{Evrard}.

As in~\cite{Evrard}, it is observed that clipping a polyhedron by a paraboloid surface is about 6 times more expensive than clipping it by a plane. In turn, the clipping of a polyhedron by a quadratic cylinder is about 4.5 times more expensive than clipping it by a paraboloid and 28 times more expensive than clipping it by a plane. This reduction in speed can be explained by the fact that the clipping by a quadratic cylinders requires making two copies of the polyhedron corresponding to its intersections with the top and bottom $z$-halfspaces, to compute new vertices and faces, and to rotate the polyhedron about the $x$-axis. These manipulations involve memory allocation and an increased amount of faces to consider, in turn increasing the unit-cost of each intersection case. Note that the timing results for when a vertex is exactly on $\mathcal{S}$ correspond with the ill-posed edge cases mentioned in Section~\ref{sec:robustness}; these cases occur extremely rarely in practical applications and thus have a negligible overall impact on computational cost.

\section{Cylindrical Interface Reconstruction Results}
\label{Results2}

In this section, the PCIC reconstruction method is thoroughly tested in a series of static and dynamic cases. Section~\ref{sec:static} presents static reconstruction results and discusses the numerical properties of the method. Section~\ref{sec:translate} discusses a simple pure translation test case, while Section~\ref{sec:deform} shows a simple sinusoidal deformation advection case; both were run in the Interface Reconstruction Library (IRL)~\cite{Chiodi2}. Section~\ref{sec:stag} demonstrates the transition from PLIC to PCIC reconstructions in the NGA2 flow solver~\cite{Desjardins} using the ligament detection approach. Finally, Section~\ref{results: hit} presents a complex drop in homogeneous isotropic turbulence simulation in NGA2 and discusses Rayleigh--Plateau breakup modeling of the PCIC ligaments. 

\subsection{Static Reconstruction}
\label{sec:static}
The numerical properties and basic performance of PCIC were studied with a static convergence analysis. A $5\times5\times5$ stencil of cells was considered with a straight circular cylinder used to initialize volume moments in the cells. From these moments, the PCIC method was used for the center cell, and the resulting reconstruction was then compared with the initial cylinder. All grid cells had non-dimensional edge lengths of $1$. The radius of the initial cylinder was varied to study convergence. $10,000$ cylinders with random origins and orientations that cut through the center cell were generated for each radius. The largest radius considered was $1$, which is the mesh size $\Delta x$, and then the radius was halved in each successive data set. To quantify error, three metrics were considered: the orientation error ($E_{\theta}$), the origin error ($E_{o}$), and the liquid barycenter error ($E_{b}$). The orientation error was calculated as the angle between the reconstructed orientation vector and the actual orientation vector
\begin{equation}
\label{theta}
\begin{aligned}
E_{\theta}=\frac{1}{n}\sum_{i=1}^{n}\arccos\left(\left|\mathbf{v}_i\cdot\mathbf{v}_{a_i}\right|\right),
\end{aligned}
\end{equation}
where $\mathbf{v}_i$ is the reconstructed orientation vector at the $i^{th}$ sample and $\mathbf{v}_{a_i}$ is the actual orientation vector at the $i^{th}$ sample. The origin error was calculated with the magnitude of the cross product between $\mathbf{v}_a$ and the difference between the reconstructed origin and the actual origin
\begin{equation}
\label{origin}
\begin{aligned}
E_{o}=\frac{1}{n}\sum_{i=1}^{n}\left|\left|\left(\mathbf{p}_i-\mathbf{p}_{a_i}\right)\times\mathbf{v}_{a_i}\right|\right|,
\end{aligned}
\end{equation}
where $\mathbf{p}_i$ is the reconstructed origin at the $i^{th}$ sample and $\mathbf{p}_{a_i}$ is the actual origin at the $i^{th}$ sample. This measures the perpendicular distance between the reconstructed origin and the cylinder's orientation axis, which is the appropriate metric. Finally, the liquid barycenter error was defined as
\begin{equation}
\label{bary}
\begin{aligned}
E_{b}=\frac{1}{n}\sum_{i=1}^{n}\left|\left|\mathbf{x}_{i}-\mathbf{x}_{a_{i}}\right|\right|,
\end{aligned}
\end{equation}
where $\mathbf{x}_i$ is the center cell's liquid barycenter from the reconstructed interface for the $i^{th}$ sample and $\mathbf{x}_{a_{i}}$ is the corresponding liquid barycenter from the actual interface.

\begin{figure} \centering
    \includegraphics{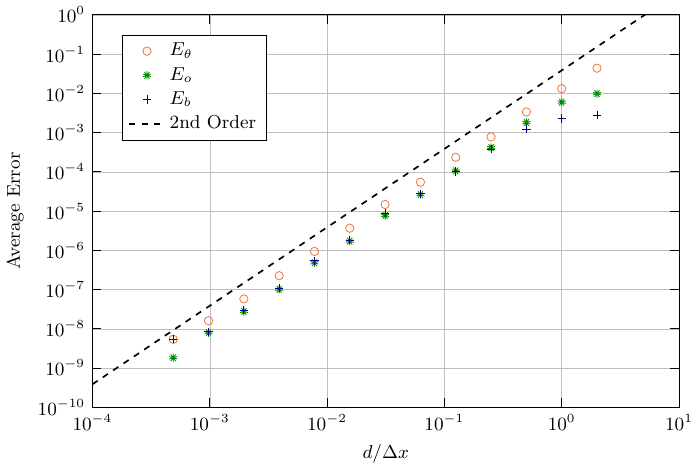}
    \caption{Convergence plot of the orientation error $E_{\theta}$, origin error $E_o$, and barycenter error $E_b$ for PCIC tested with static reconstructions. $d/\Delta x$ is the cylinder diameter divided by the grid size. Each data point is the average of \num{1e4} samples, as defined in Eqs.~\eqref{theta},~\eqref{origin}, and~\eqref{bary}.}
    \label{convergence}
\end{figure}

Figure~\ref{convergence} shows the resulting error plot. For all of the error metrics, PCIC demonstrates approximately second-order accuracy on average with respect to decreasing cylinder radius, especially at the sub-grid scale. It must be noted that, since PCIC uses a $5\times5\times5$ stencil, the largest cylinder that can be safely reconstructed regardless of origin and orientation has a radius equal to the mesh size (note: this corresponds with the choice made for CCL ligament detection in Section~\ref{Methods:lig_det}). This ensures that the entire cross section of the cylinder is contained within the stencil. Since the observed average convergence is with respect to decreasing radius, another interpretation of this result is that the method generally converges with increasing grid size instead of decreasing grid size. Although such a property would be problematic with most traditional geometric VOF reconstructions, for sub-grid ligament reconstruction this is a favorable result: the performance improves at smaller and smaller scales. Consider this in the context of a multiphase flow simulation; if a PLIC ligament is converted to a PCIC ligament at a diameter of about twice the mesh size, then the maximum error is made when the transfer occurs. This is because PLIC errors increase as the ligament approaches sub-grid scale. Then, when converted to PCIC, the reconstruction becomes more and more accurate as the ligament continues to become thinner. 

\begin{figure} \centering
    \includegraphics{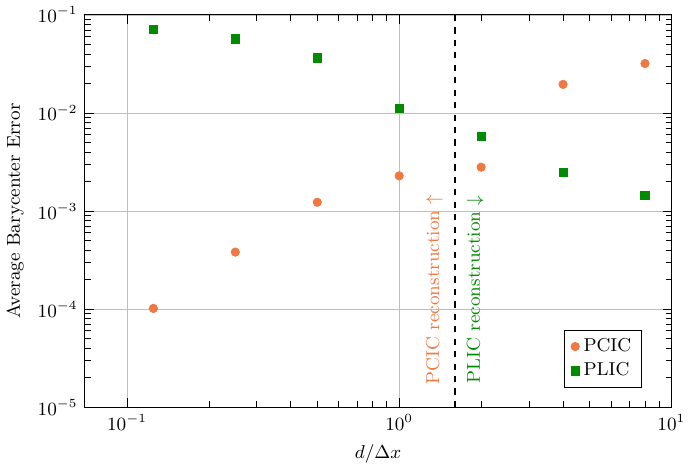}
    \caption{Comparison of the barycenter error $E_b$ between PCIC and PLIC tested with static reconstructions. $d/\Delta x$ is the cylinder diameter divided by the grid size. Each data point is the average of \num{1e4} samples, as defined in Eq.~\eqref{bary}. The vertical dashed line represents the cutoff where an interface would be switched from a PLIC to a PCIC reconstruction in a simulation.}
    \label{pcic_vs_plic}
\end{figure}

This can be quantitatively demonstrated by comparing the errors produced by PLIC and PCIC near the transition for \num{1e4} samples on the $5\times5\times5$ stencil. The barycenter error $E_{b}$ introduced above was used for this test, and PLIC-Net~\cite{Cahaly} was used for PLIC reconstructions (which uses a $3\times3\times3$ stencil). The results are shown in Fig.~\ref{pcic_vs_plic}.

At a cylinder radius of $2\Delta x$, which is larger than the maximum PCIC safe radius, the PCIC error is an order of magnitude higher than the PLIC error. As expected, this demonstrates that PLIC is superior for more resolved ligaments. At a radius of $\Delta x$, PCIC becomes more accurate than PLIC. However, it is observed that cleaner simulation results are obtained when the transition to PCIC is made at a slightly lower radius of $0.8\Delta x$. Regardless, when the transition to PCIC is made, there is a jump in accuracy as desired. At a radius of $\Delta x/2$, the trend continues as PCIC continues to improve and PLIC worsens.

For a cylinder that is mesh-aligned (i.e., one component of $\mathbf{v}$ is $1$ and the other two are $0$), the PCIC reconstruction method is exact to machine precision regardless of radius. Errors are introduced in mesh-unaligned cases due to barycenter data asymmetries caused by using a cubic stencil in all cases. The chosen stencil size of $5\times5\times5$ helps mitigate these errors and ensures satisfactory results. As shown above, on average the errors converge with second-order accuracy with decreasing radius; however, there are exceptions for individual samples, such as the mesh-aligned case which has zeroth-order convergence. Indeed, for any given cylinder, there is no expectation of a certain order of accuracy, meaning that the previous convergence results are only true for averaged data sets.

\subsection{Translation Advection}
\label{sec:translate}
\begin{figure}[b!]
   \centering
   \begin{subfigure}{.4\textwidth}
       \centering
       \includegraphics[width=\linewidth, clip]{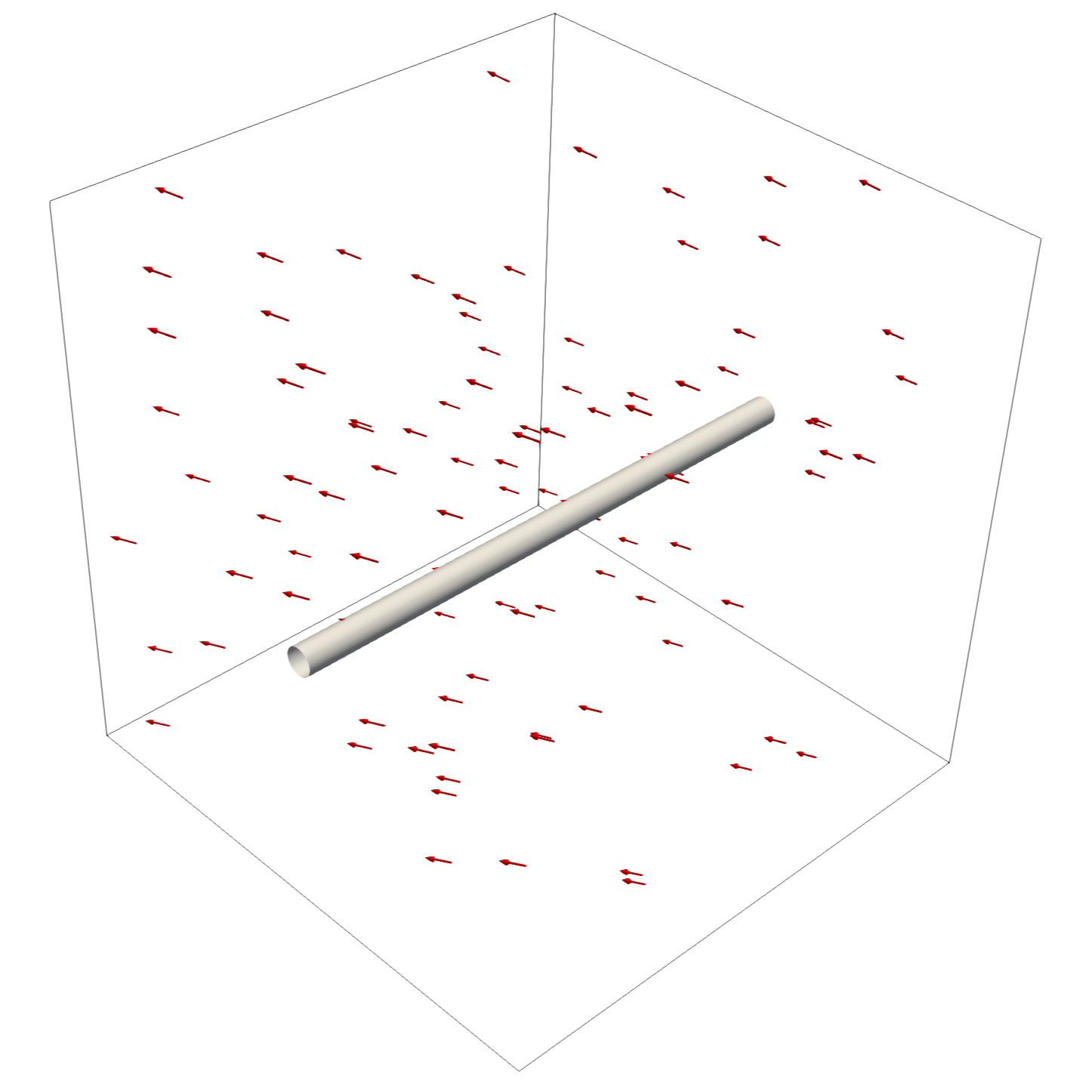}
       \caption{$r = 0.0125$}
       \label{translate1:sub1}
   \end{subfigure}%
   \begin{subfigure}{.4\textwidth}
       \centering
       \includegraphics[width=\linewidth, clip]{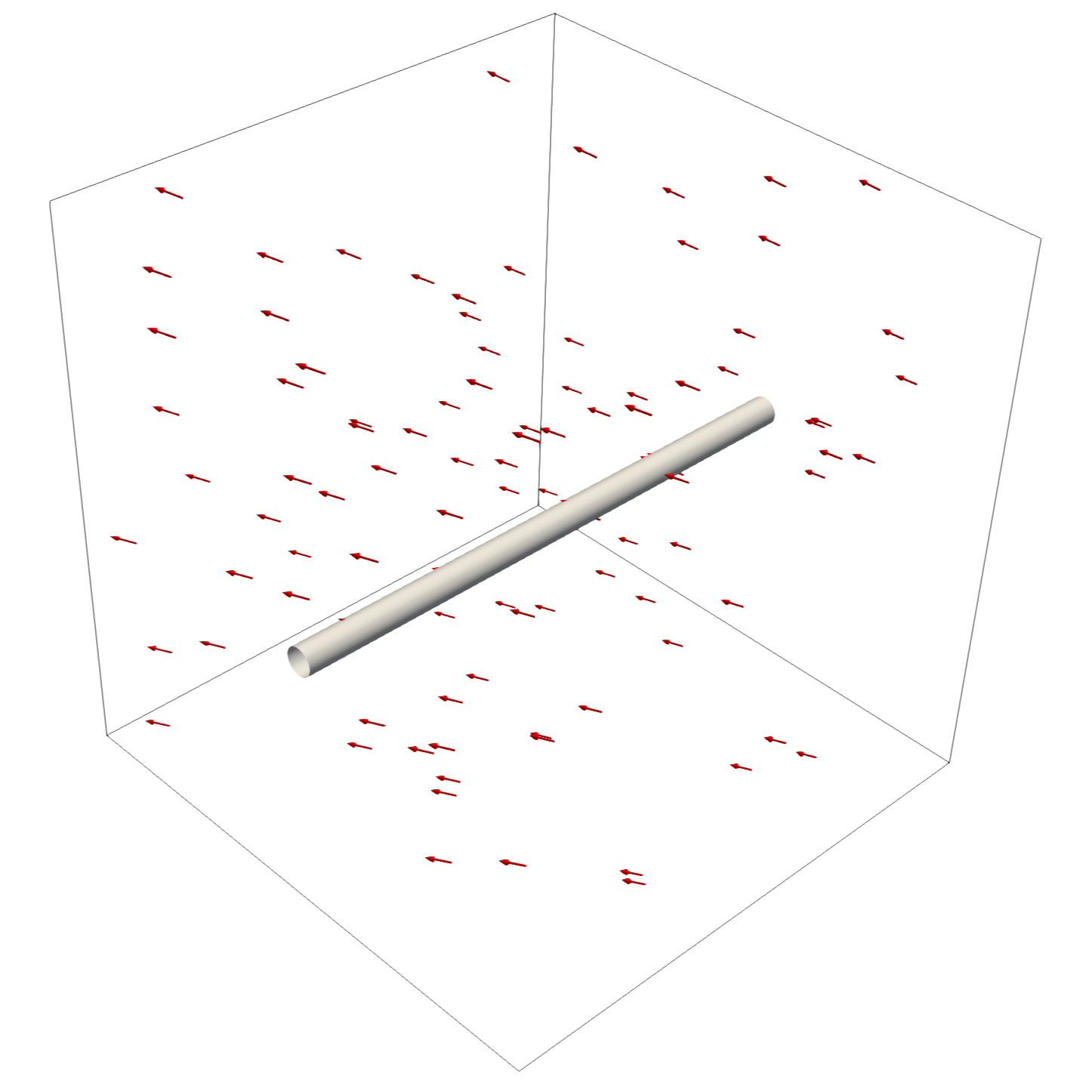}
       \caption{$r = 0.025$}
       \label{translate1:sub2}
   \end{subfigure}%
  \caption{Final time step of the horizontal advection case. The velocity vector field is shown.}
  \label{horz_translate}
\end{figure}

\begin{figure} [t!]
   \centering
   \begin{subfigure}{.31\textwidth}
       \centering
       \includegraphics[width=\linewidth, clip]{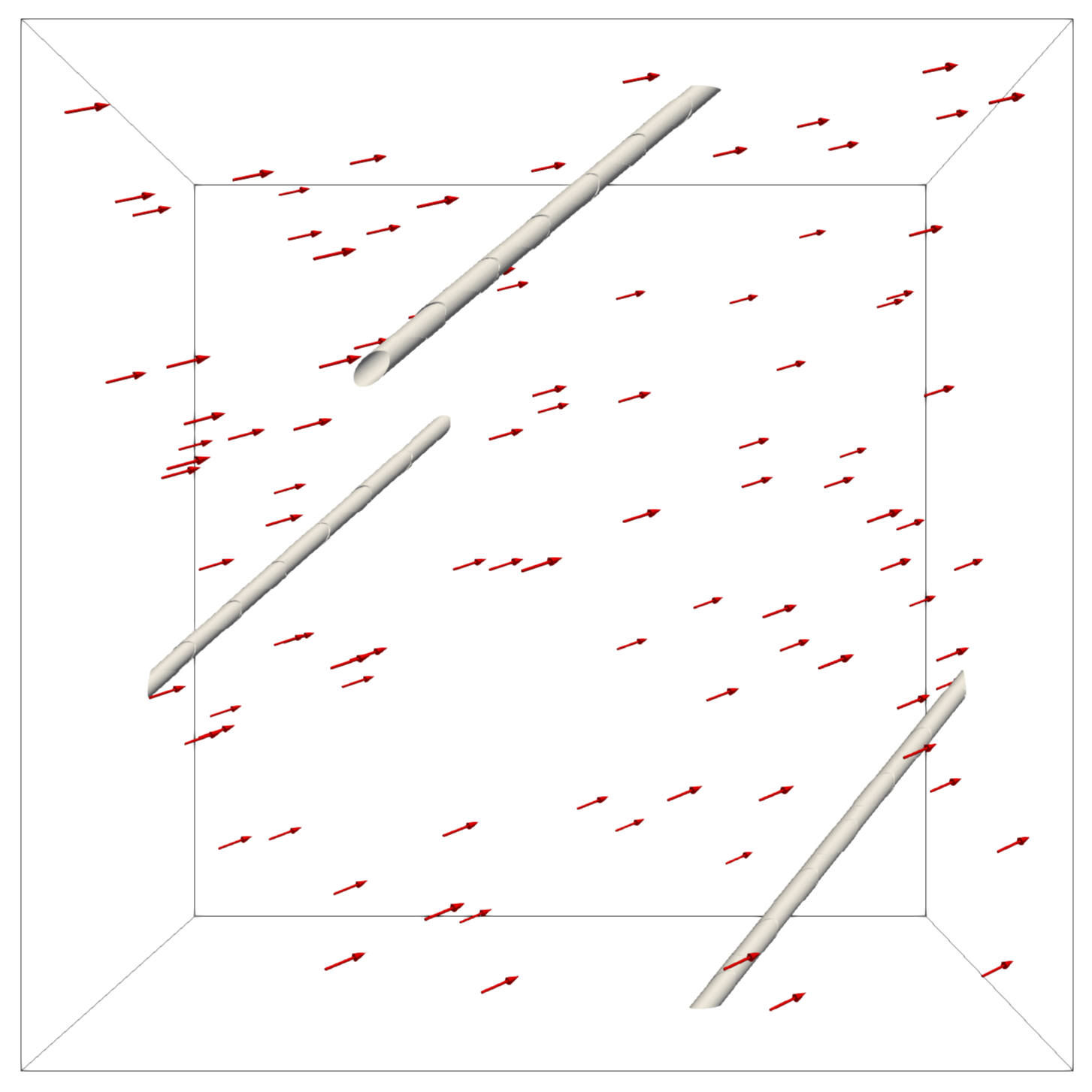}
       \caption{ $t=1$}
       \label{diag_translate:sub1}
   \end{subfigure}%
   \begin{subfigure}{.31\textwidth}
       \centering
       \includegraphics[width=\linewidth, clip]{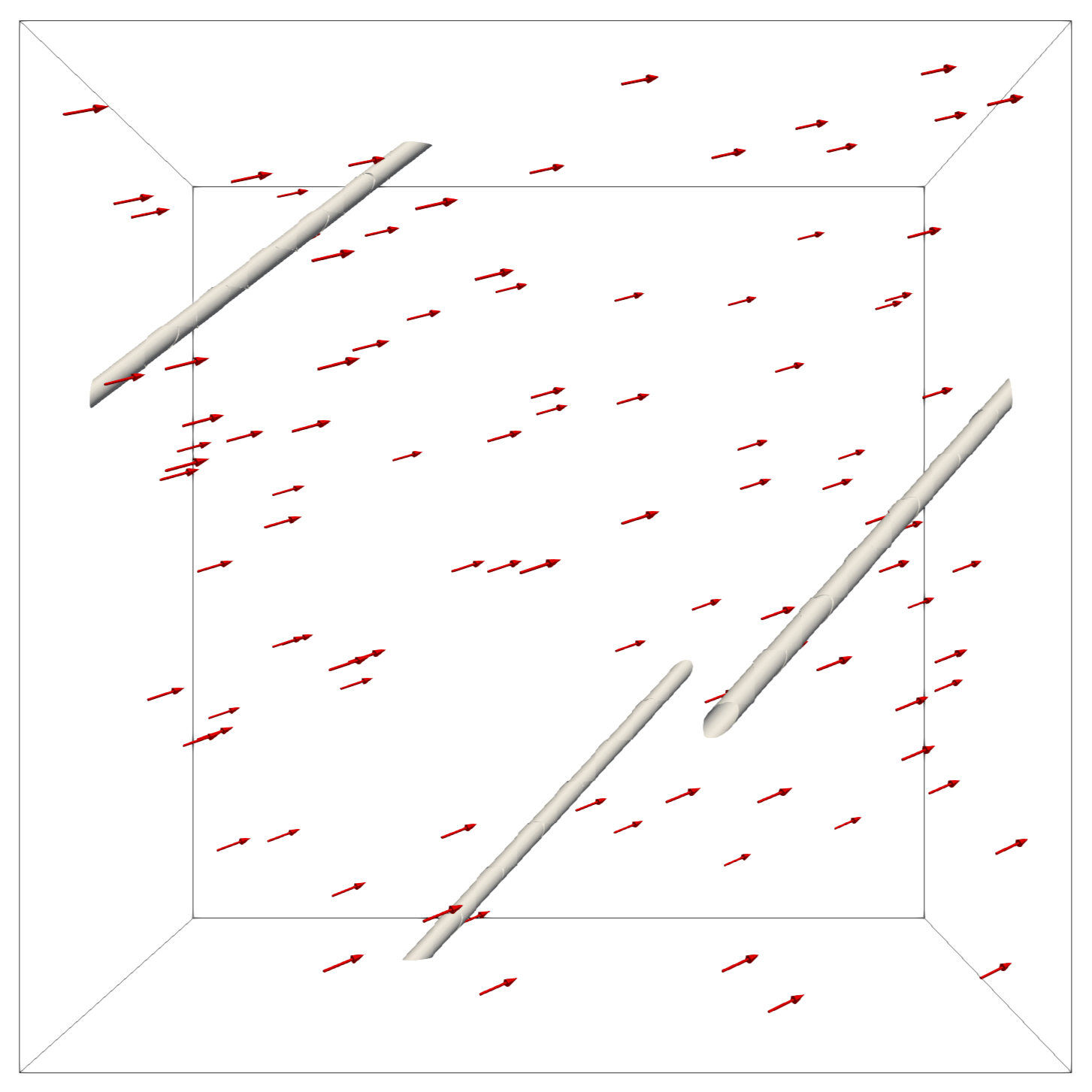}
       \caption{ $t = 2$}
       \label{diag_translate:sub2}
   \end{subfigure}%
   \begin{subfigure}{.31\textwidth}
       \centering
       \includegraphics[width=\linewidth, clip]{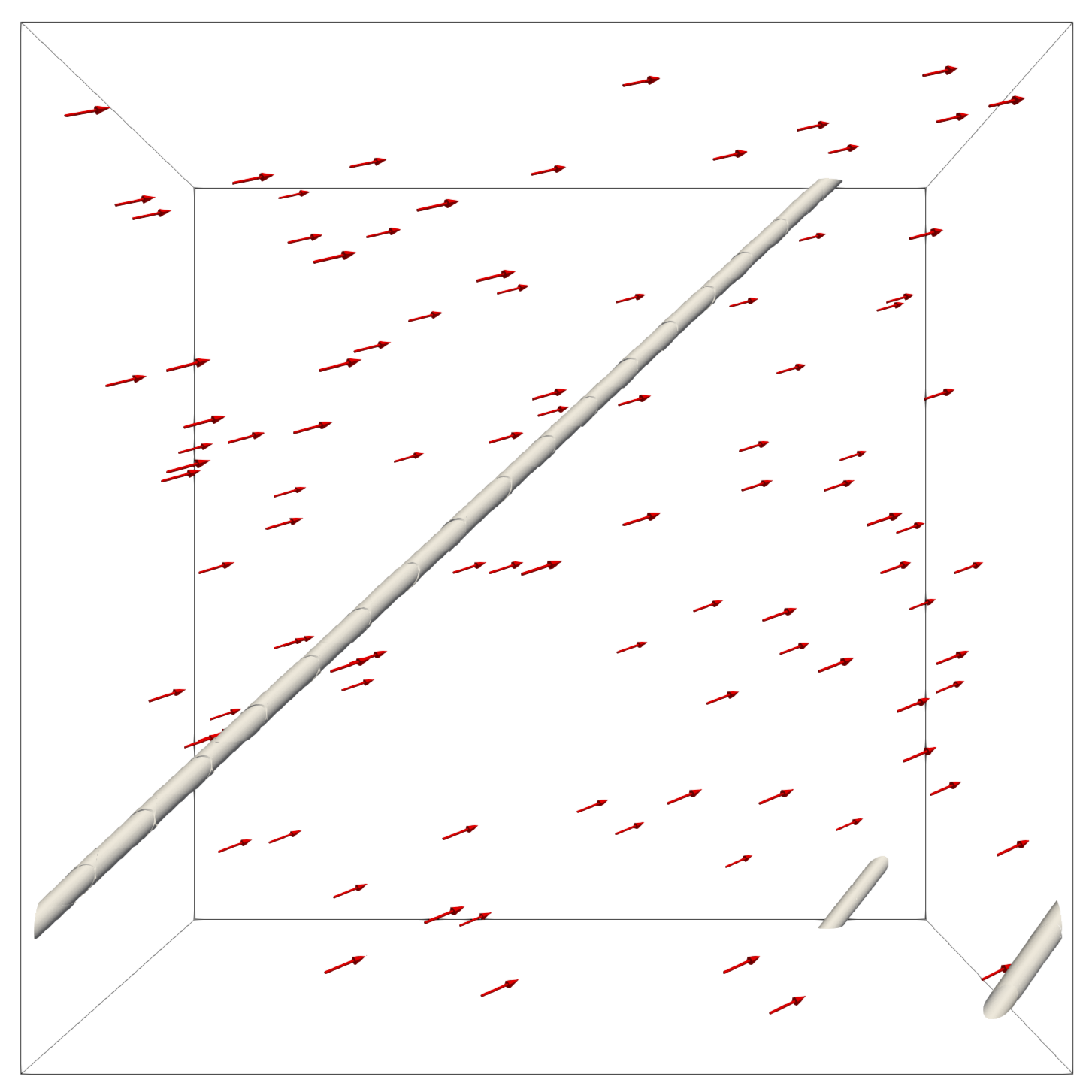}
       \caption{ $t = 2.8$}
       \label{diag_translate:sub3}
   \end{subfigure}%
   \\
   \centering
   \begin{subfigure}{.31\textwidth}
       \centering
       \includegraphics[width=\linewidth, clip]{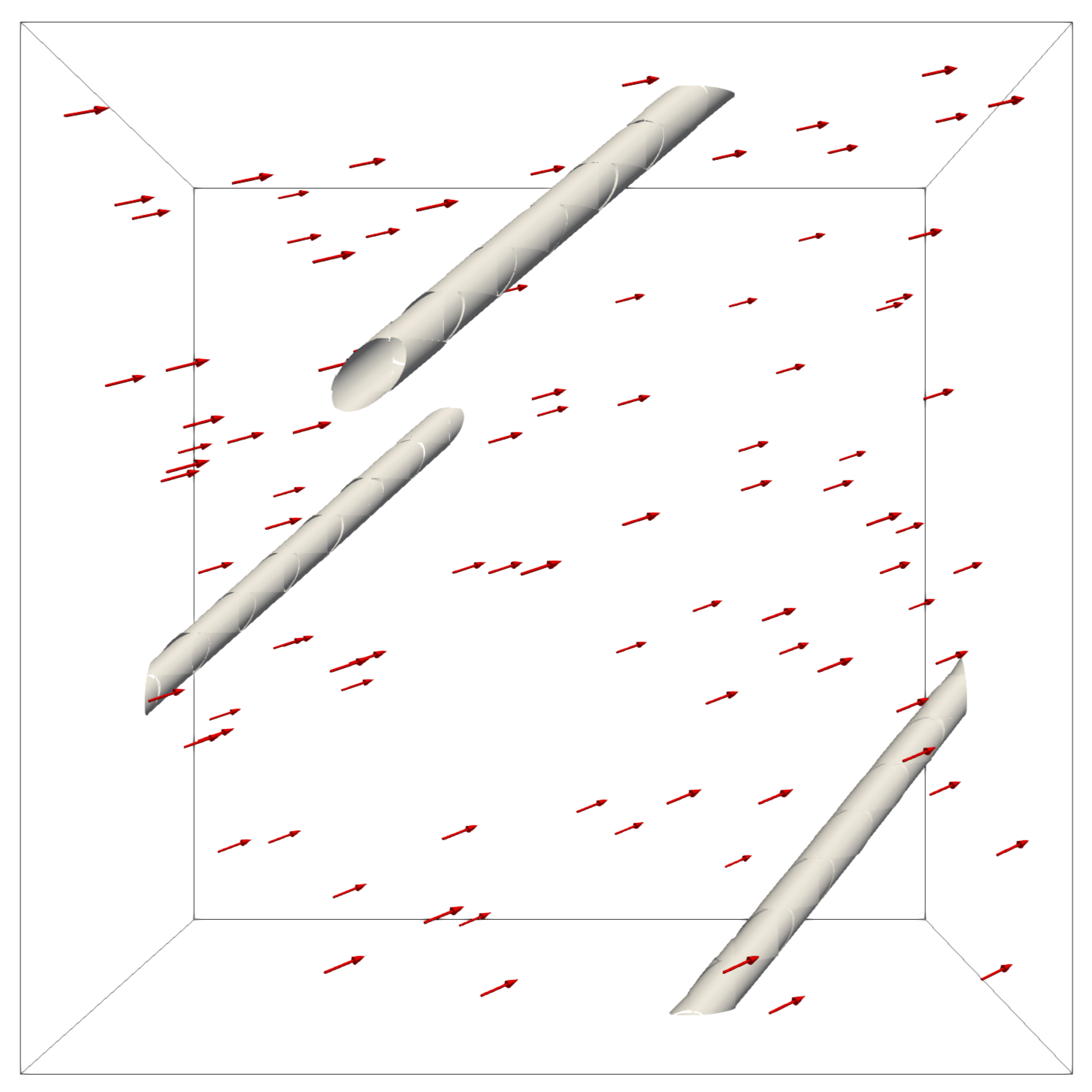}
       \caption{ $t = 1$}
       \label{diag_translate:sub4}
   \end{subfigure}%
   \begin{subfigure}{.31\textwidth}
       \centering
       \includegraphics[width=\linewidth, clip]{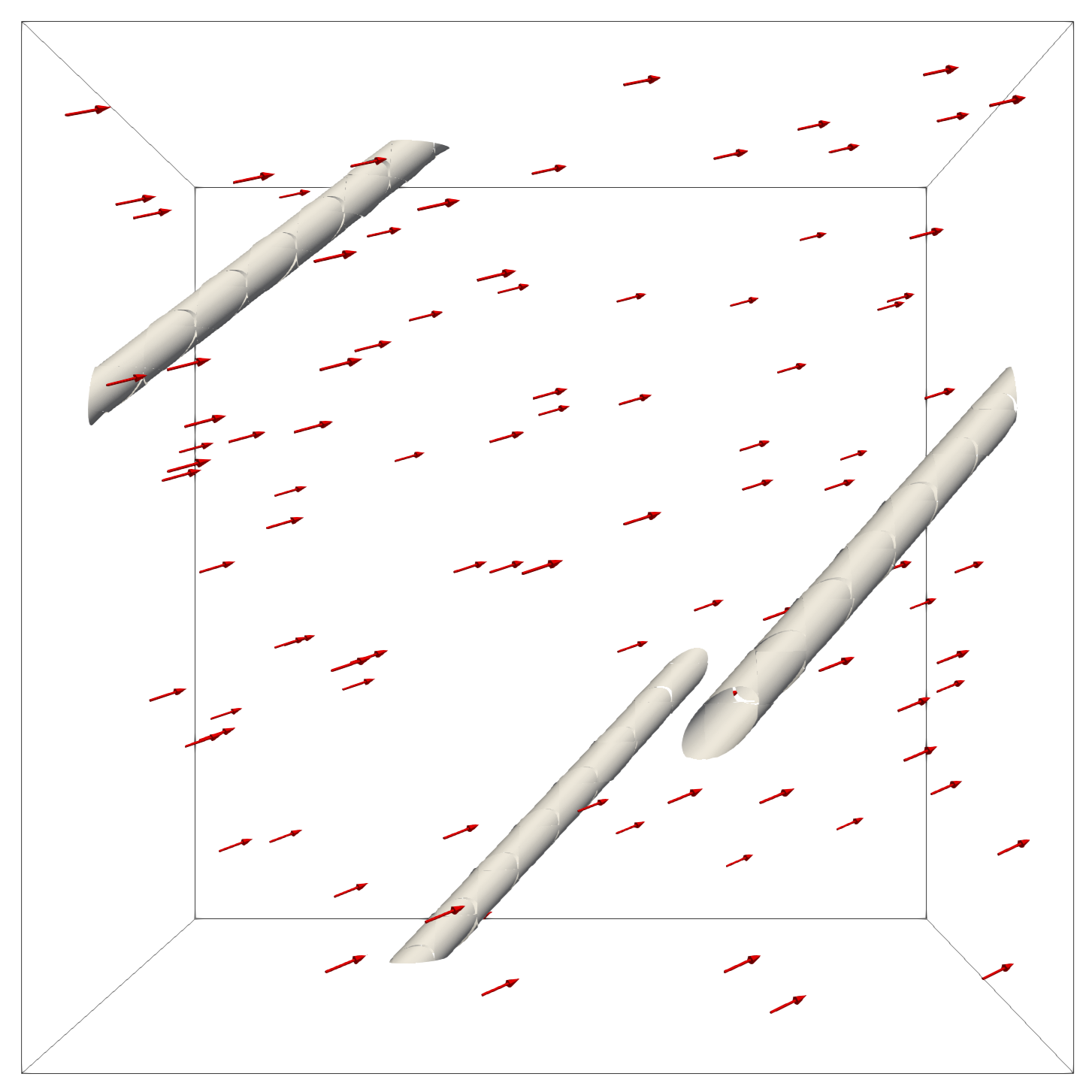}
       \caption{ $t = 2$}
       \label{diag_translate:sub5}
   \end{subfigure}%
   \begin{subfigure}{.31\textwidth}
       \centering
       \includegraphics[width=\linewidth, clip]{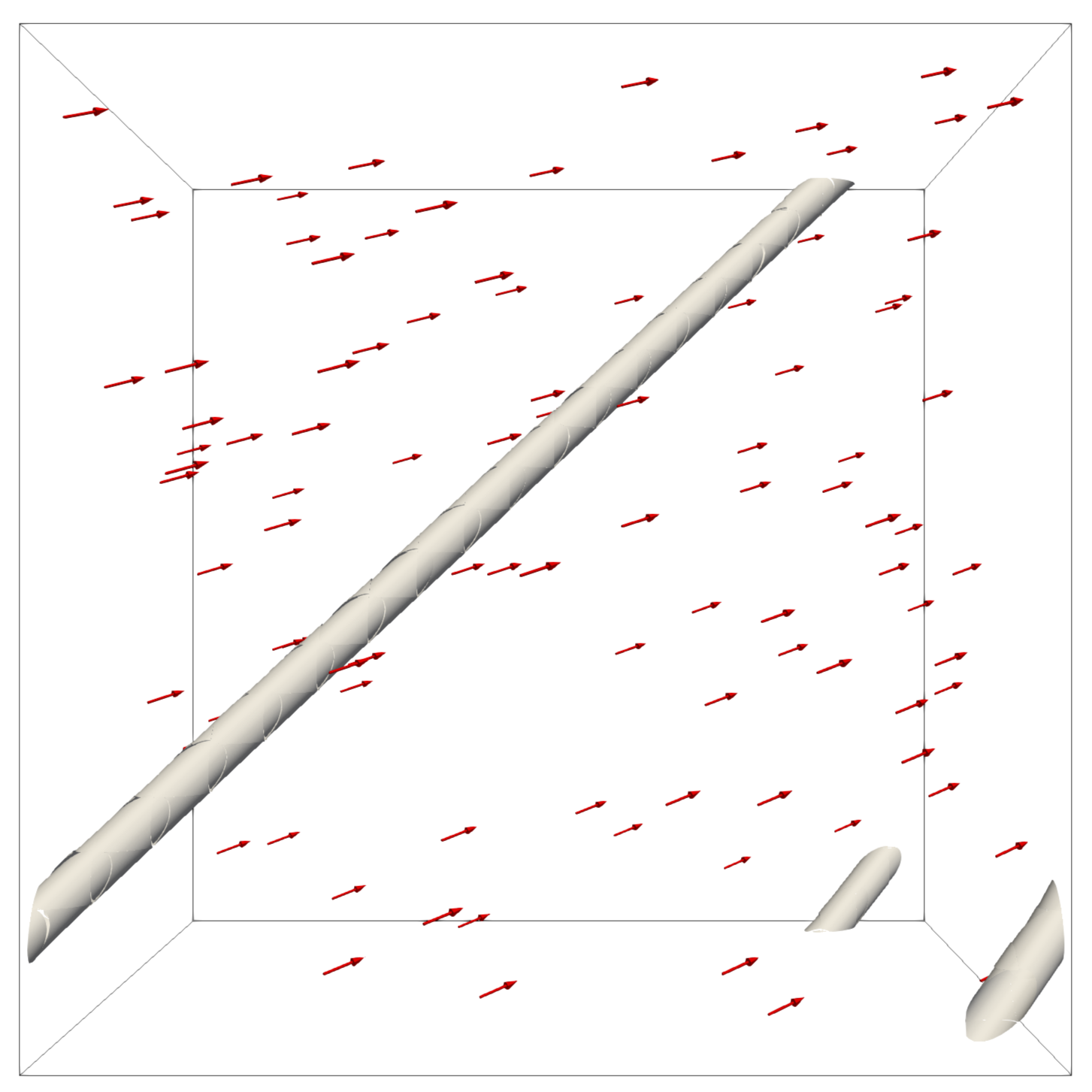}
       \caption{ $t = 2.8$}
       \label{diag_translate:sub6}
   \end{subfigure}%
  \caption{Time series of the diagonal advection case. (a), (b), and (c) have a radius $r=0.0125$ and (d), (e), and (f) have a radius $r=0.025$. All times are non-dimensional. The velocity vector field is shown.}
  \label{diag_translate}
\end{figure}

Simple advection test cases were conducted within the Interface Reconstruction Library (IRL)~\cite{Chiodi2} using a prescribed velocity to transport the volume fraction field. A pure translation case of a stright circular cylinder was first run to provide an initial benchmark of PCIC's performance in cases involving VOF transport. The domain was set to be a unit cube with the cylinder initialized in the center. One case was run with the cylinder aligned with the $x$-axis (the horizontal case), and another with a cylinder orientation of $\langle\frac{1}{\sqrt{3}},\frac{1}{\sqrt{3}},\frac{1}{\sqrt{3}}\rangle$ (the diagonal case). The cylinder was given a non-dimensional velocity in the $x$ direction of $1$, in the $y$ direction of $\frac{2}{3}$, and in the $z$ direction $\frac{1}{3}$. Periodic boundary conditions were imposed such that the cylinder would return to its initial position after a non-dimensional time of $3$. A constant CFL number of $0.8$ was used. $20$ grid cells were used in each direction, and cases were run with cylinder radii of $0.0125$ and $0.025$. The results for the horizontal cases at $t=3$ are shown in Fig.~\ref{horz_translate}, and time series of the diagonal cases are shown in Fig.~\ref{diag_translate}. In all cases, PCIC performs well. No spurious interfaces are created and the cylinder maintains its shape during the entire simulation. In the diagonal case with the larger radius, there are visually noticeable reconstruction errors in agreement with the static reconstruction results of the preceding subsection. However, the errors are at an acceptable level and performance remains robust. As expected, the case with the smaller radius shows less noticeable error, and the horizontal cases demonstrate no visible errors regardless of the radius.

\begin{figure}
   \centering
    \begin{subfigure}{.33\textwidth}
       \centering
       \includegraphics[width=\linewidth, clip]{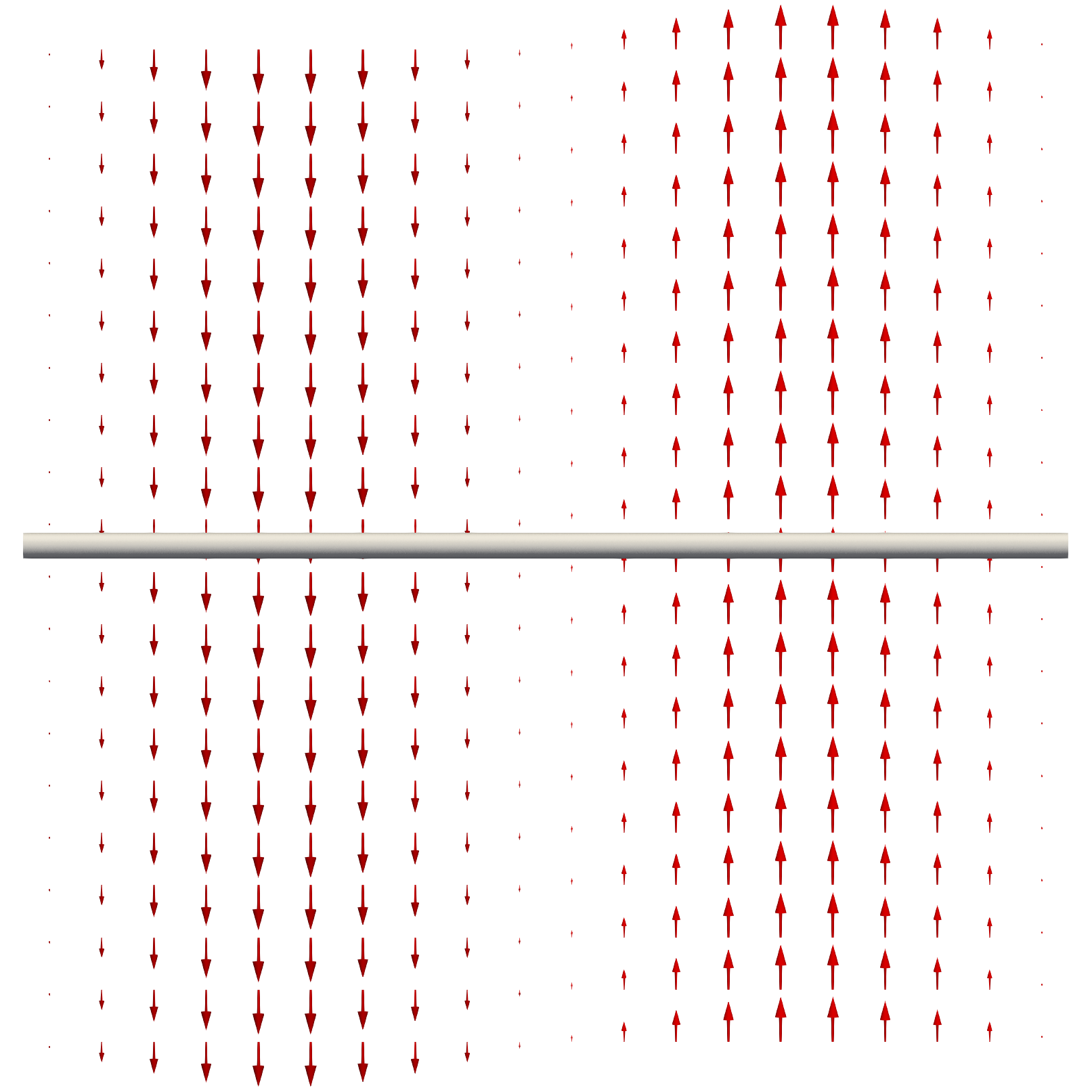}
       \caption{$t = 0$}
       \label{stag:deform0}
   \end{subfigure}%
   \hspace{0.4em}
   \begin{subfigure}{.33\textwidth}
       \centering
       \includegraphics[width=\linewidth, clip]{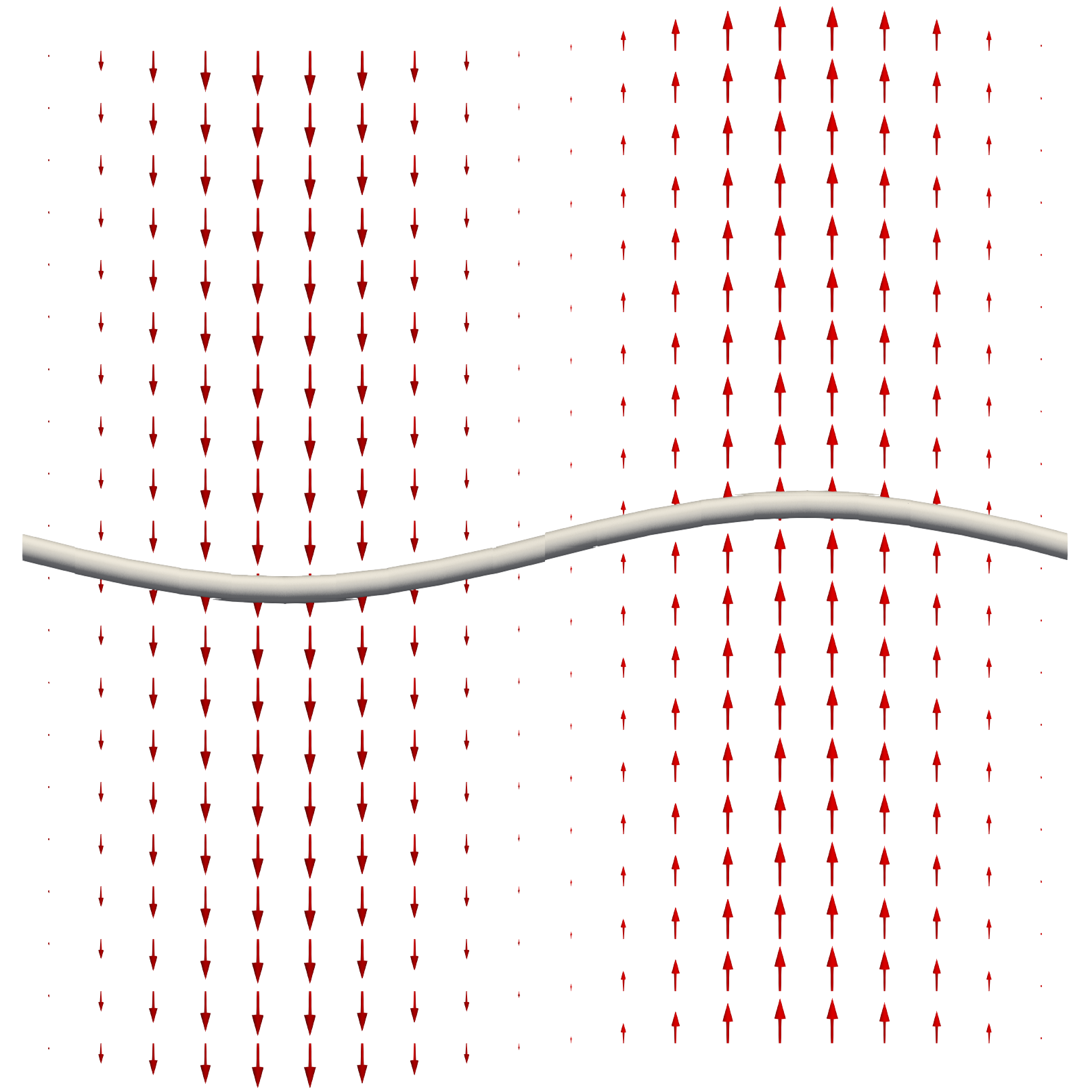}
       \caption{$t = 0.02025$}
       \label{stag:deform1}
   \end{subfigure}%
   \\
   \begin{subfigure}{.33\textwidth}
       \centering
       \includegraphics[width=\linewidth, clip]{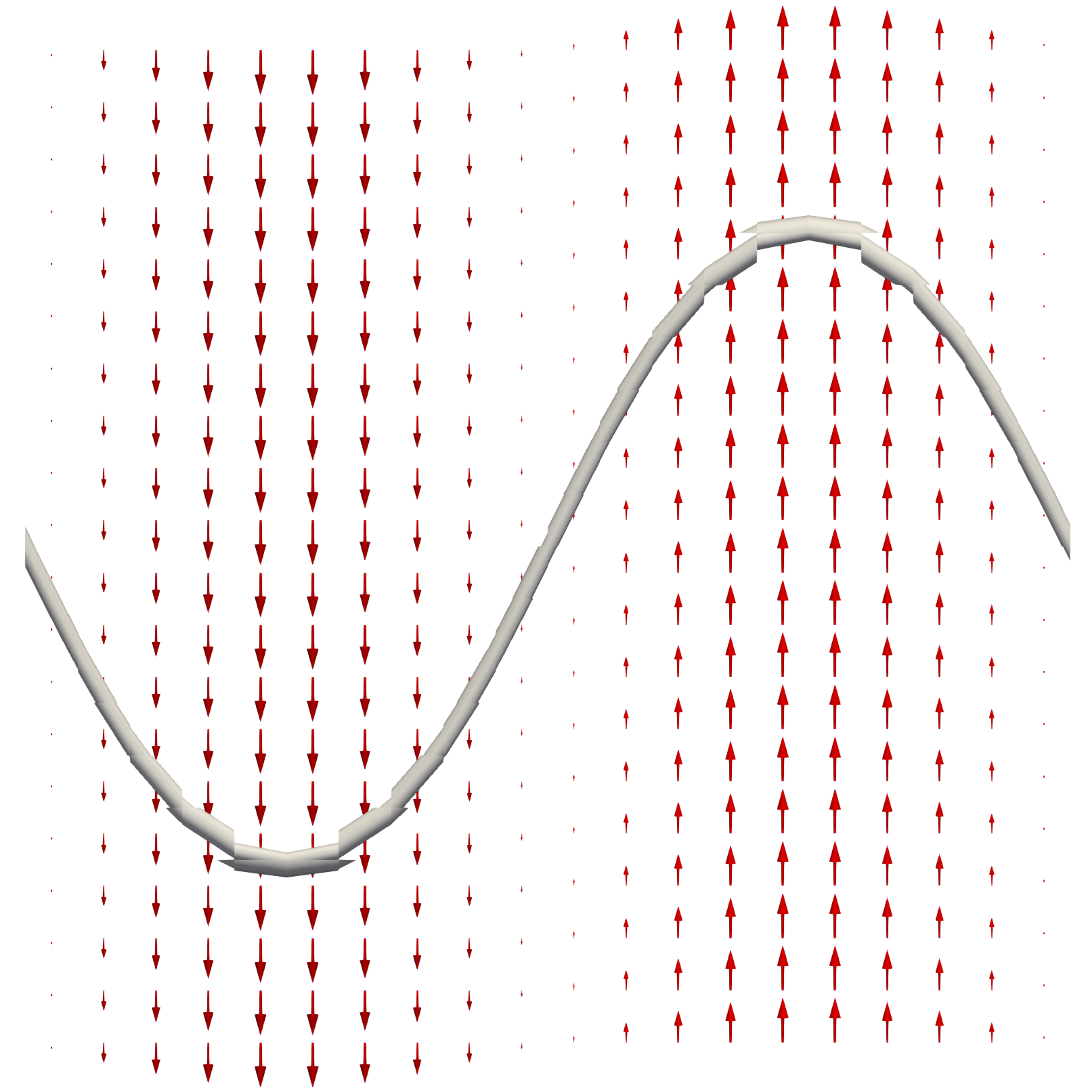}
       \caption{$t = 0.162$}
       \label{deform:sub2}
   \end{subfigure}%
   \hspace{0.4em}
   \begin{subfigure}{.33\textwidth}
       \centering
       \includegraphics[width=\linewidth, clip]{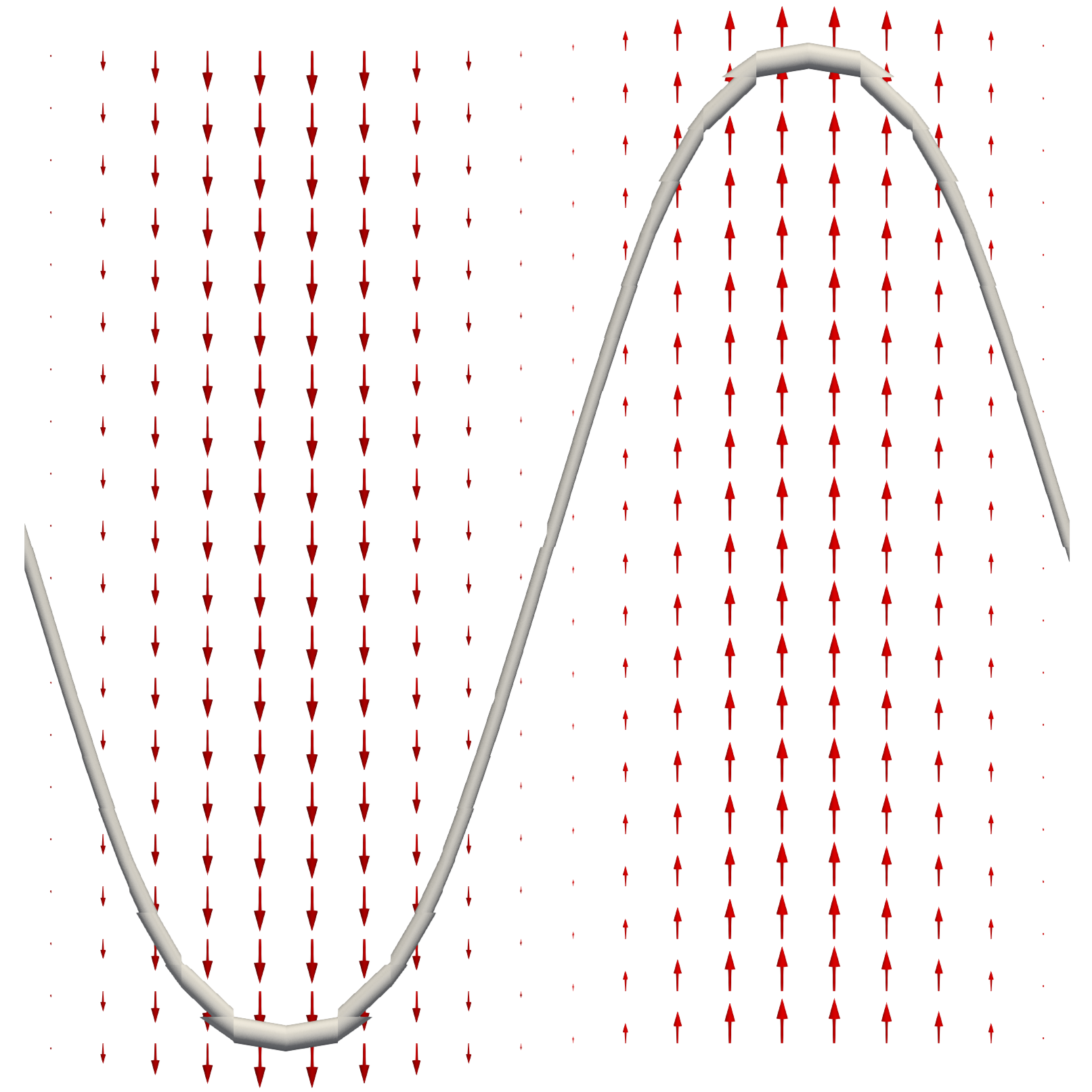}
       \caption{$t = 0.26325$}
       \label{deform:sub3}
   \end{subfigure}%
  \caption{Time series of the sinusoidal deformation case. All times are non-dimensional. The velocity vector field is shown.}
  \label{deformation_series}
\end{figure}

\begin{figure}
   \centering
    \begin{subfigure}{.48\textwidth}
       \centering
       \includegraphics[width=\linewidth, clip]{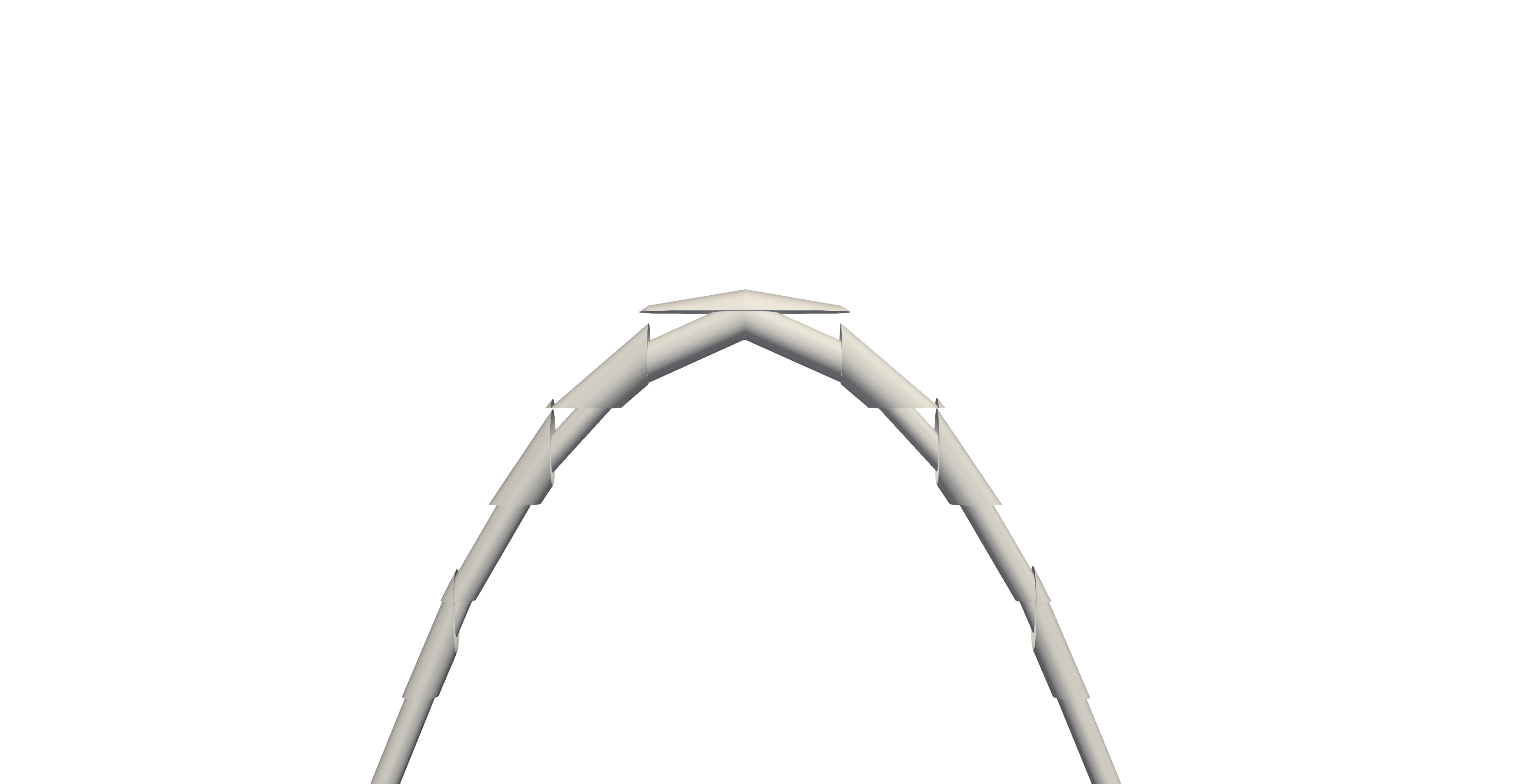}
       \caption{PCA-Only Method}
       \label{closeup:PCA}
   \end{subfigure}%
   \hspace{0.4em}
   \begin{subfigure}{.48\textwidth}
       \centering
       \includegraphics[width=\linewidth, clip]{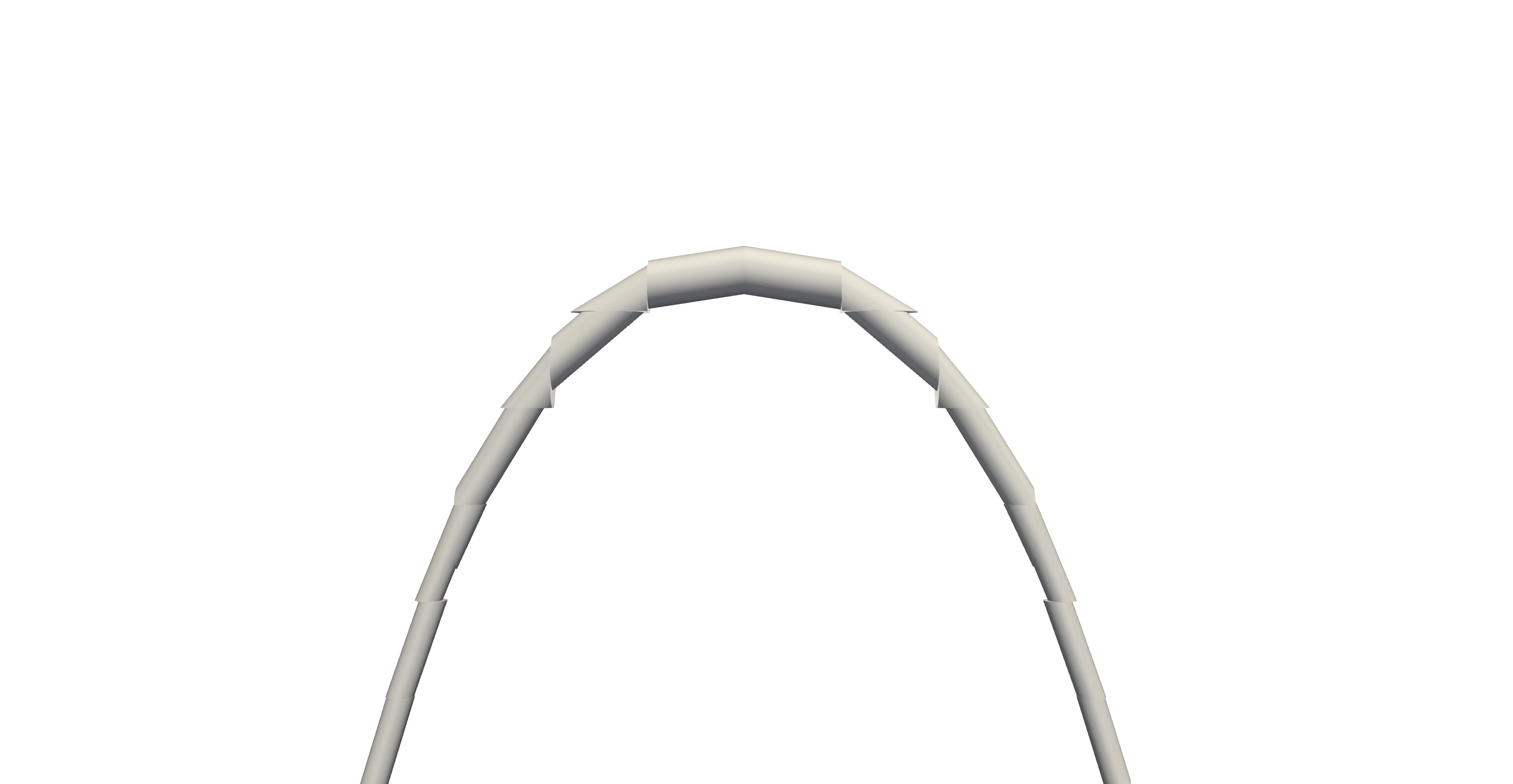}
       \caption{PCIC (quadratic curve fitting)}
       \label{closeup:global}
   \end{subfigure}%
  \caption{Comparison of a PCA-only method and the PCIC algorithm. The quadratic curve fitting provides for a cleaner reconstruction.}
  \label{closeup}
\end{figure}

\begin{figure}
  \begin{center}
   \includegraphics[width=8.8cm]{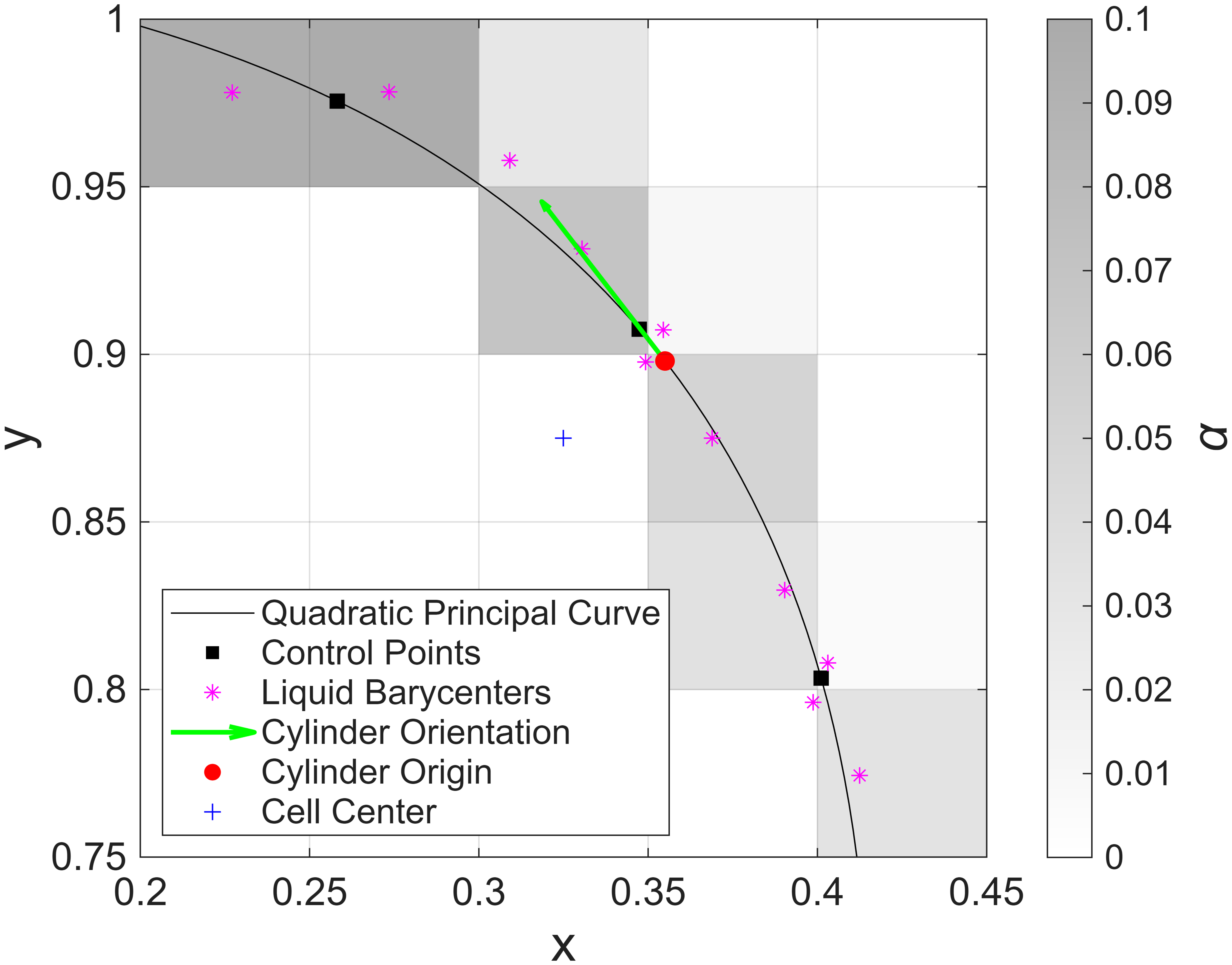}
  \end{center}
  \caption{Example of PCIC reconstruction in a stencil from the deformation advection case. The liquid barycenters are shown, along with the resulting quadratic principal curve and its control points. The cylinder origin and orientation vector are indicated. The cells are shaded according to their volume fractions.}
  \label{PCIC_fit}
\end{figure}

\subsection{Deformation Advection}
\label{sec:deform}
A more complex advection case was run in IRL involving a prescribed sinusoidal deformation of an initially straight circular cylinder. This case was designed to test PCIC's performance for ligaments undergoing high-curvature deformations. Like with the translation case, the domain was a periodic unit cube with the cylinder initialized in the center and aligned with the $x$-axis. A non-dimensional velocity of $0$ in the $x$ and $z$ directions and of $2\sin{\left(2\pi x\right)}$ in the $y$ direction was prescribed. $20$ grid cells were used in each direction, and the initial cylinder was given a radius of $0.0125$. A maximum CFL number of $0.8$ was used. A time series of the case is shown in Fig.~\ref{deformation_series}. Despite being subjected to high-curvature bending, PCIC performs well and the ligament does not break. It is observed that PCIC's quadratic principal curve fitting provides for a cleaner reconstruction over a simpler PCA-only approach, as demonstrated in Fig.~\ref{closeup}. Figure~\ref{PCIC_fit} shows an example at the final time step ($t=0.26325$) of the control points returned by the principal curve fitting, along with the cylinder origin and the orientation vector. The algorithm does an excellent job of fitting the sinusoidal curve and providing proper orientation vectors for the cylinder reconstruction.

\subsection{Drop in Stagnation Flow}
\label{sec:stag}
In order to test the transition from a PLIC reconstruction to a PCIC reconstruction, a drop in a stagnation flow was considered. The domain was given a length of $1$ unit in the $x$ and $y$ directions with $10$ grid cells, and a length of $4$ units in the $z$ direction with $40$ grid cells. A prescribed divergence-free stagnation flow was simply defined as
\begin{equation}
\begin{bmatrix} u \\ v \\ w \end{bmatrix} = \begin{bmatrix} -x \\ -y \\ 2z \end{bmatrix},
\end{equation}
where $u$, $v$, and $w$ are the velocities in the $x$, $y$, and $z$ directions, respectively, with non-dimensional units. A spherical drop with a radius of $0.2$ was initialized with its center at the flow's stagnation point (in the center of the domain). A time step of $0.02$ was used, which corresponds to a maximum CFL number of $0.8$. 

\begin{figure} [b!]\vspace{-1em}
   \centering
      \begin{subfigure}{.22\textwidth}
       \centering
       \includegraphics[width=\linewidth, clip]{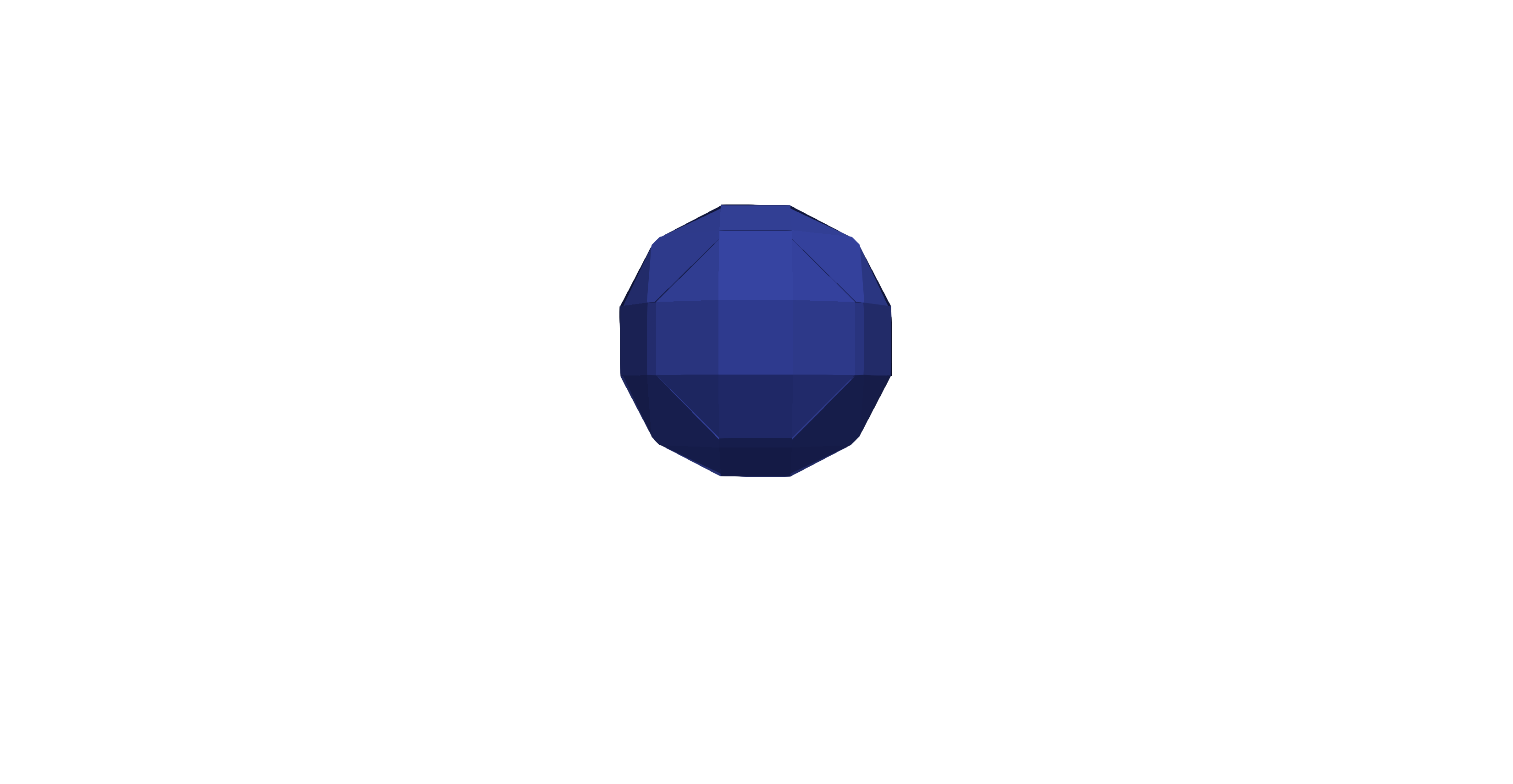}
       \caption{$t = 0$}
       \label{stag:sub0}
   \end{subfigure}%
   \hspace{0.4em}
      \begin{subfigure}{.22\textwidth}
       \centering
       \includegraphics[width=\linewidth, clip]{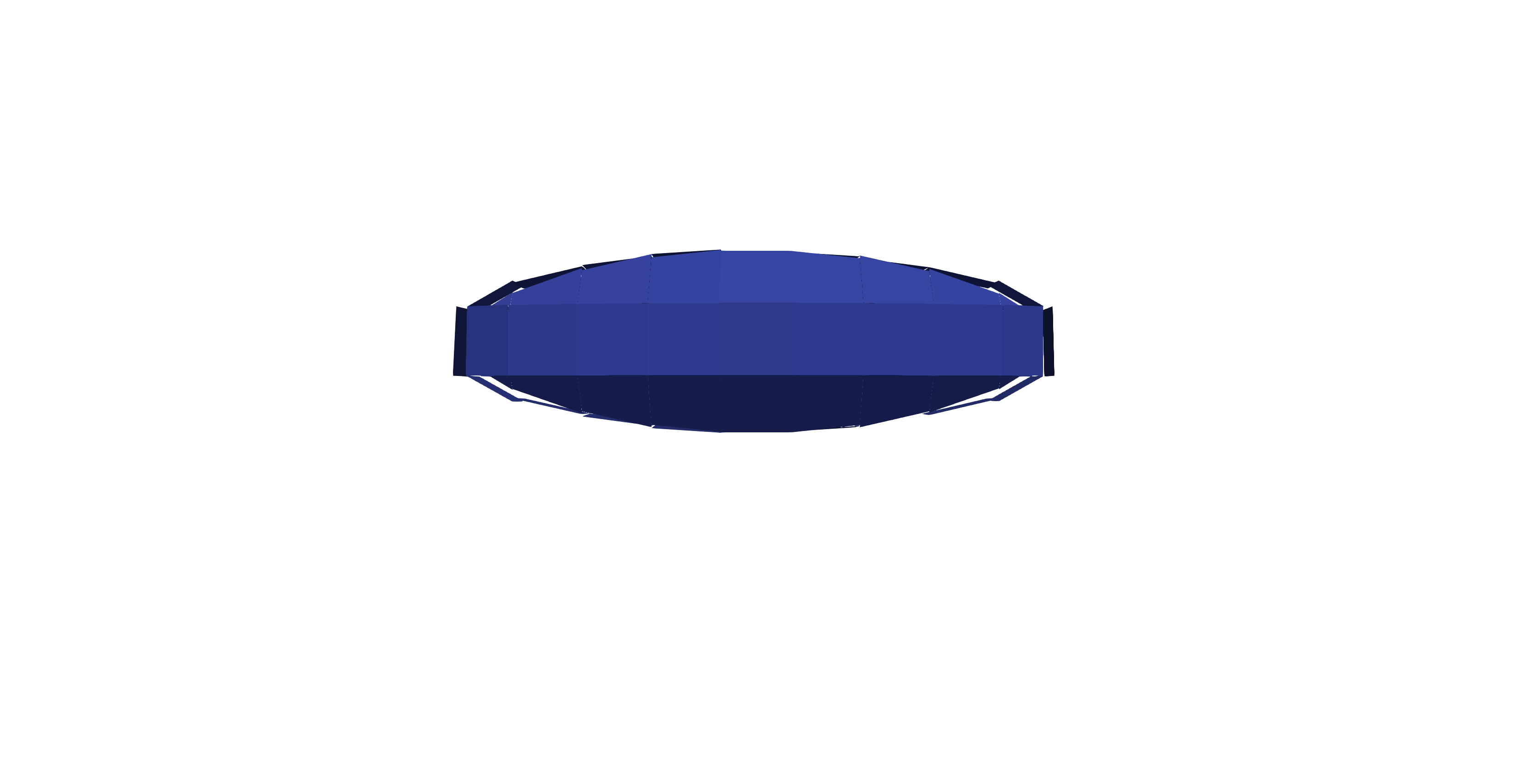}
       \caption{$t = 0.4$}
       \label{stag:sub00}
   \end{subfigure}%
   \hspace{0.4em}
   \begin{subfigure}{.22\textwidth}
       \centering
       \includegraphics[width=\linewidth, clip]{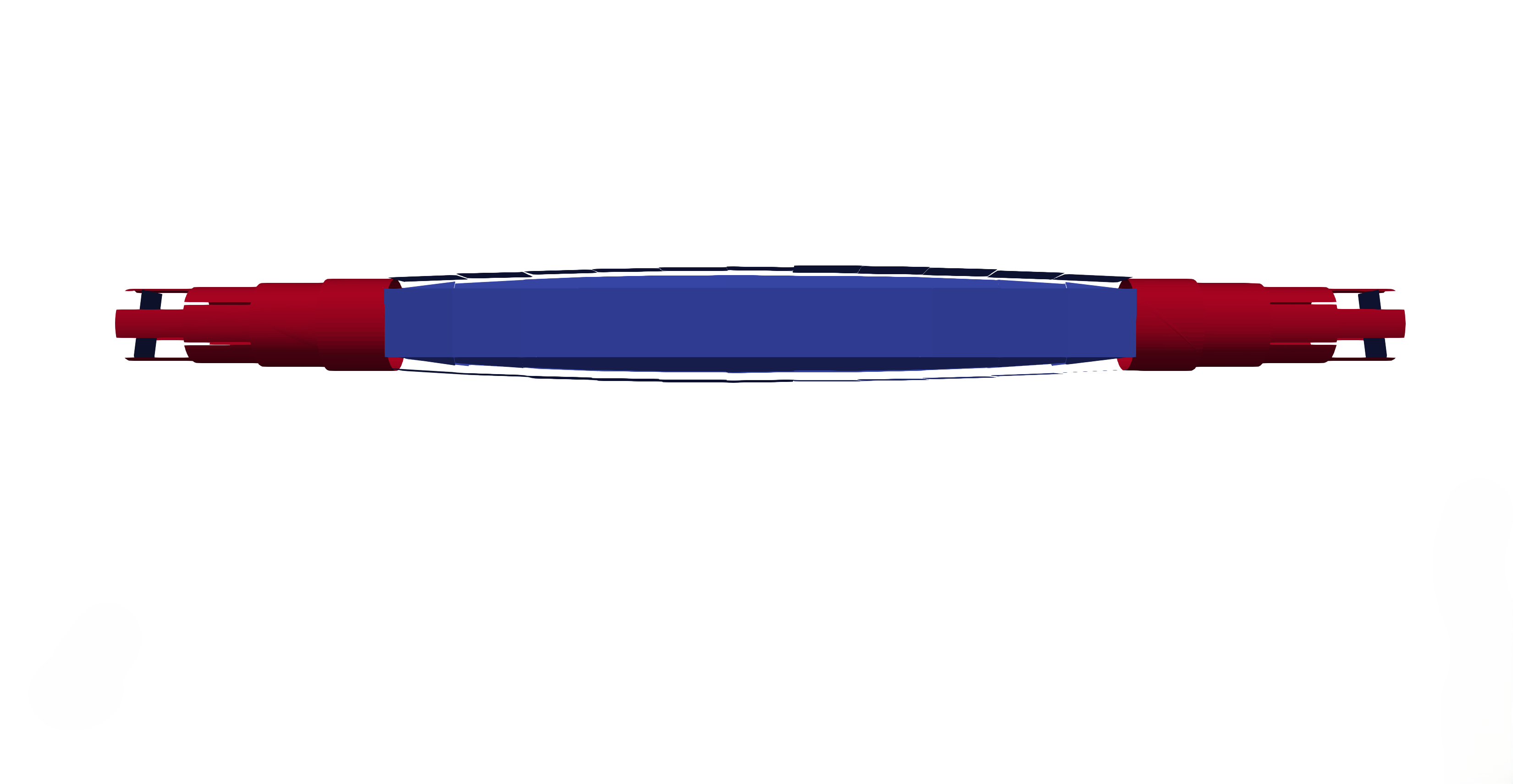}
       \caption{$t = 0.8$}
       \label{stag:sub1}
   \end{subfigure}%
   \\
   \begin{subfigure}{.22\textwidth}
       \centering
       \includegraphics[width=\linewidth, clip]{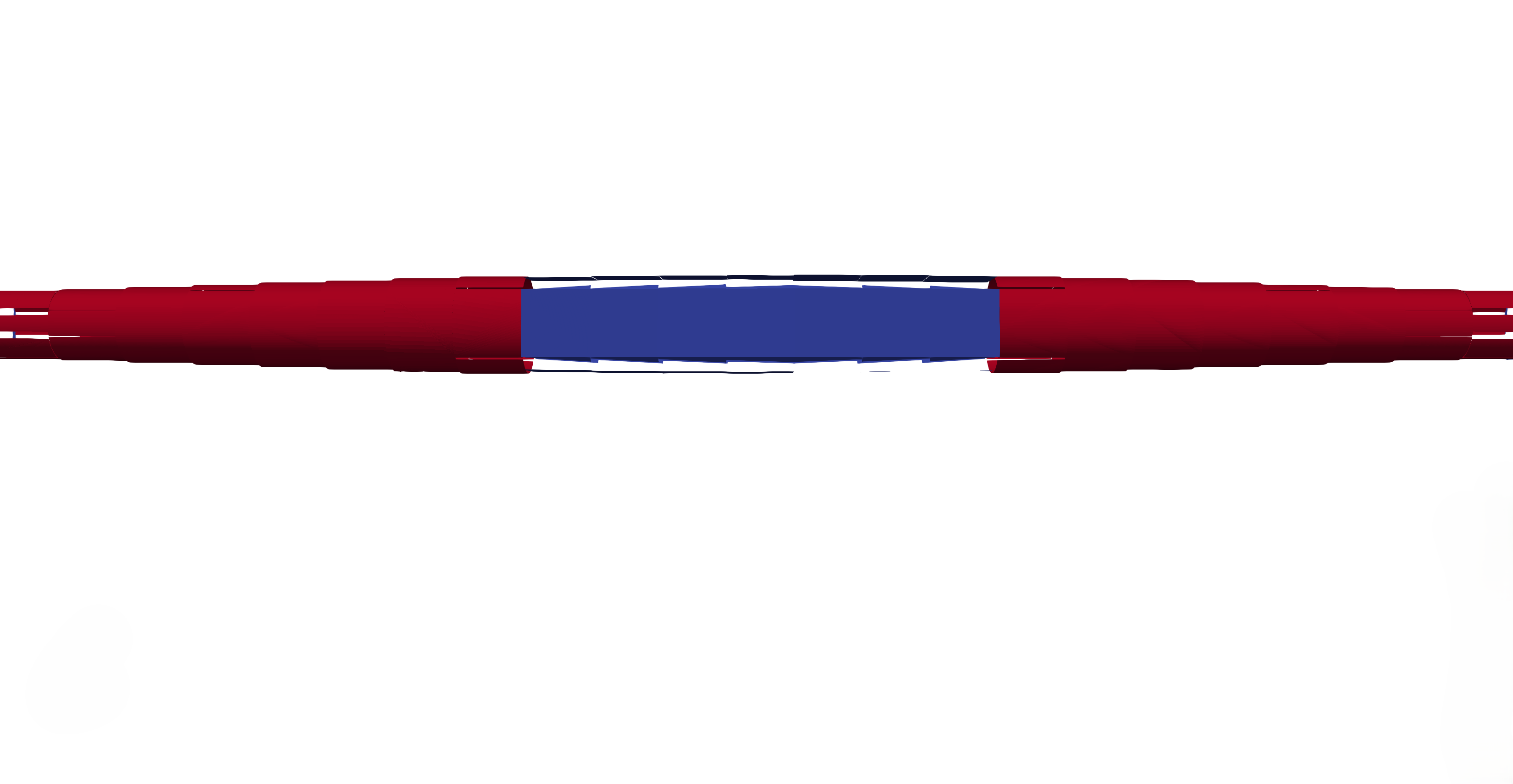}
       \caption{$t = 0.94$}
       \label{stag:sub2}
   \end{subfigure}%
   \hspace{0.4em}
   \begin{subfigure}{.22\textwidth}
       \centering
       \includegraphics[width=\linewidth, clip]{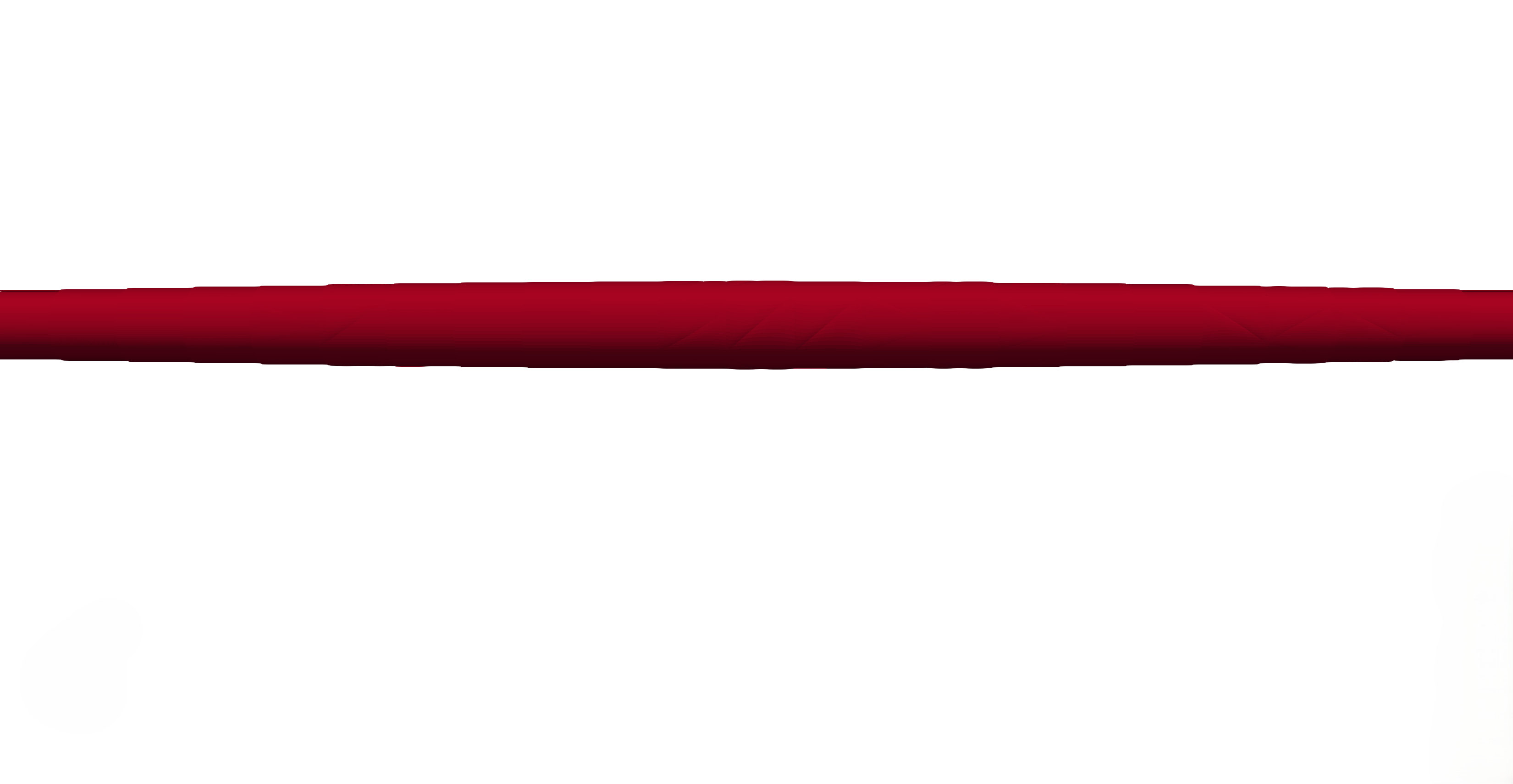}
       \caption{$t = 1.1$}
       \label{stag:sub3}
   \end{subfigure}%
   \hspace{0.4em}
   \begin{subfigure}{.22\textwidth}
       \centering
       \includegraphics[width=\linewidth, clip]{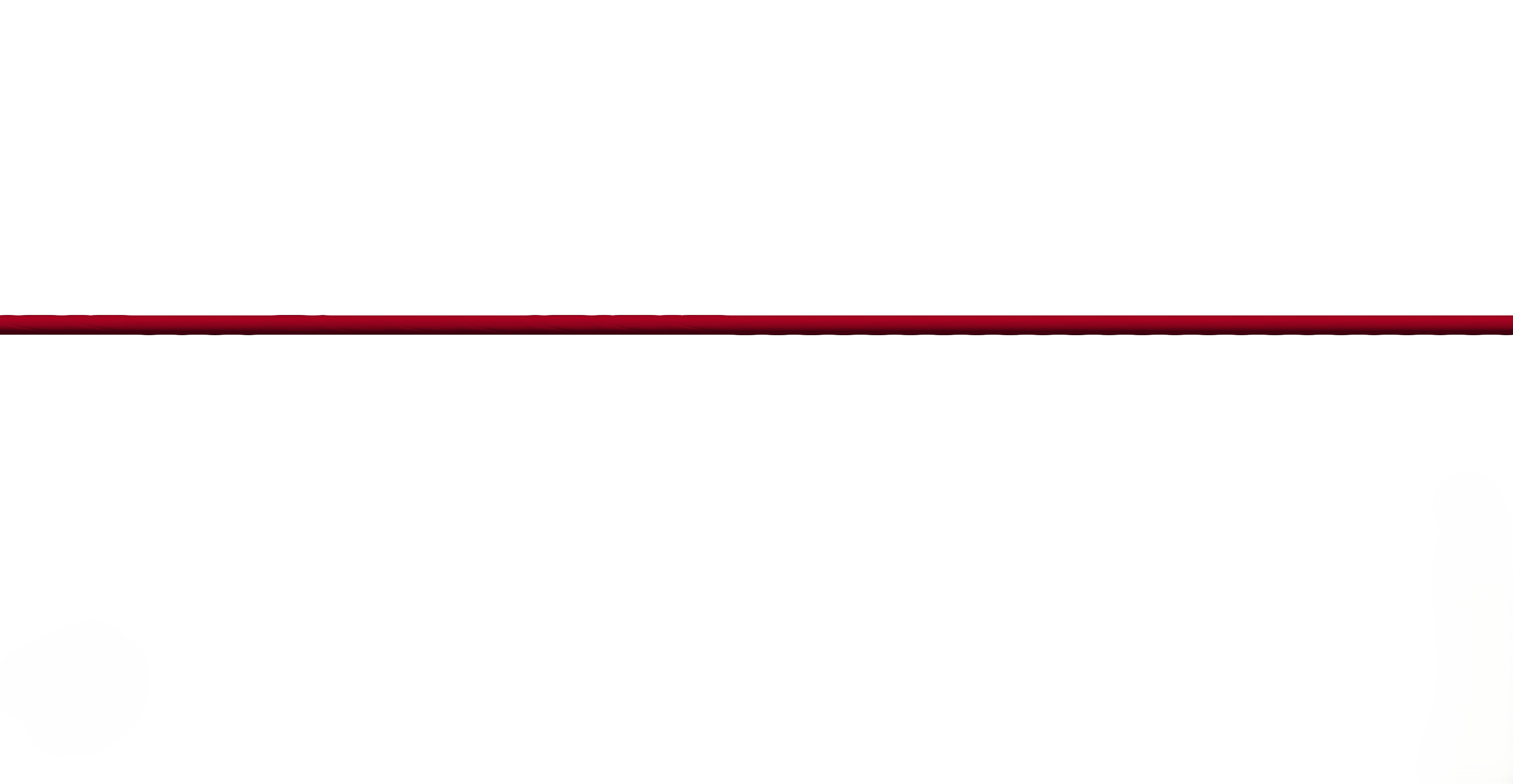}
       \caption{$t = 2.6$}
       \label{stag:sub4}
   \end{subfigure}%
  \caption{Time series of the stagnation flow case. Blue interfaces are PLIC reconstructions, while red interfaces are PCIC reconstructions. All times are non-dimensional.}
  \label{stagnation_series}
\end{figure}

\begin{figure}[t!]
   \centering\vspace{-3em}
   \includegraphics[width=0.3\linewidth, clip]{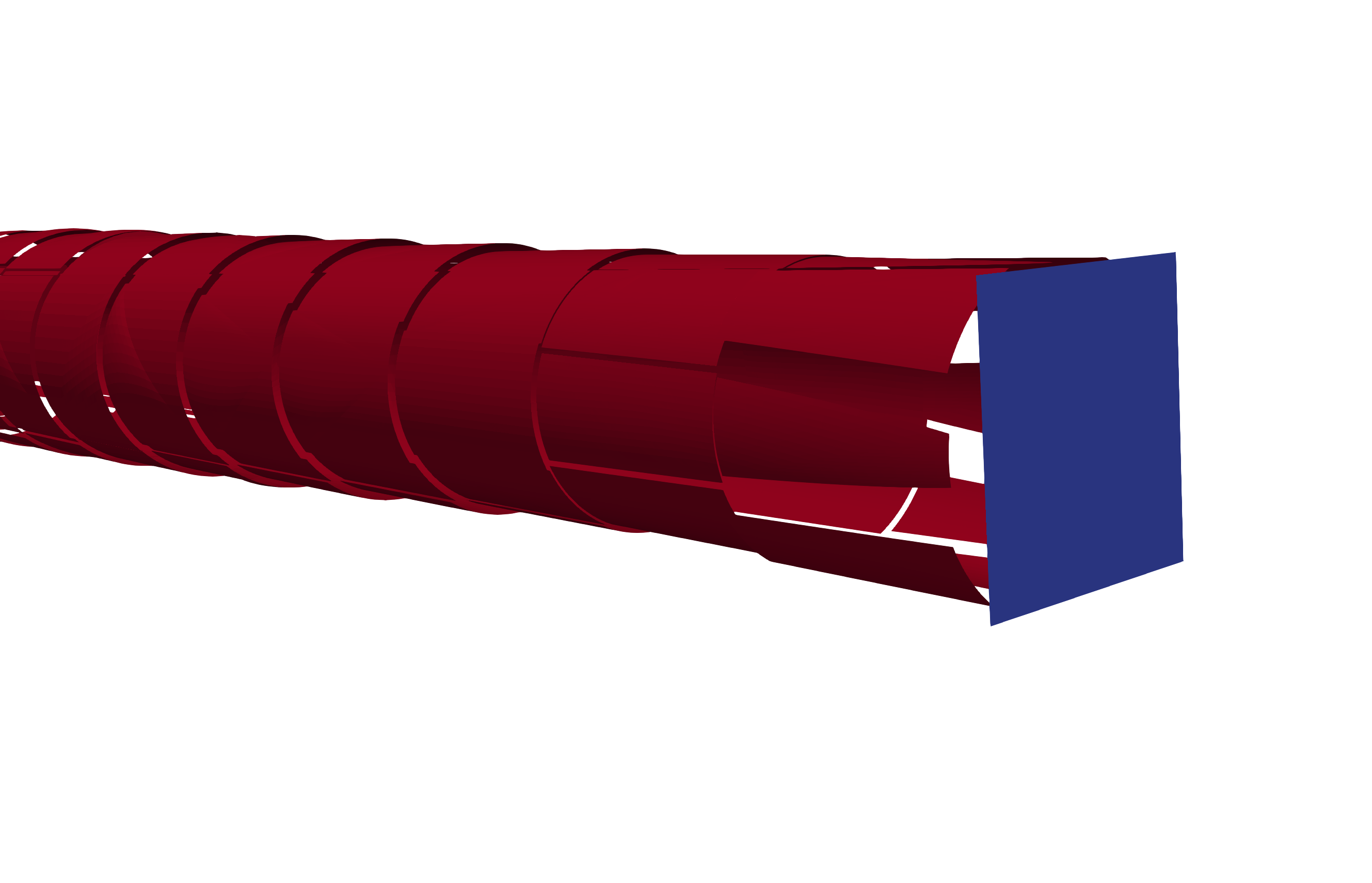}
  \caption{Example of the ligament tip detection. PCIC interfaces are red, and the PLIC cap is highlighted in blue.}
  \label{ligament_cap}
\end{figure}

This flow is designed to stretch the initial spherical drop into a long, thin ligament. The simulation was run in the NGA2 flow solver~\cite{Desjardins}, although a prescribed velocity field was still used. The ligament detection and ligament tip detection methods described in Sections~\ref{Methods:lig_det} and~\ref{Methods:tip_det} were enabled. PLIC-Net was used for the PLIC reconstructions. Figure~\ref{stagnation_series} shows a time series of the case. As the drop is stretched and begins to approach the sub-grid scale, the ligament detection triggers towards the tips of the growing ligament, which are the thinnest regions. These regions are converted to PCIC, while the tip detection maintains a PLIC cap at the actual end of the ligament, as shown in Fig.~\ref{ligament_cap}. As the ligament continues to become thinner, it is progressively converted completely to a PCIC representation, which eventually becomes entirely sub-grid.

These results demonstrate that the ligament detection and ligament tip detection work as intended. Under-resolved ligament-like regions are reliably detected, and the transition from PLIC to PCIC works well and does not result in any spurious interfaces or an artificial breakup. The ligament tip detection is also reliable and successfully prevents any spuriously thin cylinders from forming and propagating at the ends of the structure. Critically, volume is conserved throughout the entire simulation, as it is in every simulation shown in this work.

\subsection{Drop in HIT Flow}
\label{results: hit}

\begin{figure}[b!]
   \centering
   \begin{subfigure}{.31\textwidth}
       \centering
       \includegraphics[width=\linewidth, clip]{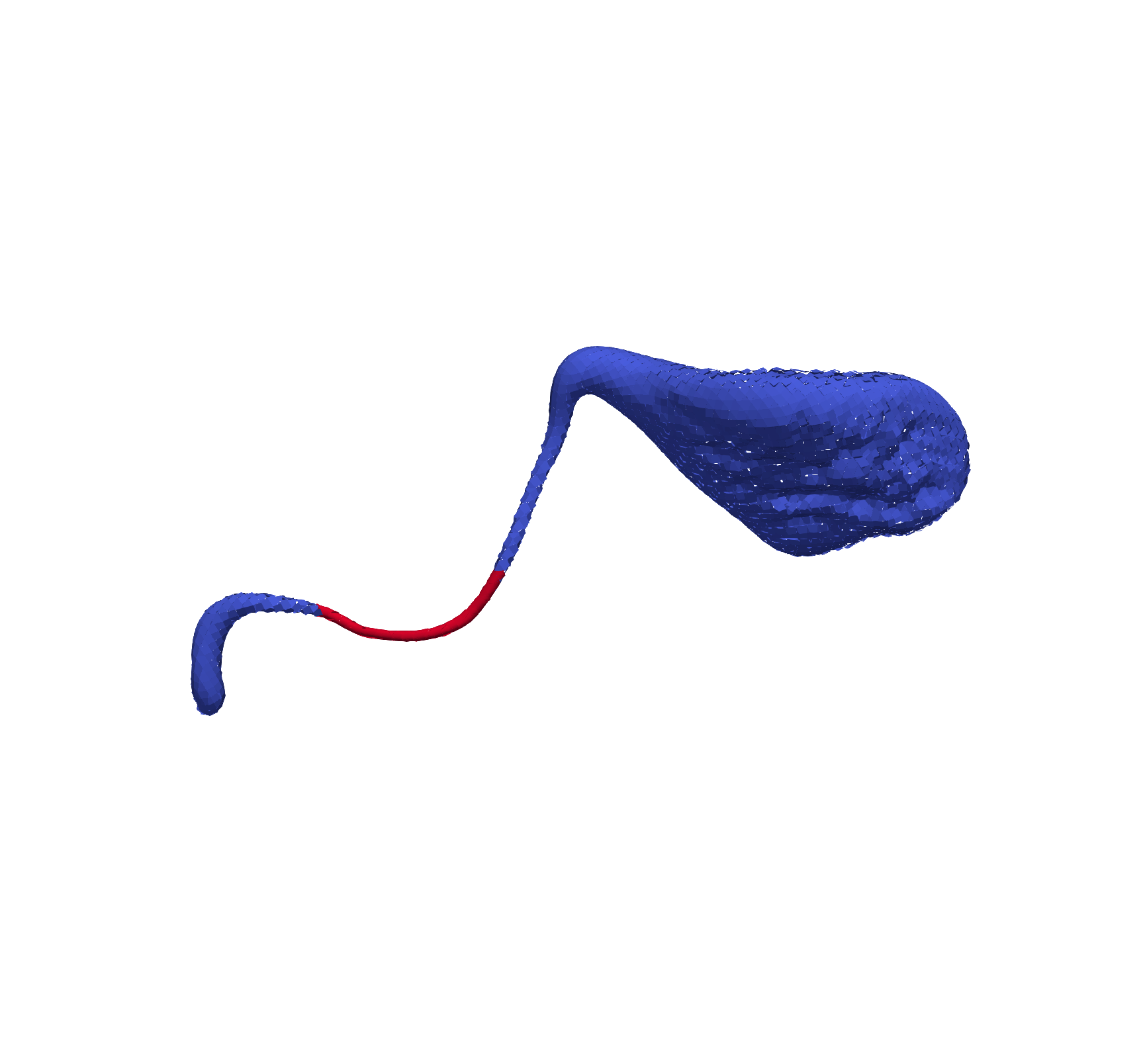}
       \caption{$t = 3.5$}
       \label{hit:sub3}
   \end{subfigure}%
   \hspace{0.4em}
   \begin{subfigure}{.31\textwidth}
       \centering
       \includegraphics[width=\linewidth, clip]{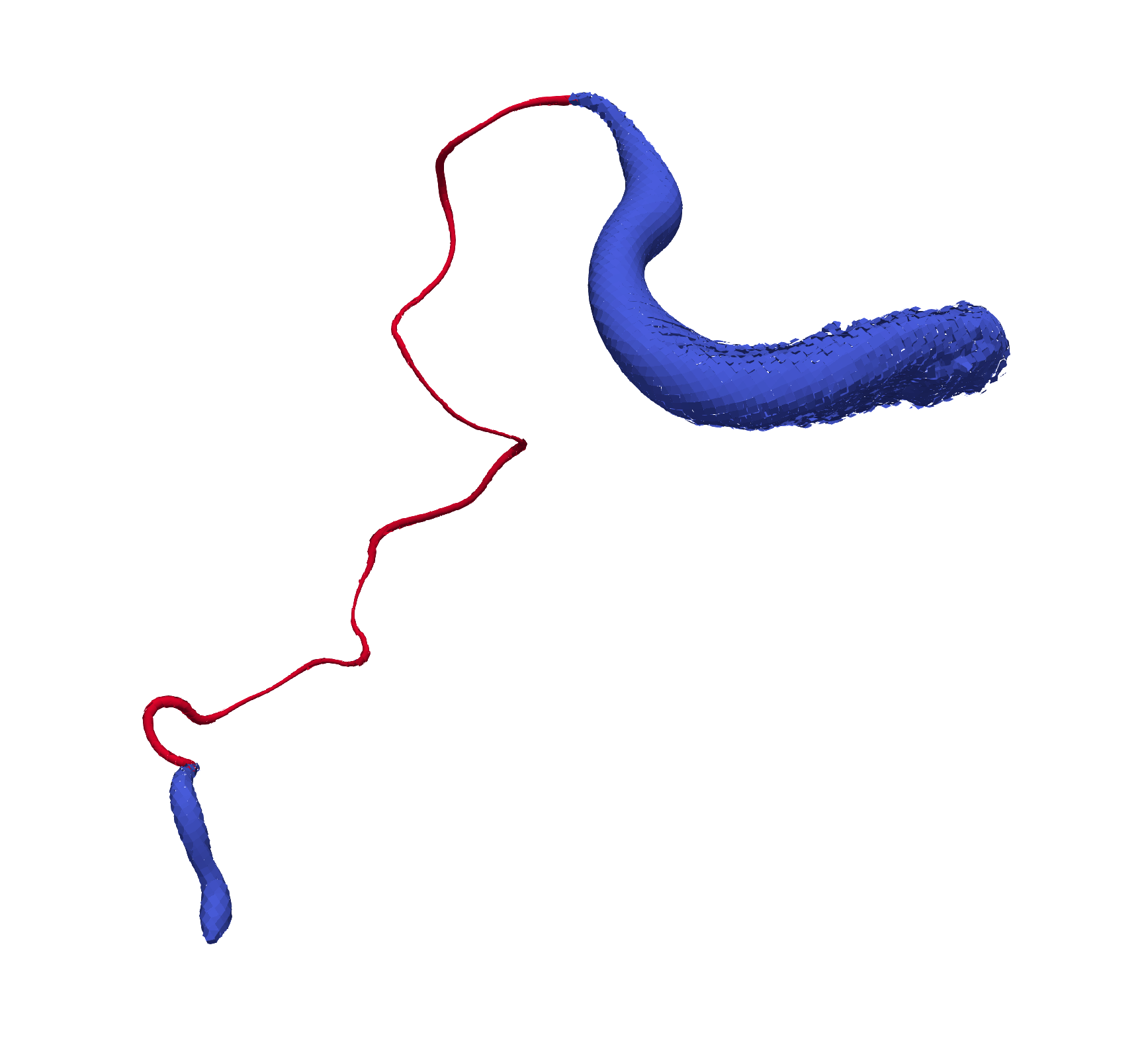}
       \caption{$t = 4$}
       \label{hit:sub4}
   \end{subfigure}%
   \hspace{0.4em}
   \begin{subfigure}{.31\textwidth}
       \centering
       \includegraphics[width=\linewidth, clip]{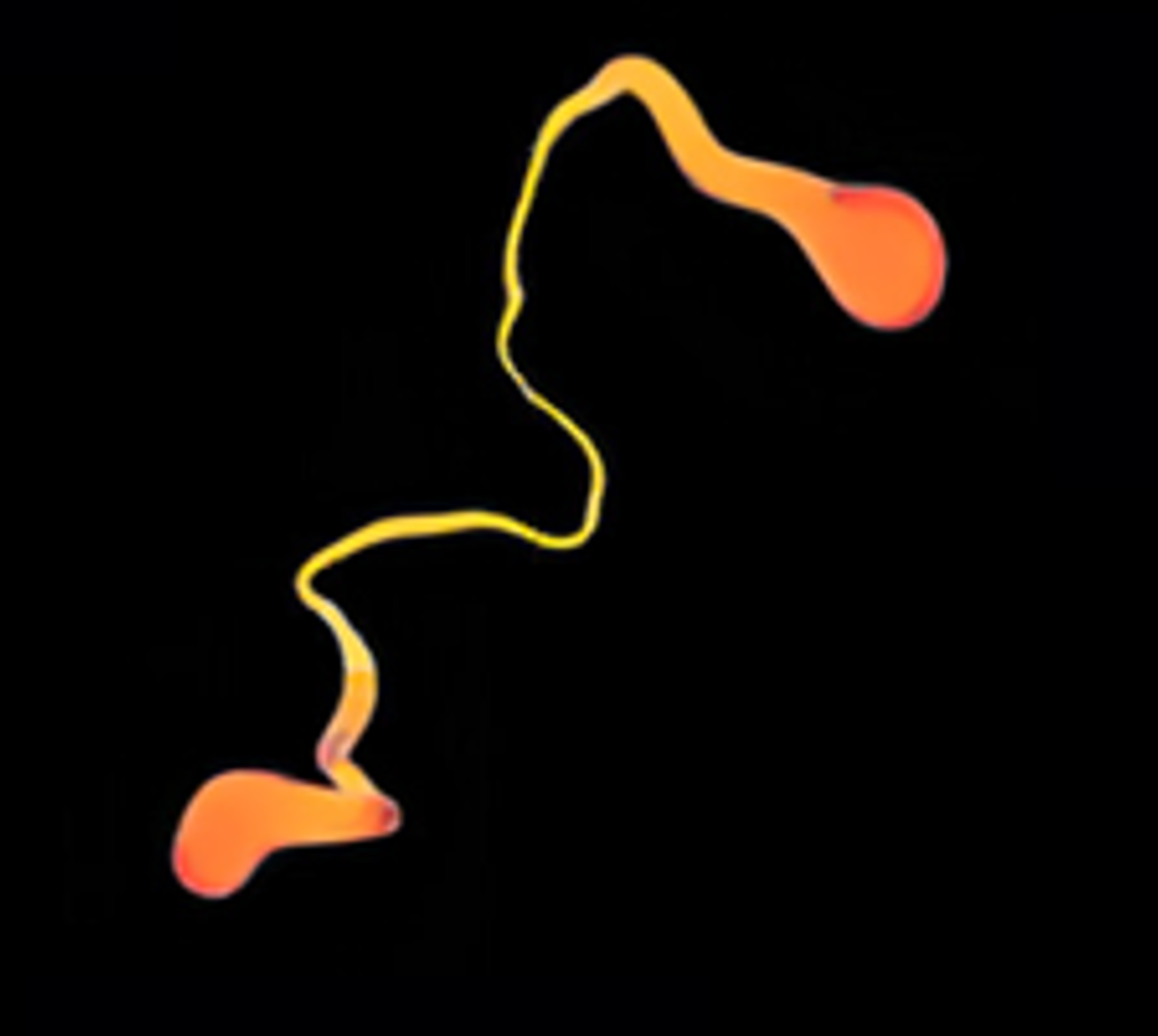}
       \caption{Reference simulation \citep{Farsoiya}}
       \label{hit:sub5}
   \end{subfigure}%
  \caption{(a-b) Time series of the viscous drop in HIT case with $We=20$ and $\mu_r=150$. Blue interfaces are PLIC reconstructions, while red interfaces are PCIC reconstructions. (c) A comparison to the reference simulation from \citet{Farsoiya} is shown. All times are non-dimensionalized with the eddy turnover time, and $t=0$ corresponds to drop injection.}
  \label{hit}
\end{figure}

\begin{figure}
   \centering
       \includegraphics[width=0.4\linewidth, clip]{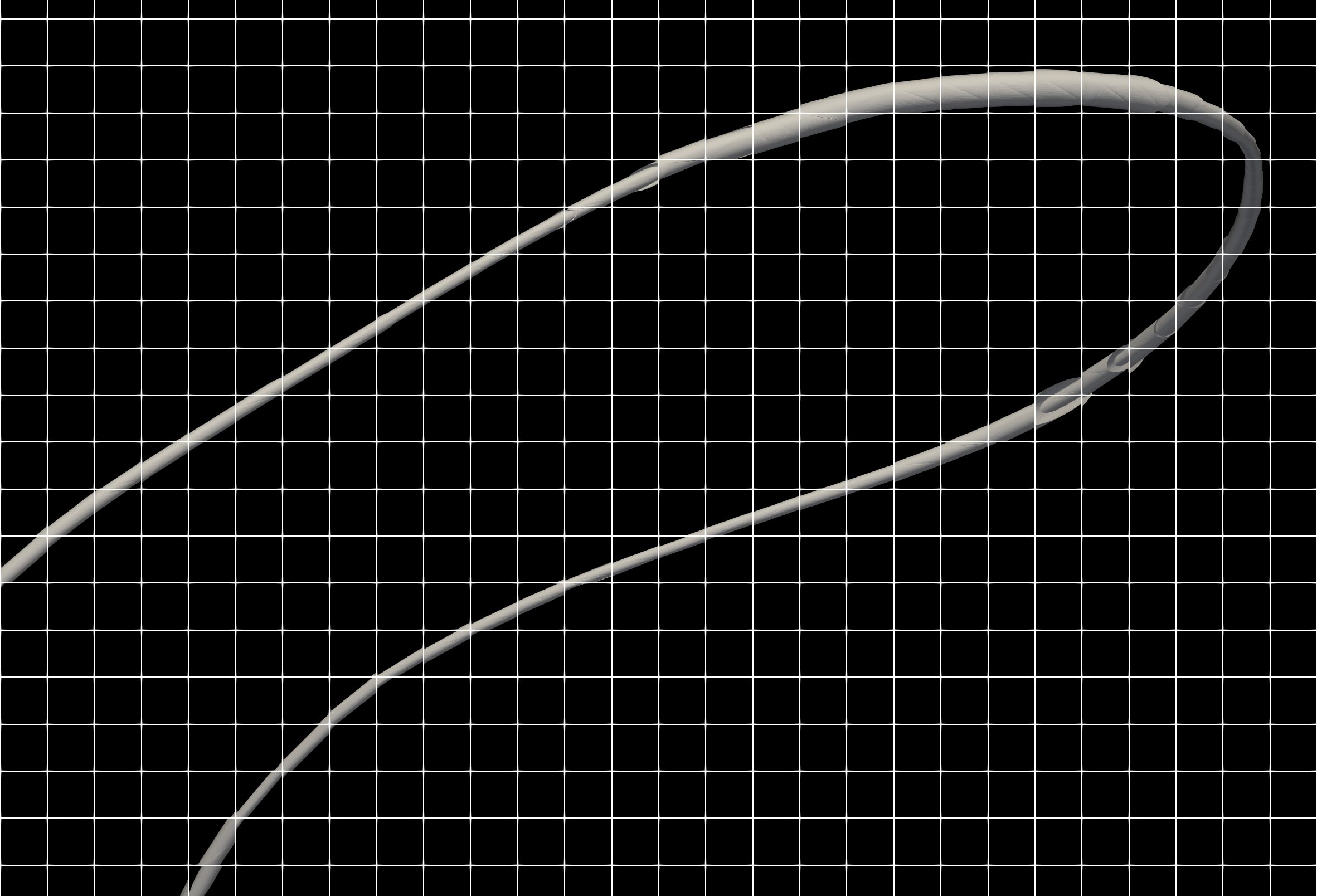}
  \caption{Example of PCIC reconstructing a sub-grid ligament at $t=4.25$ in the viscous drop in a HIT flow case with $We=20$ and $\mu_r=150$. The grid is shown.}
  \label{sub_grid}
\end{figure}

\begin{figure}[t]
   \centering
   \begin{subfigure}{.46\textwidth}
       \centering
       \includegraphics[width=\linewidth, clip]{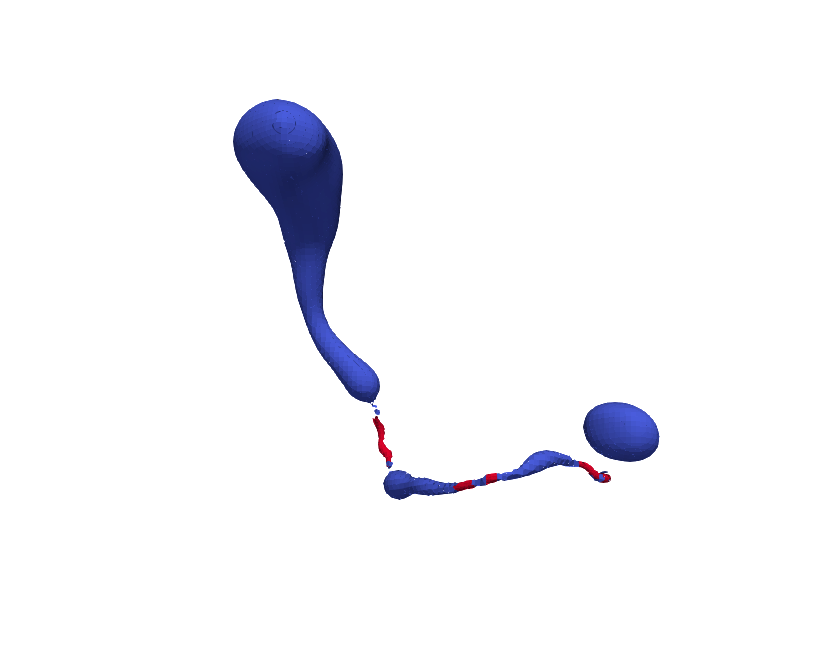}
       \caption{$t = 3.3$}
       \label{no_rp:sub2}
   \end{subfigure}%
   \hspace{0.4em}
   \begin{subfigure}{.46\textwidth}
       \centering
       \includegraphics[width=\linewidth, clip]{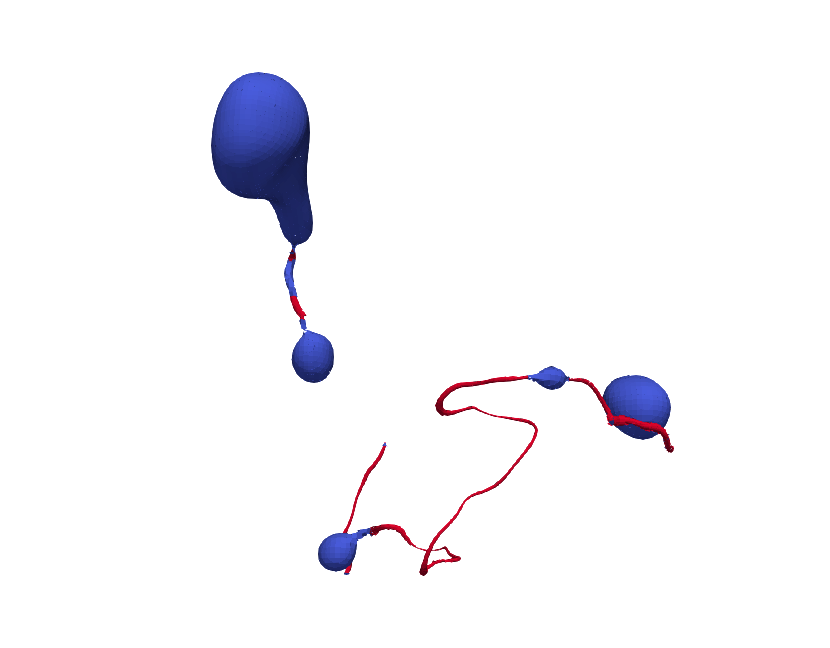}
       \caption{$t = 3.6$}
       \label{no_rp:sub3}
   \end{subfigure}%
  \caption{Time series of the viscous drop in HIT case with $We=4$ and $\mu_r=10$ and no Rayleigh--Plateau modeling. Blue interfaces are PLIC reconstructions, while red interfaces are PCIC reconstructions. All times are non-dimensionalized with the eddy turnover time, and $t=0$ corresponds to drop injection.}
  \label{no_rp}
\end{figure}

\begin{figure}[!ht]
   \centering
   \begin{subfigure}{.31\textwidth}
       \centering
       \includegraphics[width=\linewidth, clip]{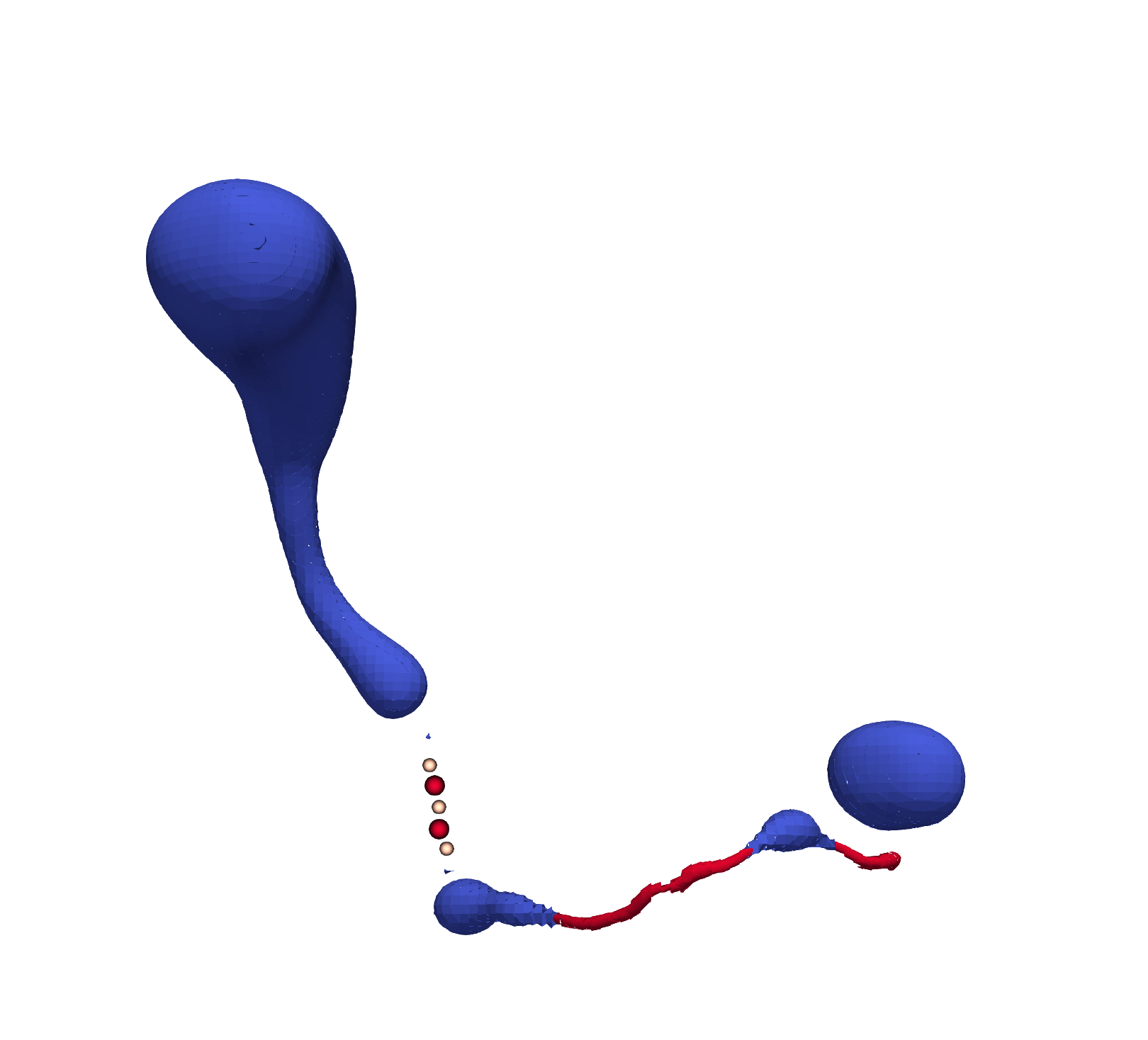}
       \caption{$t = 3.35$}
       \label{rp:sub2}
   \end{subfigure}%
   \hspace{0.4em}
   \begin{subfigure}{.31\textwidth}
       \centering
       \includegraphics[width=\linewidth, clip]{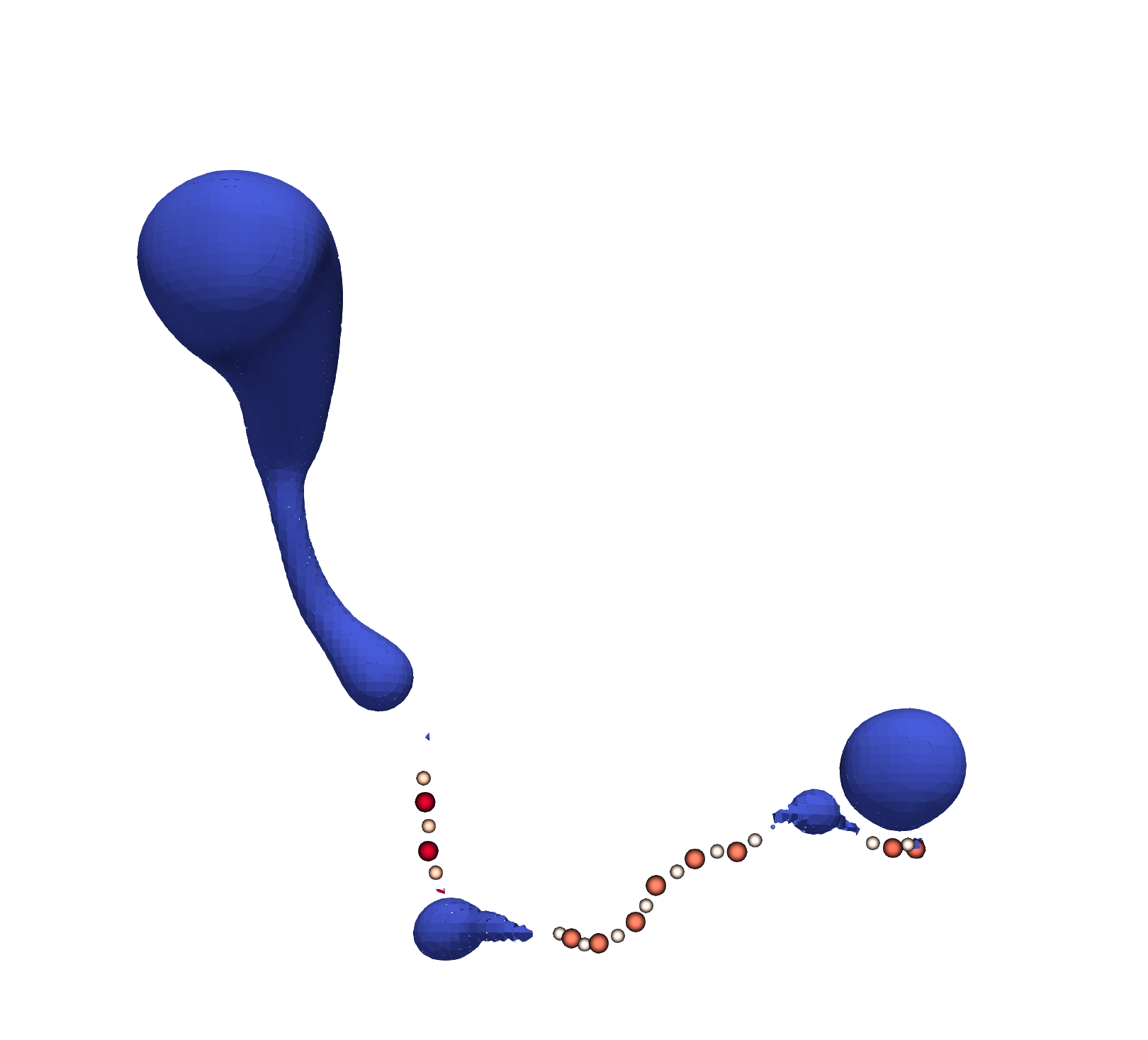}
       \caption{$t = 3.4$}
       \label{rp:sub3}
   \end{subfigure}%
   \hspace{0.4em}
   \begin{subfigure}{.31\textwidth}
       \centering
       \includegraphics[width=\linewidth, clip]{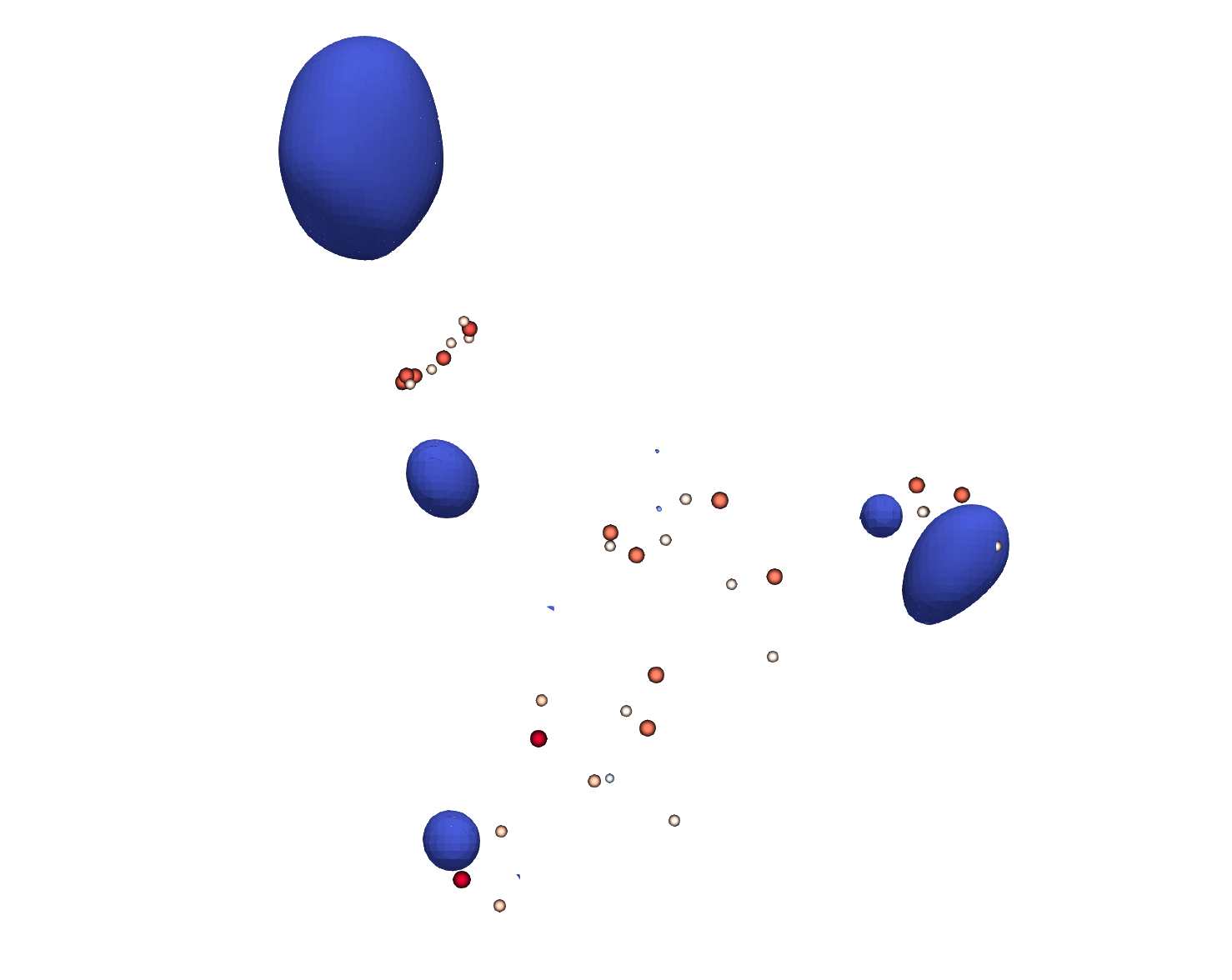}
       \caption{$t = 3.8$}
       \label{rp:sub3_1}
   \end{subfigure}%
   \\
   \vspace{1.8em}
   \begin{subfigure}[t]{.33\textwidth}
       \centering
       \includegraphics[width=\linewidth, clip]{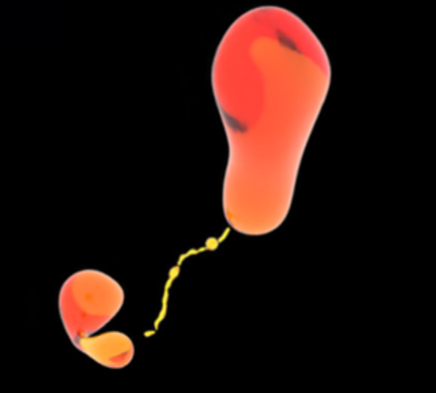}
       \caption{Reference simulation undergoing Rayleigh--Plateau breakup \citep{Farsoiya}.}
       \label{rp:sub4}
   \end{subfigure}%
   \hspace{0.4em}
   \begin{subfigure}[t]{.33\textwidth}
       \centering
       \includegraphics[width=\linewidth, clip]{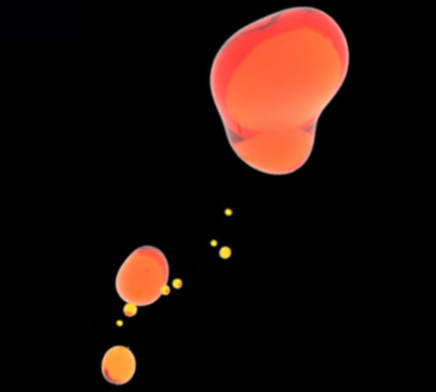}
       \caption{Reference simulation after Rayleigh--Plateau breakup \citep{Farsoiya}.}
       \label{rp:sub5}
   \end{subfigure}%
  \caption{(a-c) Time series of the viscous drop in HIT case with $We=4$ and $\mu_r=10$ and Rayleigh--Plateau modeling. Blue interfaces are PLIC reconstructions, while red interfaces are PCIC reconstructions. Lagrangian drops are colored by radius. Note that the displayed size ratio between the mother and satellite Lagrangian drops deviates from the inviscid Rayleigh-Plateau theory for visualization purposes. All times are non-dimensionalized with the eddy turnover time, and $t=0$ corresponds to drop injection. (c-d) A comparison with the reference simulation \citep{Farsoiya} is shown.}
  \label{rp}
\end{figure}

PCIC was tested in a complex multiphase flow simulation involving a viscous drop exposed to a homogeneous isotropic turbulent (HIT) flow. This case was studied by \citet{Farsoiya}, in which the researchers used a highly resolved AMR simulation to investigate the effects of Weber number, Reynolds number, and liquid/gas viscosity ratio on the breakup behavior of the drop. The simulation is constructed as follows. First, the HIT flow is developed in a cubic domain with only the gas phase present. Once the turbulence statistics have converged, the viscous drop is injected into the center of the domain; the velocity inside the drop is initially set to $0$. The HIT flow then perturbs the drop and causes it to deform and stretch, creating thin ligaments and/or triggering breakup events. In the work by \citet{Farsoiya}, a case with a Taylor microscale Reynolds number ($Re_{\lambda}$) of $77$, a Weber number ($We$) of $20$, and a liquid/gas viscosity ratio ($\mu_r$) of 150 resulted in a thin, meandering ligament. Thus, this reference case was chosen to test PCIC. The Taylor microscale Reynolds number is defined as $Re_{\lambda}=u'\lambda_t/\nu_g$, where $u'$ is the root mean square of the velocity fluctuations, $\lambda_t$ is the Taylor microscale, and $\nu_g$ is the gas kinematic viscosity. The Weber number is defined with a characteristic length of the drop's initial diameter. The density ratio is set to unity. The reference AMR simulation used an effective grid size of up to $1024^3$, but in this work the simulation was run without AMR at a resolution of $192^3$. At this resolution, the maximum allowable $Re_{\lambda}$ that still resolves the turbulent scales can be calculated as $Re_{\lambda}=72.2$~\cite{Pope}. This is slightly lower than the $Re_{\lambda}$ from the reference case, but $We$ and $\mu_r$ were both set to match the reference. The domain had a non-dimensional size of $5^3$, and the drop diameter was initially $0.67$ at the time of injection. A maximum CFL number of 0.9 was used, which included the capillary time step constraint. PLIC-Net was used for PLIC reconstructions, and PCIC was used for sub-grid cylinder reconstructions with its ligament detection and ligament tip detection capabilities.

A time series of the simulation is shown in Fig.~\ref{hit} along with a comparison to the reference case. Despite using a much coarser mesh, the results show a striking resemblance to the reference; a long, meandering ligament forms between two resolved droplets. This is possible because of PCIC's ability to maintain the ligament well below the mesh size, as demonstrated in Fig.~\ref{sub_grid}. As observed in the previous simpler test cases, PCIC is a robust reconstruction method without visible errors, despite being exposed to turbulence-induced sharp bending. The ligament detection also works reliably as it efficiently converts the PLIC interfaces to PCIC without noticeable error or interface breakup. By producing results comparable to a $1024^3$ AMR simulation with only $192^3$ cells, PCIC allows for a significant reduction in computational resources without sacrificing simulation accuracy. It must be noted that, since this work's objective was to demonstrate PCIC's reconstruction abilities, no surface tension modeling was done on the cylindrical interfaces. Future work will address such sub-grid surface tension modeling to enable further applications of PCIC.

This simulation showed how PCIC can be used in cases involving very thin ligaments that are maintained for a long time below the mesh size. However, PCIC can also be used when a faster breakup of ligaments occurs, such as when the Weber number and viscosity ratio are instead set to $We=4$ and $\mu_r=10$. If PCIC is naively applied to this case, then unphysical ligaments would form regardless, as shown in Fig.~\ref{no_rp}. In this case, separation between the resolved drops and ligaments is observed due to sharper interface curvatures between the PLIC and PCIC regions. The presence of the unphysical ligaments can be addressed by introducing a Rayleigh--Plateau breakup model derived from the work of \citet{Kim}. To demonstrate this capability, a simple version of this model is tested on the HIT breakup case with $We=4$ and $\mu_r=10$. For this test, the inviscid model of \citet{Kim} is implemented to convert sub-grid ligament structures into Lagrangian droplets. A CCL approach derived from the one described in Section~\ref{Methods:lig_det} is again used to detect ligaments, and the breakup is triggered when a CCL structure has a minimum diameter of at most $\Delta x$. The resulting Lagrangian droplets are evenly distributed along a cubic spline~\cite{Dierckx} fitted to the cylinder origins in their CCL structure. The time series of this case is shown in Fig.~\ref{rp} along with a comparison to the reference case from \citet{Farsoiya}. With the addition of the breakup modeling, the simulation now produces a physical result since the PCIC ligaments are effectively converted into Lagrangian droplets before unphysical thin structures form.

The above results successfully demonstrate PCIC's abilities under two very different conditions, which provides confidence that the sub-grid reconstruction aspect of the ligament modeling challenge is addressed well by PCIC. Therefore, future work will continue to develop physics-based breakup models for inviscid, viscous, and non-Newtonian cases so that PCIC can be used in a variety of cases of interest to produce accurate flow statistics.

\section{Conclusion}
\label{Conclusion}
In this work, a novel cylindrical reconstruction method named PCIC was introduced for sub-grid ligament modeling in VOF simulations. The method was enabled by the ability to compute the exact volume moments of polyhedra clipped by quadratic cylinders. These moments were calculated from closed-form expressions, which were derived with successive applications of the divergence theorem to simplify the three-dimensional integrals into one-dimensional integrals. Then, a reconstruction based on these volume moments was presented. Data from a $5\times5\times5$ stencil of cells was used to fit a straight circular cylinder in the stencil's center cell. A quadratic principal curve was fitted to the liquid barycenters in the stencil, which was used to provide the orientation and origin of the cylinder. The cylinder's radius was then set to conserve volume through iterative solution of a root-finding problem. Ligament detection was carried out with a CCL algorithm designed to detect thin ligaments before they become thinner than the grid size. Finally, the tips of the ligaments were detected so that they could be capped with a PLIC reconstruction. 

PCIC was thoroughly tested in pure advection and flow solver test cases. In all cases, the reconstruction was clean and robust, and was capable of handling ligaments with high curvatures without breaking or shedding spurious interfaces. The ligament detection algorithm also worked reliably and managed a clean transition from PLIC to PCIC regions, while the tip detection algorithm prevented spuriously thin ligaments from growing. When deployed in a complex case involving a viscous drop in homogeneous isotropic turbulence, PCIC enabled results comparable to highly resolved AMR simulations at a fraction of the resolution. Thus, this work successfully addressed the interface reconstruction problem for sub-grid fluid ligaments. Future work will develop sub-grid models for surface tension and breakup models for a wide range of scenarios, including viscous and non-Newtonian cases.

\section*{Acknowledgments}
This work was supported by the National Science Foundation, Division of Chemical, Bioengineering, Environmental and Transport Systems (CBET) under Awards \#2438853 and \#2438854. A.~C.~has received support from the International Fine Particle Research Institute. V.~W.~has received support from the University Research School EUR-MINT (French government support managed by the National Research Agency for Future Investments program bearing the reference ANR-18-EURE-0023).

\appendix
\section{Third Contribution to the Volume Moments} 
\begingroup
\allowdisplaybreaks
\label{apdx:M3}
The vectorial operator $\boldsymbol{\mathcal{B}}^{(3)}$, introduced in Eq.~\eqref{eq:m3k}, can be expressed as
\begin{equation}
    \boldsymbol{\mathcal{B}}^{(3)} \left(w, \mathbf{x}_a, \mathbf{x}_b, \mathbf{x}_c\right) = \diag\left(\boldsymbol{\mathcal{E}}(w)\right)\left( \boldsymbol{\mathcal{C}} (\mathbf{x}_a, \mathbf{x}_b, \mathbf{x}_c) \boldsymbol{\mathcal{K}} \boldsymbol{\mathcal{D}}(w) \right) \, ,
\end{equation}
where 
$\boldsymbol{\mathcal{C}}$ is a $4 \times 6$ matrix 
whose non-zero coefficients contributing to the volume of the clipped polyhedron are given as
\begin{align}
    {\mathcal{C}}_{1,1}  & = x_a + x_b \, ,\nonumber \\
    {\mathcal{C}}_{1,2}  & = x_c \, .
\end{align}
The non-zero coefficients of $\boldsymbol{\mathcal{C}}$ contributing to the $x$-component of the first moment of the clipped polyhedron are given as
\begin{align}
    {\mathcal{C}}_{2,3}  & = (x_a + x_b)^2 \, ,\nonumber \\
    {\mathcal{C}}_{2,4}  & = (x^2_a + x^2_b) \, ,\nonumber \\
    {\mathcal{C}}_{2,5}  & = (x_a + x_b) \, x_c \, ,\nonumber \\
    {\mathcal{C}}_{2,6}  & = x^2_c \, ,
\end{align}
the non-zero coefficients of $\boldsymbol{\mathcal{C}}$ contributing to the $y$-component of the first moment of the clipped polyhedron are given as
\begin{align}
    {\mathcal{C}}_{3,3}  & = 2(x_a + x_b) (y_a + y_b) \, ,\nonumber \\
    {\mathcal{C}}_{3,4}  & = 2(x_a \, y_a + x_b \, y_b) \, ,\nonumber \\
    {\mathcal{C}}_{3,5}  & = (x_a + x_b) \, y_c +  (y_a + y_b) \, x_c\, ,\nonumber \\
    {\mathcal{C}}_{3,6}  & = 2 \, x_c \, y_c \, ,
\end{align}
and the non-zero coefficients of $\boldsymbol{\mathcal{C}}$ contributing to the $z$-component of the first moment of the clipped polyhedron are given as
\begin{align}
    {\mathcal{C}}_{4,3}  & = 2(x_a + x_b) (z_a + z_b) \, ,\nonumber \\
    {\mathcal{C}}_{4,4}  & = 2(x_a \, z_a + x_b \, z_b) \, ,\nonumber \\
    {\mathcal{C}}_{4,5}  & = (x_a + x_b) \, z_c +  (z_a + z_b) \, x_c\, ,\nonumber \\
    {\mathcal{C}}_{4,6}  & = 2 \, x_c \, z_c \, .
\end{align}
The constant matrix $\boldsymbol{\mathcal{K}}$, of size $6\times6$, is given as
{\renewcommand*{\arraystretch}{1.2}\begin{equation}
   {\text{\normalsize} \boldsymbol{\mathcal{K}}} = {\text{\small}\begin{bmatrix}
        1 & -\frac{5}{6} & 0 & \frac{1}{3} & 0 & 0 \\ 
        0 & \frac{2}{3} & -2 & \frac{1}{3} & 0 & 0 \\ 
        -\frac{3}{16} & \frac{23}{96} & -\frac{1}{8} & -\frac{1}{8} & 0 & \frac{1}{24} \\
        -\frac{1}{8} & \frac{5}{48} & \frac{1}{8} & -\frac{7}{48} & 0 & \frac{1}{24} \\
        0 & -\frac{1}{3} & \frac{5}{4} & -\frac{3}{8} & 0 & \frac{1}{12} \\
        0 & 0 & -\frac{1}{4} & \frac{13}{24} & -1 & \frac{1}{12} \\
    \end{bmatrix}} \, ,
\end{equation}}
$\boldsymbol{\mathcal{D}}$ is the vector given as
\begin{align}
    \boldsymbol{\mathcal{D}}(w) & = \begin{bmatrix}[1]0 & w^2 & 0 & w^4 & 0 & w^6\end{bmatrix}^\intercal + \Theta(w) \begin{bmatrix}[1] w & 0 & w^3 & 0 & w^5 & 0 \end{bmatrix}^\intercal \, 
\end{align}
and $\boldsymbol{\mathcal{E}}$ is the vector given as
\begin{equation}
    \boldsymbol{\mathcal{E}}(w) = \begin{bmatrix}[1] \Lambda(w)^2 & \Lambda(w)^3 & \Lambda(w)^3 & \Lambda(w)^3 \end{bmatrix}^\intercal \, ,
\end{equation}
with
\begin{align}
    \Theta(w) & = \left\{ \begin{array}{ll} \text{arctan}\left(\dfrac{1 - w}{\sqrt{1 - w^2}}\right)\dfrac{1}{\sqrt{1-w^2}} & 0 < w < 1 \\ & \\
            \text{arctanh}\left(\dfrac{w - 1}{\sqrt{w^2 - 1}}\right)\dfrac{1}{\sqrt{w^2-1}} & 1 \le w \end{array} \right.  \, , \\
            \Lambda(w) & = \frac{1}{(w - 1) (w + 1)}  \, .
\end{align}
As in \citet{Evrard}, the Taylor series expansion of $\smash{\boldsymbol{\mathcal{B}}^{(3)}}$ to order $40$ is used for $w \in [0.35,1.7]$ so as to avoid introducing significant round-off errors when $w$ is in the vicinity of $1$. This implementation has been employed for producing the results presented in Sections~\ref{Results1} and \ref{Results2}.
\endgroup

\bibliographystyle{model1-num-names}
\bibliography{bibliography.bib}

\end{document}